\documentclass[11pt]{article}
\pdfoutput=1
\usepackage[margin=1in,marginparwidth=.95in]{geometry}

\usepackage[utf8x]{inputenc}
\usepackage{amsmath,amssymb,amsthm}
\usepackage[pdftex]{graphicx}   
\usepackage{color}
\usepackage{ifthen}
\usepackage{mathtools}
\usepackage{todonotes}
\usepackage{verbatim}

\usepackage{yfonts}

\usepackage{tikz}
\usetikzlibrary{arrows,matrix}
\tikzset{arrow/.style={semithick,>=stealth',shorten >=1pt,shorten <=1pt}}

\usepackage[all]{xy}

\usepackage[pdftex,     
            plainpages=false,   
            breaklinks=true,    
            pdftitle={Transfer for free loop spaces of fusion systems and covering maps},
            pdfauthor={Sune Precht Reeh, Tomer M. Schlank, Nathaniel Stapleton}
           ]{hyperref}
\usepackage{thumbpdf}

\usepackage[abbrev,msc-links]{amsrefs}

\numberwithin{equation}{section}
\numberwithin{figure}{section}
\numberwithin{table}{section}
\allowdisplaybreaks

\theoremstyle{plain}
\newtheorem{theorem}{Theorem}[section]
\newtheorem{prop}[theorem]{Proposition}
\newtheorem{lemma}[theorem]{Lemma}
\newtheorem{cor}[theorem]{Corollary}

\theoremstyle{definition}
\newtheorem{definition}[theorem]{Definition}
\newtheorem{convention}[theorem]{Convention}

\newtheoremstyle{remark}%
{5pt}
{11pt}
{}
{}
{\itshape}
{.}
{ }
{}

\theoremstyle{remark}
\newtheorem{remark}[theorem]{Remark}

\newtheorem{example}[theorem]{Example}

\newcommand{\ph}{\varphi}
\newcommand{\e}{\varepsilon}
\renewcommand{\hat}{\widehat}
\renewcommand{\tilde}{\widetilde}

\DeclareMathOperator{\Hom}{Hom}
\DeclareMathOperator{\Map}{Map}

\DeclareMathOperator{\Ho}{Ho}
\DeclareMathOperator*{\colim}{colim}

\newcommand{\cE}{\mathcal{E}}
\newcommand{\cF}{\mathcal{F}}
\newcommand{\cG}{\mathcal{G}}

\newcommand{\ZZ}{\mathbb{Z}}

\newcommand{\x}{\times}

\newcommand{\into}{\hookrightarrow}

\newcommand{\xto}{\xrightarrow}

\newcommand{\cmp}{\mathbin{\odot}}

\newcommand{\cntuples}[2]{{#2}^{[#1]}}

\newcommand{\cnptuples}[2]{{#2}^{[#1],p}}
\newcommand{\freeO}[1]{\ensuremath \Loops^{#1}}
\newcommand{\tup}{\underline}
\newcommand{\ptup}{\partial\tup}
\newcommand{\change}{\zeta}
\newcommand{\rep}{\PackageError{main}{\\rep not currently implemented}{Use \\reptup instead}\textswab}
\newcommand{\reptup}[1]{\tup{\textswab{#1}}}

\DeclareMathOperator{\wind}{wind}

\newcommand{\F}{\mathbb{F}}

\newcommand{\Z}{\mathbb{Z}}

\newcommand{\Loops}{\mathcal{L}}
\newcommand{\twistO}[1]{L^{\dagger}_{#1}}
\newcommand{\untwistO}[1]{L^{\dagger}_{#1}|_{\mathrm{lift}}}
\newcommand{\twistT}[1]{T^\dagger_{#1}}

\DeclareMathOperator{\Id}{Id}

\DeclareMathOperator{\Res}{Res}

\DeclareMathOperator{\id}{id}
\DeclareMathOperator{\incl}{incl}
\DeclareMathOperator{\iso}{iso}
\DeclareMathOperator{\tel}{tel}
\DeclareMathOperator{\tr}{tr}
\DeclareMathOperator{\ev}{ev}
\DeclareMathOperator{\pev}{\partial ev}

\DeclarePairedDelimiter{\gen}{\langle}{\rangle}
\DeclarePairedDelimiter{\abs}{\lvert}{\rvert}

\newcommand{\AF}{\mathbb{AF}_p}
\newcommand{\AG}{\mathbb{AG}}

\newcommand{\Sp}{\mathord{\mathrm{Sp}}}
\newcommand{\Spp}{\mathord{\mathrm{Sp}}_{p}}

\newcommand{\Cov}{\mathord{\mathrm{Cov}}}

\newcommand{\Sinf}{\Sigma^{\infty}}
\newcommand{\Sinfp}{\Sigma^{\infty}_{+}}
\newcommand{\Sinfpp}{\hat\Sigma^{\infty}_{+}}

\newcommand{\lc}[1]{\prescript{#1}{}\!}

\title{Evaluation maps and transfers for free loop spaces II}
\author{Sune Precht Reeh \and Tomer M. Schlank \and Nathaniel Stapleton}
\date{}

\begin{document}

\maketitle

\begin{abstract}
In \cite{RSS_Bold1}, we constructed and studied a functorial extension of the evaluation map $S^1 \times \Loops X \to X$ to transfers along finite covers. In this paper, we show that this induces a natural evaluation map on the full subcategory of the homotopy category of spectra consisting of $p$-completed classifying spectra of finite groups. To do this, we leverage the close relationship between this full subcategory and the Burnside category of fusion systems.
\end{abstract}

\tableofcontents

\section{Introduction}
In \cite{RSS_Bold1}, we produced a functorial extension of the evaluation map $S^1 \times \Loops X \to X$ to transfers along finite covers and showed that this induces a natural evaluation map on the Burnside category of finite groups. The homotopy category of classifying spectra of finite groups admits an algebraic interpretation in terms of completions of Burnside modules of finite groups, but in \cite{RSS_Bold1} we showed that the natural evaluation map does not extend to the homotopy category of classifying spectra of finite groups. For a prime $p$, the homotopy category of $p$-completed classifying spectra of finite groups admits an algebraic interpretation in terms of Burnside modules for fusion systems. Our goal in this paper is to show that the natural evaluation map does extend to the Burnside category of fusion systems and also to the homotopy category of $p$-completed classifying spectra of finite groups. 

The notion of a saturated fusion system is an axiomatization of the properties that can be detected about a finite group $G$ from a choice of Sylow $p$-subgroup $S \leq G$ equipped with conjugation data. There are saturated fusion systems that do not come from a choice of finite group and Sylow $p$-subgroup. Every saturated fusion system has a classifying space that is a natural generalization of the $p$-completion of the classifying space of a finite group. 

The classifying spectrum of a saturated fusion system is built from the classifying space of the saturated fusion system in a simple manner. There is a close relationship between the classifying spectrum of a saturated fusion system and the $p$-completed classifying spectrum of the underlying Sylow $p$-subgroup. There is an idempotent endomorphism of the classifying spectrum of the underlying Sylow $p$-subgroup that splits off the classifying spectrum of the saturated fusion system. This idempotent is known as the characteristic idempotent of the saturated fusion system. If the saturated fusion system comes from a finite group $G$, this construction recovers the $p$-completion of the classifying spectrum of $G$. 

The set of homotopy classes of maps between two $p$-completed classifying spectra of finite $p$-groups is isomorphic to the $p$-completion of the Burnside module of bisets between the two $p$-groups. The $p$-complete Burnside category of $p$-groups (in which the Burnside modules have been completed at $p$) is therefore equivalent to the homotopy category of $p$-completed classifying spectra of finite $p$-groups. The characteristic idempotent for a saturated fusion system appears as an idempotent in the $p$-completion of the double Burnside ring of the underlying Sylow $p$-subgroup. Similarly, there is a $p$-complete Burnside category for all saturated fusion systems, and this is equivalent to the homotopy category of all classifying spectra of fusion systems, which includes $p$-completed classifying spectra of all finite groups.

We extend the natural evaluation map of \cite{RSS_Bold1} from the $p$-complete Burnside category of $p$-groups to the $p$-complete Burnside category of saturated fusion systems by analyzing its effect on the characteristic idempotent. For formal reasons, this provides us with an evaluation map whose domain is a direct summand of the $p$-completed classifying spectrum of $B\Z/p^k \times \Loops BS$. However, work needs to be done to relate this summand to the free loop space of the classifying space of the saturated fusion system.

\subsection*{In more detail...}
Let $\AG$ be the Burnside category of finite groups. The objects of $\AG$ are finite groups and the abelian group of morphisms between two groups $G$ and $H$ is the Burnside module of $H$-free $(G,H)$-bisets. We may formally extend this category to include formal coproducts of finite groups. In this case, a map is given by a matrix of virtual bisets. The homotopy category of classifying spaces of finite groups faithfully embeds in $\AG$. For a finite group $G$, let $\Loops G = \coprod_{[g]} C_G(g)$, the formal coproduct over conjugacy classes of elements in $G$ of the centralizers. In \cite{RSS_Bold1}*{Section 3}, we constructed and studied a functor $\twistO n \colon \AG \to \AG$ with the property that $\twistO n(G) = (\Z/k)^n \times \Loops^n(G)$ (for $k$ large enough) and a natural transformation $\twistO n \Rightarrow \Id_{\AG}$ that extends the evaluation map natural transformation.

Fix a prime $p$ and let $\AF$ be the Burnside category of saturated fusion systems. The objects of $\AF$ are saturated fusion systems and, given two saturated fusion systems $\cF$ on a $p$-group $S$ and $\cG$ on a $p$-group $T$, the morphisms $\AF(\cF,\cG)$ is the submodule of $\AG(S,T)^{\wedge}_{p}$ consisting of bistable elements. Just as in \cite{RSS_Bold1}*{Section 3}, we may formally extend this category to include formal coproducts of fusion systems, which we call fusoids.

The category $\AF$ is a full subcategory of the homotopy category of $p$-complete spectra. Let $\Sinfpp BS$ be the $p$-completion of the classifying spectrum $\Sinfp BS$. Given a saturated fusion system $\cF$ on a $p$-group $S$, there is a characteristic idempotent 
\[
\omega_{\cF} \in \AG(S,S)^{\wedge}_{p} \cong [\Sinfpp BS,\Sinfpp BS].
\]
We will denote the summand of $\Sinfpp BS$ split off by this idempotent by $\Sinfpp B\cF$. There is a canonical isomorphism
\[
\AF(\cF,\cG) \cong [\Sinfpp B\cF, \Sinfpp B\cG].
\]
Recall that if the fusion system $\cF$ comes from a finite group $G$ (with Sylow $p$-subgroup $S$), then $\Sinfpp B\cF$ is equivalent to the $p$-completion of the classifying spectrum of $G$. Thus $\AF$ contains the homotopy category of $p$-completed classifying spectra of finite groups as a full subcategory.

Just as with groups, there is a notion of a centralizer fusion system of an element in the underlying $p$-group of a fusion system and we define
\[
\Loops \cF = \coprod_{[g]} C_\cF(g). 
\]
It is important to note that $\Loops \cF$ generally has fewer components than $\Loops S$ because $\cF$ generally has fewer conjugacy classes. Consequently $\Loops S$ is not the ``Sylow $p$-subgroup'' of $\Loops \cF$ even though $S$ is the Sylow subgroup of $\cF$. Other technical challenges that we address stem from the fact that the notion of a centralizer fusion system only behaves well for certain choices of representatives for the conjugacy classes in $\cF$.

For $e$ large enough, the evaluation map is a map of fusion systems
\[
\Z/p^e \times (\coprod_{[g]} C_\cF(g)) \to \cF.
\]
Taking $p$-completed classifying spectra, we get a map
\[
\Sinfpp (B\Z/p^e \times \Loops B\cF) \to \Sinfpp B\cF.
\]
Here the free loop space $\Loops B\cF$ is modeled algebraically by $\Loops \cF$ as stated in Proposition \ref{propFreeLoopSpaceDef}, which is a result essentially due to Broto--Levi--Oliver \cite{BLO2}.

We show that the domain of the evaluation map above arises by two other constructions. First, applying $\Id_{\Z/p^e} \times \Loops(-)$ to the characteristic idempotent gives an idempotent
\[
\Sinfpp (B\Z/p^e \times \Loops BS) \xrightarrow{\Id_{\Z/p^e} \times \Loops \omega_{\cF}} \Sinfpp (B\Z/p^e \times \Loops BS).
\]
Secondly, we may apply $\twistO 1$ to the characteristic idempotent $\omega_{\cF}$ to get an idempotent
\[
\Sinfpp (B\Z/p^e \times \Loops BS) \xrightarrow{\twistO 1 \omega_{\cF}} \Sinfpp (B\Z/p^e \times \Loops BS).
\]
Neither of these idempotents are the characteristic idempotent for $\Z/p^e\x \Loops \cF$, but we prove that both of these idempotents still split off the spectrum
\[
\Sinfpp (B\Z/p^e \times \Loops B\cF).
\]

Since each map between fusion systems arises as a map between classifying spectra of $p$-groups that commutes with the characteristic idempotent, we get a natural transformation 
\[
\twistO 1 (-) \Rightarrow \Id_{\AF}
\]
of functors from $\AF$ to $\AF$. Or, said another way, we have extended the functoriality of the evaluation map to the homotopy category of classifying spectra of fusion systems and, in particular, to the homotopy category of $p$-completed classifying spectra of finite groups.

From here, Theorem \ref{thmIntroMain} is proved for fusion systems in a reasonably formal manner as Theorem \ref{thmFusionMain}. We also derive formulas for $\twistO n$ and the evaluation map. These formulas are similar to those described for finite groups in \cite{RSS_Bold1}*{Section 3}.

\begin{theorem}\label{thmIntroMain}
We construct a family of endofunctors $\twistO n \colon \AF\to \AF$ for $n\geq 0$ with the following properties:
\begin{enumerate}
\renewcommand{\theenumi}{$(\roman{enumi})$}\renewcommand{\labelenumi}{\theenumi}
\item[$(\emptyset )$] Let $L^{\dagger,\AG}_n\colon \AG \to \AG$ be the functor constructed in \cite{RSS_Bold1}*{Section 3}. When restricted to the full subcategories of $\AG$ and $\AF$ spanned by formal unions of finite $p$-groups, the functor $\twistO n\colon \AF \to \AF$ is the $\Z_p$-linearization of $L^{\dagger,\AG}_n$.
\item\label{itemIntroFusionLnZero} $\twistO 0$ is the identity functor on $\AF$.
\item\label{itemIntroFusionLnObjects} On objects, $\twistO n$ takes a saturated fusoid $\cF$ to the saturated fusoid
\[\twistO n(\cF) = (\Z/p^e)^n\x \freeO n (\cF)=\coprod_{\tup a\in \cntuples n\cF} (\Z/p^e)^n\x C_\cF(\tup a).\]
\item\label{itemIntroFusionEquivariant} The group $\Sigma_n$ acts on $\freeO n\cF=\coprod_{\tup a\in \cntuples n\cF} C_\cF(\tup a)$ by permuting the coordinates of the $n$-tuples $\tup a$. Explicitly, if $\sigma\in \Sigma_n$ and if $\widetilde{\sigma(\tup a)}$ is the representative for the $\cF$-conjugacy class of $\sigma(\tup a)$, then $\sigma\colon \freeO n\cF\to \freeO n\cF$ maps $C_\cF(\tup a)=C_\cF(\sigma(\tup a))$ to $C_\cF(\widetilde{\sigma(\tup a)})$ via the isomorphism $\change_{\sigma(\tup a)}^{\widetilde{\sigma(\tup a)}}\in \AF(C_\cF(\tup a),C_\cF(\widetilde{\sigma(\tup a)}))$. 

The functor $\twistO n$ is equivariant with respect to the $\Sigma_n$-action on $(\Z/p^e)^n\x \freeO n(-)$ that permutes the coordinates of both $(\Z/p^e)^n$ and $\freeO n(-)$, i.e. for every $\sigma\in\Sigma_n$ the diagonal action of $\sigma$ on $(\Z/p^e)^n\x \freeO n(-)$ induces a natural isomorphism $\sigma\colon \twistO n\overset\cong\Rightarrow \twistO n$.
\item\label{itemIntroFusionLnForwardMaps} Let $\cE$ and $\cF$ be saturated fusion systems on $R$ and $S$ respectively. For forward maps, i.e. transitive bisets $[R,\ph]_\cE^\cF\in \AF(\cE,\cF)$ with $\ph\colon R\to S$ fusion preserving, the functor $\twistO n$ coincides with $(\Z/p^e)^n \x \freeO n(-)$ so that 
\[\twistO n([R,\ph]_\cE^\cF)=(\Z/p^e)^n\x \freeO n([R,\ph]_\cE^\cF).\]
In addition, $\freeO n([R,\ph]_\cE^\cF)$ is the biset matrix that takes a component $C_\cE(\tup a)$ of $\freeO n\cE$ to the component $C_\cF(\tup b)$ of $\freeO n\cF$ by the biset 
\[
[C_R(\tup a),\ph]_{C_\cE(\tup a)}^{C_S(\ph(\tup a))}\cmp \change_{\ph(\tup a)}^{\tup b}\in \AF( C_\cE(\tup a), C_\cF(\tup b)),
\]
where $\tup b$ represents the $\cF$-conjugacy class of $\ph(\tup a)$.
\item\label{itemIntroFusionEvalSquare} For all $n \geq 0$, the functor $\twistO n$ commutes with evaluation maps, i.e. the evaluation maps $\ev_\cF\colon (\Z/p^e)^n\x \freeO n(\cF)\to \cF$ form a natural transformation $\ev\colon \twistO n \Rightarrow \Id_{\AF}$.
\item\label{itemIntroFusionLnPartialEvaluation} For all $n \geq 0$, the partial evaluation maps $\pev_\cF\colon \Z/p^e\x \freeO {n+1}(\cF) \to \freeO n (\cF)$ given as fusion preserving maps $\pev_{\tup a}\colon \Z/p^e \x C_\cF(\tup a) \to C_\cF(a_1,\dotsc, a_n)$ in terms of the formula
\[
\pev_{\tup a}(t,z)= (a_{n+1})^t\cdot z\in C_S(a_1,\dotsc,a_n), \quad \text{for } t\in \Z/p^e,  z\in C_S(\tup a),
\] form natural transformations $(\Z/p^e)^n\x \pev\colon \twistO {n+1} \Rightarrow \twistO n$.
\item\label{itemIntroFusionIterateLn} For all $n,m\geq 0$, and any saturated fusoid $\cF$ on $S$, the formal union $(\Z/p^e)^{n+m}\x \freeO {n+m} \cF$ embeds into $(\Z/p^e)^m\x \freeO m((\Z/p^e)^n\x \freeO n \cF)$ as the components corresponding to the commuting $m$-tuples in $(\Z/p^e)^n\x \freeO n \cF$ that are zero in the $(\Z/p^e)^n$-coordinate, i.e. the embedding takes each component $(\Z/p^e)^{n+m}\x C_\cF(\tup x,\tup y)$ to the component $(\Z/p^e)^m \x C_{(\Z/p^e)^n\x C_\cF(\tup x)}(\tup 0\x \tup y)$, for $\tup x\in \cntuples n \cF$ and $\tup y\in \cntuples m \cF$, via the fusion preserving map given by
\[
((\tup s,\tup r),z) \mapsto (\tup r,(\tup s,z)),
\]
for $\tup s\in (\Z/p^e)^n$, $\tup r\in (\Z/p^e)^m$, and $z\in C_S(\tup x,\tup y)$.

These embeddings $(\Z/p^e)^{n+m}\x \freeO {n+m} \cF\to (\Z/p^e)^m\x \freeO m((\Z/p^e)^n\x \freeO n \cF)$ then form a natural transformation $\twistO {n+m}(-)\Rightarrow \twistO m(\twistO n(-))$.
\end{enumerate}
\end{theorem}

For any finite group $G$, if we work with the free loop space $\freeO{} G$ in the context of a fixed prime $p$, we may want to restrict our view to only those components of $\freeO{} G$ that correspond to conjugacy classes of $p$-power order elements in $G$.
In general for $n\geq 0$, we let $\freeO n_p G$ consist of the components in $\freeO n G$ corresponding to conjugacy classes of commuting $n$-tuples of $p$-power order elements in $G$.

The functors $\twistO n\colon \AG \to \AG$ from \cite{RSS_Bold1}*{Section 3} restrict naturally from $(Z/\ell)^n\times \freeO n G$ to $(\Z/\ell)^n\times \freeO n_p G$ and give us functors $\twistO {n,p} \colon \AG\to \AG$ for $n\geq 0$.

Separately, consider the canonical functor $\AG \to \Ho(\Sp)$ taking finite groups to their classifying spectra and post-compose this with the $p$-completion functor for spectra $(-)^\wedge_p\colon\Ho(\Sp)\to \Ho(\Spp)$. The resulting functor factors through $\AF$ giving us a functor $(-)^\wedge_p\colon \AG\to \AF$ corresponding to the $p$-completion functor for classifying spectra.
The $(-)^\wedge_p$ takes a finite group to its associated fusion system at the prime $p$, and explicit formulas for the functor $(-)^\wedge_p$ were previously described in \cite{RSS_p-completion}.

Our final result describes the interplay between these functors. We prove the functor $\twistO n\colon \AF\to \AF$, when applied to a fusion system coming from a finite group, is in essence the $p$-completion of the functor $\twistO {n,p}\colon \AG\to \AG$.
\begin{theorem}\label{thmIntroMainTheoremPCompletion}
We have $(-)^\wedge_p \circ \twistO {n,p} = \twistO n \circ (-)^\wedge_p$
as functors $\AG\to \AF$ for all $n\geq 0$.
\end{theorem}

\subsection*{Outline}
Section \ref{secFusionRecollections} recalls notation and terminology for saturated fusion systems, their Burnside modules, and characteristic idempotents.
Section \ref{secFusionSystems} introduces the technology needed for working with bisets between free loop spaces of fusion systems.
Section \ref{secFusionFreeLoopFunctor} proves that the mapping telescope of the idempotent $\freeO n(\omega_\cF)$ is equivalent to the classifying spectrum for the free loop space of the fusion system.
Section \ref{secFusionEvaluation} proves a similar result for the mapping telescope of $\twistO n(\omega_\cF)$. The proof reduces to the results of Section \ref{secFusionFreeLoopFunctor} by showing that $\twistO n(\omega_\cF)$ and $(\Z/p^e)^n\times \freeO n(\omega_\cF)$ induce equivalent mapping telescopes. Section \ref{secFusionEvaluation} provides explicit formulas for $\twistO n$ applied to the Burnside category of saturated fusion systems.
Section \ref{secFusionMainTheorem} proves Theorem \ref{thmIntroMain} as Theorem \ref{thmFusionMain}, and
finally, Section \ref{secFusionPCompletionCommutes} proves Theorem \ref{thmIntroMainTheoremPCompletion} as Theorem \ref{thmTwistedLoopsAndPCompletion}.

\subsection*{Acknowledgements} It is a pleasure to thank Mike Hopkins and Haynes Miller for their suggestions. Their helpful comments led to a complete revision of these papers. Further, we thank Cary Malkiewich for generously sharing his thoughts regarding transfers and free loop spaces after his work with Lind in \cite{LindMalkiewich}. Section 2 in the prequel is based on ideas coming from discussions with him. Finally, we thank Matthew Gelvin and Erg\"un Yal\c{c}{\i}n for their comments.

While working on this paper the first author was funded by the Independent Research Fund Denmark (DFF–4002-00224) and later on by BGSMath and the María de Maeztu Programme (MDM–2014-0445). The second and third author were jointly supported by the US-Israel Binational Science Foundation under grant 2018389. The second author was partially supported by the Alon Fellowship and ISF 1588/18. The third author was partially supported by NSF grant DMS-1906236 and the SFB 1085 \emph{Higher Invariants} at the University of Regensburg. All three authors thank the HIM and the MPIM for their hospitality.

\section{Recollections about fusion systems}\label{secFusionRecollections}
We first recall the basics of the definition of a saturated fusion system. For additional details see \cite{RagnarssonStancu}*{Section 2} or \cite{AKO}*{Part I}. 

\begin{definition}
A fusion system on a finite $p$-group $S$ is a category $\cF$ with the subgroups of $S$ as objects and where the morphisms $\cF(P,Q)$ for $P,Q\leq S$ satisfy
\begin{enumerate}
\renewcommand{\theenumi}{$(\roman{enumi})$}\renewcommand{\labelenumi}{\theenumi}
\item Every morphism $\ph\in \cF(P,Q)$ is an injective group homomorphism $\ph\colon P\to Q$.
\item Every map $\ph\colon P\to Q$ induced by conjugation in $S$ is in $\cF(P,Q)$.
\item Every map $\ph\in\cF(P,Q)$ factors as $P\xrightarrow{\ph} \ph(P) \xto{\incl} Q$ in $\cF$ and the inverse isomorphism $\ph^{-1}\colon \ph(P)\to P$ is also in $\cF$.
\end{enumerate}
In addition, the underlying group $S$ is considered part of the structure of the fusion system, so a fusion system is really a pair $(S,\cF)$ of a $p$-group equipped with a category as above.
\end{definition}

We think of the morphisms in $\cF$ as being conjugation maps induced by some, possibly non-existent, ambient group. Consequently, we say that two subgroups $P,Q\leq S$ are conjugate in $\cF$ if there is an isomorphism between them in $\cF$.

A saturated fusion system satisfies some additional axioms that we will not go through as they play almost no direct role in this paper (see e.g. \cite{RagnarssonStancu}*{Definition 2.5} instead). The important aspect of saturated fusion systems in this paper is that these are the fusion systems that have characteristic idempotents and classifying spectra as described below.

Given fusion systems $\cF_1$ and $\cF_2$ on $p$-groups $S_1$ and $S_2$, respectively, a group homomorphism $\phi\colon S_1\to S_2$ is said to be fusion preserving if whenever $\psi\colon P\to Q$ is a map in $\cF_1$, there is a corresponding map $\rho\colon \phi(P)\to \phi(Q)$ in $\cF_2$ such that $\phi|_{Q}\circ \psi = \rho\circ \phi|_P$. Note that each such $\rho$ is unique if it exists.

\begin{example}
Let $G$ be a finite group with Sylow $p$-subgroup $S$. This data determines a fusion system $\cF_G$ on $S$. The maps in $\cF_G(P,Q)$ for subgroups $P,Q\leq S$ are precisely the homomorphisms $P\to Q$ induced by conjugation in $G$, i.e. if $g^{-1} P g\leq Q$ for some $g\in G$, then $c_g(x) = g^{-1}xg$ defines a homomorphism $c_g\in \cF_G(P,Q)$. Note that different group elements $g$ and $g'$ can give rise to the same homomorphism $c_g=c_{g'}$ in $\cF_G(P,Q)$.
The fusion system $\cF_G$ associated to a finite group at a prime $p$ is always saturated.
\end{example}

Every saturated fusion system $\cF$ has a classifying spectrum originally constructed by Broto-Levi-Oliver in \cite{BLO2}*{Section 5}. The most direct way of constructing this spectrum, due to Ragnarsson \cite{Ragnarsson}, is as the mapping telescope
\[
\Sinfp B\cF = \colim ( \Sinfp BS \xto{\omega_{\cF}} \Sinfp BS \xto{\omega_{\cF}} \dotsb),
\]
where $\omega_{\cF} \in \AG(S,S)^\wedge_p$ is the characteristic idempotent of $\cF$ (see the characterization below). By construction, $\Sinfp B\cF$ is a wedge summand of $\Sinfp BS$. 

As remarked in Section 5 of \cite{BLO2}, the spectrum $\Sinfp B\cF$ constructed this way is in fact the suspension spectrum for the classifying space $B\cF$ defined in \cites{BLO2, Chermak}. One way to see this is to note that $H^*(B\cF,\F_p)$ coincides with $\omega_\cF \cdot H^*(BS,\F_p)$ as the $\cF$-stable elements, and that the suspension spectrum of $B\cF$ is $H\F_p$-local.

\begin{definition}\label{defFCharacteristic}
Let $\cF$ be a fusion system on $S$.
A virtual biset $X\in \AG(S,S)^\wedge_p$ is said to be \emph{$\cF$-characteristic} if it satisfies the following properties:
\begin{itemize}
\item $X$ is \emph{$\cF$-generated}, meaning that $X$ is a linear combination of transitive bisets $[P,\ph]_S^S$ with $\ph\in \cF(P,S)$.
\item $X$ is \emph{left $\cF$-stable}, meaning that for all $P\leq S$ and $\ph\in \cF(P,S)$ the restriction of $X$ along $\ph$ on the left,
    \[X_{P,\ph}^S = [P,\ph]_P^S \cmp X \in \AG(P,S)^\wedge_p,\]
    is isomorphic to the restriction $X_P^S= [P,\id]_P^S \cmp X\in \AG(P,S)^\wedge_p$ along the inclusion.
\item $X$ is \emph{right $\cF$-stable}, meaning the dual property to the one above:
    \[X\cmp [\ph P, \ph^{-1}]_S^P = X\cmp [P, \id]_S^P,\]
    for all $P\leq S$ and $\ph\in \cF(P,S)$.
\item $\abs{X}/\abs{S}$ is invertible in $\ZZ_p$.
\end{itemize}

Whenever $\cF$ is the fusion system generated by a finite group $G$ with Sylow $p$-subgroup $S$, we can consider $G$ itself as an $(S,S)$-biset $G_S^S$. The biset $G_S^S$ is always $\cF$-characteristic and is the motivating example for Definition \ref{defFCharacteristic}.

According to \cite{RagnarssonStancu}, a fusion system is saturated if and only if it has a characteristic virtual biset, in which case there is a unique idempotent $\omega_\cF$ among all the characteristic virtual bisets. An explicit description and construction of $\omega_\cF$, and a classification of all the characteristic virtual bisets, can be found in \cite{ReehIdempotent}.
\end{definition}

\begin{definition}\label{defFusionBisetModule}
Let $\cF_1$ and $\cF_2$ be saturated fusion systems over $p$-groups $S_1$ and $S_2$, respectively. We then let $\AF(\cF_1,\cF_2)$ denote the submodule of $\AG(S_1,S_2)^\wedge_p$ consisting of all virtual $(S_1,S_2)$-bisets that are both left $\cF_1$-stable and right $\cF_2$-stable.
Here $\cF_1$- and $\cF_2$-stability can be defined as for the characteristic idempotents above. Alternatively, a virtual biset $X\in \AG(S_1,S_2)^\wedge_p$ is $(\cF_1,\cF_2)$-stable if and only if
\[\omega_{\cF_1} \cmp X \cmp \omega_{\cF_2} = X.\]
\end{definition}

As a special case, for the trivial fusion systems induced by $S_1,S_2$ on themselves, we have $\AF(S_1,S_2) = \AG(S_1,S_2)^\wedge_p$. Hence $\AF(S_1,S_2)$ is the free $\Z_p$-module spanned by the transitive bisets $[P,\ph]_{S_1}^{S_2}$ for $P\leq S_1$ and $\ph\colon P\to S_2$ up to conjugation. Similarly, $\AF(\cF_1,\cF_2)$ is the free $\Z_p$-module on basis elements of the form
\[
[P,\ph]_{\cF_1}^{\cF_2} := \omega_{\cF_1} \cmp [P,\ph]_{S_1}^{S_2} \cmp \omega_{\cF_2},
\]
where $P\leq S_1$ and $\ph\colon P\to S_2$ are taken up to pre- and postcomposition with isomorphisms in $\cF_1$ and $\cF_2$ respectively (see \cite{Ragnarsson}*{Proposition 5.2}).

Suppose we have a third saturated fusion system $\cF_3$ on $S_3$ and virtual bisets $X\in \AF(\cF_1,\cF_2)$ and $Y\in \AF(\cF_2,\cF_3)$. The composition $X\cmp Y = X\x_{S_2} Y \in \AF(S_1,S_3)$ is then always $(\cF_1,\cF_3)$-stable, so composition gives a well-defined map
\[
\cmp\colon \AF(\cF_1,\cF_2) \x \AF(\cF_2,\cF_3) \to \AF(\cF_1,\cF_3).
\]
Here we use ``right-composition'' similar as for bisets of finite groups in \cite[Definition 3.1]{RSS_Bold1}.

\begin{remark}
Given finite groups $G$ and $H$ with Sylow subgroups $S$ and $T$ and a $(G,H)$-biset $X$, we may restrict the $G$-action to $S$ and the $H$-action to $T$ to get an $(S,T)$-biset $X_S^T$. Let $\cF_G$ be the fusion system associated to $G$ on $S$ and $\cF_H$ the fusion system associated to $H$ on $T$. The restricted biset  $X_S^T$ is always a stable biset and so we may further consider it as an $(\cF_G, \cF_H)$-biset $X_{\cF_G}^{\cF_H}\in \AF(\cF_G,\cF_H)$.
\end{remark}

\begin{convention}\label{conventionFusionOrbitDecomposition}
As with \cite[Convention 3.2]{RSS_Bold1} for groups, we allow flexibility when decomposing virtual bisets into basis elements for fusion systems.

Let $\cF$ and $\cG$ be saturated fusion systems over $p$-groups $S$ and $T$ respectively.
Given $X\in \AF(\cF,\cG)$, we can write $X$ as a $\Z_p$-linear combination of basis elements:
\[\sum_{(R,\ph)} c_{R,\ph}\cdot [R,\ph]_\cF^\cG.\]
The summation runs over all $R\leq S$ and $\ph\colon R\to T$ (not taken up to conjugacy). The coefficient function $c_{(-)}$ is a choice of function from the set of all pairs $(R,\ph)$ to $\Z_p$ such that the sum of coefficients $c_{R',\ph'}$ over the pairs $(\cF,\cG)$-conjugate to $(R,\ph)$ is the number of copies of the basis element $[R,\ph]_\cF^\cG$ in $X$.

As with \cite[Convention 3.2]{RSS_Bold1}, the linear combination is \emph{not} unique as several conjugate pairs $(R,\ph)$ can contribute to the sum at the same time. If we require $c_{(-)}$ to be concentrated on chosen representatives for the conjugacy classes of pairs, then the linear combination is unique.
\end{convention}

\begin{remark}
An advantage of the flexibility in the linear combinations above is that we can use the same coefficients for \cite[Convention 3.2]{RSS_Bold1} and Convention \ref{conventionFusionOrbitDecomposition} at the same time:
Given $X\in \AF(\cF,\cG)$, we first consider the $(S,T)$-biset ${}_S X_T$ and write this as a linear combination according to \cite[Convention 3.2]{RSS_Bold1},
\[
{}_S X_T = \sum_{(R,\ph)} c_{R,\ph}\cdot [R,\ph]_S^T.
\]
Recalling that $X$ is taken to be $(\cF,\cG)$-stable, we can compose with the characteristic idempotents from each side without changing $X$:
\[
X= \omega_\cF \cmp {}_S X_T \cmp \omega_\cG = \sum_{(R,\ph)} c_{R,\ph}\cdot (\omega_\cF\cmp [R,\ph]_S^T\cmp \omega_\cG) = \sum_{(R,\ph)} c_{R,\ph}\cdot [R,\ph]_\cF^\cG.
\]
Hence the coefficients $c_{R,\ph}$ chosen when decomposing ${}_S X_T$ as an $(S,T)$-biset also work when decomposing $X$ as an $(\cF,\cG)$-biset.
\end{remark}

Suppose $\cE$ and $\cF$ are saturated fusion systems on $R$ and $S$ respectively. Because we construct $\Sinfpp B\cF$ for a saturated fusion system $\cF$ on $S$ as the colimit 
\[
\Sinfpp B\cF = \colim(\Sinfpp BS \xto{\omega_\cF} \Sinfpp BS \xto{\omega_\cF} \dotsb)
\]
with respect to the idempotent $\omega_\cF\in \AF(S,S)$, the stable maps from $\Sinfpp B\cE$ to $\Sinfpp B\cF$ are given by
\[
[\Sinfpp B\cE, \Sinfpp B\cF] \cong \omega_{\cE}\cmp \AF(R,S)\cmp \omega_{\cF} =\AF(\cE,\cF).
\]

\section{Free loop spaces for saturated fusion systems}\label{secFusionSystems}

In order to have a framework in which to work with free loop spaces for fusion systems, we introduce formal unions of fusion systems:
\begin{definition}
Suppose we have a formal union of finite $p$-groups $S=S_1\sqcup \dotsb \sqcup S_k$. We define a \emph{fusoid} $\cF$ on $S$ to be a collection of fusion systems $\cF_i$ on $S_i$ for $1\leq i\leq k$, and we write $\cF=\cF_1\sqcup \dotsb \sqcup \cF_k$.

Given another fusoid $\cE$ with underlying union $R$ of $p$-groups, we define a fusion preserving map $\cE\to \cF$ to be a collection of homomorphisms that take each component $R_i$ of $R$ to some component $S_j$ of $S$ and where the homomorphism $R_i\to S_j$ is fusion preserving from $\cE_i$ to $\cF_j$.

A fusoid $\cF$ is \emph{saturated} if each component of $\cF$ is saturated. Furthermore each saturated fusoid has as classifying space $B\cF$ the disjoint union of the classifying spaces of the components, and $\cF$ has a classifying spectrum
$\Sinfpp B\cF = \colim (\Sinfpp BS\xto{\omega_{\cF}} \Sinfpp BS \xto{\omega_{\cF}} \dotsb)$, where $\omega_\cF\in \AG(S,S)^\wedge_p$ is the diagonal matrix with entries $\omega_{\cF_i}\in \AG(S_i,S_i)^\wedge_p$.

The classifying spectrum $\Sinfpp B\cF$ is also the sum under $\mathbb S^\wedge_p$ of the classifying spectra $\Sinfpp B\cF_i$ for the components of $\cF$ where the copy of $\mathbb S^\wedge_p$ in each $\Sinfpp B\cF_i$ coming from the disjoint basepoints are identified with each other.
\end{definition}

\begin{definition}
For saturated fusoids $\cE$ and $\cF$ over unions $R$ and $S$ of $p$-groups, we define $\AF(\cE,\cF)$ similarly as for unions of groups, so that each $X\in \AF(\cE,\cF)$ is a matrix of virtual bisets with entries $X_{i,j} \in \AF(\cE_i,\cF_j)$ for the corresponding components of $\cE$ and $\cF$. The $\Z_p$-module $\AF(\cE,\cF)$ is a $\Z_p$-submodule of the module $\AF(R,S)$ of matrices for the underlying unions of $p$-groups.

Again we define composition $\cmp$ in terms of matrix multiplication, and we let \emph{the biset category of fusoids} $\AF$ be the category with objects the saturated fusoids at the prime $p$ and morphism set from $\cE$ to $\cF$ given by $\AF(\cE,\cF)$.

The identity map $\id_{\cF}\in \AF(\cF,\cF)$ for a saturated fusoid is just the diagonal matrix $\omega_{\cF}$ with entry $\omega_{\cF_i}$ for each component $\cF_i$ of $\cF$.
\end{definition}

There is a functor $\AF\to \Ho(\Spp)$ that takes a fusoid $\cF$ to the $p$-completed classifying spectrum $\Sinfpp B\cF \simeq \mathbb S^\wedge_p \vee \Sinf B\cF$, and on morphisms it is the Segal map for fusoids, which is an isomorphism:
\[[\Sinfpp B\cE,\Sinfpp B\cF] \cong \AF(\cE,\cF).\]
As such, $\AF\to \Ho(\Spp)$ is fully faithful.

A word of caution: While the restriction of actions ${}_G X_H\mapsto {}_{\cF_G}X_{\cF_H}$ defines a map $\AG(G,H)\to \AF(\cF_G,\cF_H)$ for any finite groups $G$ and $H$, this does not define a functor $\AG\to \AF$. See \cite{RSS_p-completion}*{Theorem 1.1} for a description of the functor $\AG\to \AF$ that corresponds to $p$-completion of spectra.

We wish to give an algebraic model for the $n$-fold free loop space $\freeO n(B\cF)$, when $\cF$ is a saturated fusion system or fusoid. In order to do this, we first need to specify what we mean by commuting $n$-tuples in $\cF$, their conjugacy classes, and their centralizer fusion systems.

\begin{definition}
Let $\cF$ be a saturated fusoid or fusion system on a union $S$ of finite $p$-groups. For each $n\geq 1$, we consider
$n$-tuples $\tup a = (a_1,\dotsc, a_n)$ of commuting elements in $S$. Note that the elements of a tuple $\tup a$ are required to lie in the same component of the formal union $S$. We say that two $n$-tuples $\tup a$ and $\tup b$ are \emph{$\cF$-conjugate} if they lie in the same component of $S$ and there is a map in $\cF$ sending $\tup a$ to $\tup b$, i.e. a map
\[\ph\colon \gen{\tup a} = \gen{a_1,\dotsc,a_n} \to \gen{\tup b}=\gen{b_1,\dotsc,b_n} \text{ in $\cF$}\]
such that $\ph(a_i)=b_i$.

We let $\cntuples n\cF$ denote the collection of equivalence classes of commuting $n$-tuples in $S$ up to $\cF$-conjugation, for $n>0$. For $n=0$, we consider $\cntuples 0\cF$ to consist of a single empty/trivial $0$-tuple.

For a finite group $G$ and $n\geq 0$, let $\cnptuples nG$ denote the classes of commuting $n$-tuples of \emph{elements with $p$-power order} in $G$ up to $G$-conjugation.
\end{definition}

\begin{definition}
Let $\cF$ be a saturated fusion system or fusoid over $S$. We say that an $n$-tuple $\tup a$ in $S$ is fully $\cF$-centralized if $\abs{C_S(\tup a)}\geq \abs{C_S(\tup a')}$ for all $n$-tuples $\tup a'$ conjugate to $\tup a$ in $\cF$, which is the case if and only if the subgroup $\gen{a_1,\dotsc,a_n}$ is fully $\cF$-centralized in the normal terminology of fusion systems.

When $\cF$ is a saturated fusion system, and $\tup a$ is fully centralized in $\cF$, we define the centralizer fusion system $C_\cF(\tup a)$ to be the fusion system over the $p$-group $C_S(\tup a)$ with maps
\begin{multline*}
\Hom_{C_\cF(\tup a)}(Q,P)=\{\ph\in \cF(Q,P) \mid \text{$\ph$ extends to a map $\tilde\ph\in \cF$ defined on }\\\text{ $Q$ and each $a_i$, $1\leq i\leq n$, such that $\tilde\ph|_Q = \ph$ and $\tilde\ph(a_i)=a_i$}\}
\end{multline*}
for subgroups $Q,P\leq C_S(\tup a)$.
The fusion system $C_\cF(\tup a)$ coincides with the usual notion of the centralizer fusion system for the subgroup $\gen{a_1,\dotsc,a_n}\leq S$.

When $\cF$ is a fusoid, and $\tup a$ is fully $\cF$-centralized, we define the centralizer fusion system $C_\cF(\tup a)$ to be the centralizer inside the component of $\cF$ containing $\tup a$. As such, the centralizer $C_\cF(\tup a)$ is always a fusion system and not a fusoid.
\end{definition}

\begin{lemma}\label{lemmaSaturatedCentralizers}
Suppose $\tup a$ is a fully $\cF$-centralized $n$-tuple, then the centralizer fusion system $C_\cF(\tup a)$ is saturated.
\end{lemma}

\begin{proof}
This is just \cite{BLO2}*{Proposition A.6} applied to the subgroup $\gen{a_1,\dotsc,a_n}$ and the centralizer system $C_\cF(\gen{a_1,\dotsc,a_n})$.
\end{proof}

\begin{remark}\label{remarkExtendingToCentralizers}
The following fact, which we will need for the next lemma, is a special case of the second saturation axiom \cite{RagnarssonStancu}*{Definition 2.5(II)}.

If a commuting $n$-tuple $\tup a$ is fully $\cF$-centralized in a saturated fusion system $\cF$, and if $\ph$ in $\cF$ takes any other $n$-tuple $\tup b$ to $\tup a$, then the map $\ph\colon \gen{b_1,\dotsc,b_n}\to \gen{a_1,\dotsc,a_n}$ extends to a map between centralizers $\tilde\ph\colon C_S(\tup b)\to C_S(\tup a)$ with $\tilde\ph(b_i)=a_i$.
\end{remark}

\begin{lemma}\label{lemmaIteratedRepr}
Let $\cF$ be a saturated fusion system on $S$, and let $\tup a\in S$ be a commuting $(n+1)$-tuple. Write $\tup a= (a_1,\dotsc,a_{n+1})$. Suppose $\ptup a= (a_1,\dotsc,a_n)$ is fully centralized in $\cF$, and suppose further that $a_{n+1}\in C_S(\ptup a)$ is fully centralized in $C_\cF(\ptup a)$, then $\tup a$ is fully centralized in $\cF$.

Furthermore, each commuting $(n+1)$-tuple $\tup b$ in $S$ is $\cF$-conjugate to a fully centralized tuple of the form above.
\end{lemma}

\begin{proof}
First of all, for any $(n+1)$-tuple $\tup a$ the element $a_{n+1}$ commutes with $\ptup a$ if and only if $a_{n+1}\in C_S(\ptup a)$, and in that case
\[C_S(\tup a) =C_{C_S(\ptup a)}(a_{n+1}).\]

Suppose $\tup a$ is a commuting $(n+1)$-tuple with $\ptup a$ fully $\cF$-centralized and $a_{n+1}$ fully $C_\cF(\ptup a)$-centralized.
Consider any other $(n+1)$-tuple $\tup b$ that is $\cF$-conjugate to $\tup a$, and suppose $\ph\in \cF$ sends $\tup b$ to $\tup a$. Let $\ptup b=(b_1,\dotsc,b_n)$ and let $\ph_\partial$ be the restriction of $\ph$ to the subgroup generated by $\ptup b$, so that $\ph_\partial(\ptup b)=\ptup a$. Since $\ptup a$ is fully $\cF$-centralized, Remark \ref{remarkExtendingToCentralizers} implies that $\ph_\partial$ extends to the centralizers as a map $\psi\colon C_S(\ptup b) \to C_S(\ptup a)$. In particular $\psi(b_{n+1})\in C_S(\ptup a)$.
All elements of $C_S(\tup b)\leq C_S(\ptup b)$ centralize $b_{n+1}$, so after applying $\psi$ we get
\[\psi(C_S(\tup b)) \leq C_{C_S(\ptup a)}(\psi(b_{n+1})).\]
We proceed to look at the composite
\[\ph\circ \psi^{-1} \colon \gen{a_1,\dotsc,a_n,\psi (b_{n+1})} \to \gen{b_1,\dotsc,b_{n+1}} \to \gen{a_1,\dotsc,a_{n+1}}.\]
Note that $\ph\circ \psi^{-1}$ maps $\psi(b_{n+1})$ to $a_{n+1}$ and at the same time maps $\ptup a$ identically to itself. The composite $\ph\circ \psi^{-1}$ therefore defines a map in $C_\cF(\ptup a)$ from $\psi(b_{n+1})$ to $a_{n+1}$. Hence $a_{n+1}$ and $\psi(b_{n+1})$ are conjugate in $C_\cF(\ptup a)$, wherein $a_{n+1}$ was assumed to be fully centralized, so we conclude that
\[\abs{C_S(\tup a)} = \abs{C_{C_S(\ptup a)}(a_{n+1})} \geq \abs{C_{C_S(\ptup a)}(\psi(b_{n+1}))} \geq \abs{\psi(C_S(\tup b))} =\abs{C_S(\tup b)}.\]
This completes the proof that $\tup a$ is in fact fully centralized in $\cF$.

For the second part of the lemma, let $\tup b$ be any commuting $(n+1)$-tuple in $S$. Consider the truncated $n$-tuple $\ptup b=(b_1,\dotsc,b_n)$, and choose any preferred fully $\cF$-centralized conjugate $\tup a$ of $\ptup b$. Let $\ph$ be a map in $\cF$ from $\ptup b$ to $\tup a$, then by Remark \ref{remarkExtendingToCentralizers} as above we get an extension of $\ph$ to $\tilde\ph\colon C_S(\ptup b)\to C_S(\tup a)$. We have $b_{n+1}\in C_S(\ptup b)$ and $\ph(b_{n+1})\in C_S(\tup a)$. Choose any $z\in C_S(\tup a)$ that is fully $C_\cF(\tup a)$-centralized and conjugate to $\ph(b_{n+1})$ inside $C_\cF(\tup a)$. Any $\psi \in C_\cF(\tup a)$ that takes $\ph(b_{n+1})$ to $z$ then $\psi$ extends trivially onto $\tup a$, hence $\psi\circ \ph$ takes the entire tuple $\tup b$ to $(a_1,\dotsc,a_n,z)$. Furthermore $(a_1,\dotsc,a_n,z)$ has the requested form -- and hence is fully $\cF$-centralized by the lemma.
\end{proof}

We wish to take the formal union of centralizer fusion systems $C_\cF(\tup a)$ with $[\tup a]\in \cntuples n\cF$ as an algebraic model for the $n$-fold free loop space of $B\cF$. To show that the centralizers do not depend on the choice of representatives, we give the following analogue of \cite[Lemma 3.12]{RSS_Bold1}.

\begin{lemma}\label{lemmaCentralizerInclusion}
Let $\tup a$ be an $n$-tuple of commuting elements in $S$, and suppose that $\reptup a$ is fully $\cF$-centralized and $\cF$-conjugate to $\tup a$. Any map $\ph$ in $\cF$ that takes $\tup a$ to $\reptup a$, then induces a fusion preserving injective map $C_\cF(\tup a)\into C_\cF(\reptup a)$. In $\AF$, all such inclusions give rise to the same virtual bifree biset
\[\change_{\tup a}^{\reptup a} \in \AF(C_S(\tup a), C_\cF(\reptup a))\]
that is left $C_\cF(\tup a)$-stable in addition to being right $C_\cF(\reptup a)$-stable.

If $\reptup a'$ is another fully $\cF$-centralized tuple that is conjugate to $\tup a$ (and therefore to $\reptup a$), then
\[\change_{\reptup a}^{\reptup a'}\in \AF(C_\cF(\reptup a), C_\cF(\reptup a')),\]
and the chosen bisets are compatible with composition
\[\change_{\tup a}^{\reptup a'} = \change_{\tup a}^{\reptup a} \cmp \change_{\reptup a}^{\reptup a'}.\]
\end{lemma}

\begin{proof}
The $\cF$-conjugation from $\tup a$ to $\reptup a$ is given by a unique map $\ph\colon \gen{\tup a} \to \gen{\reptup a}$ in $\cF$. 
By Remark \ref{remarkExtendingToCentralizers} the map $\ph$ extends to a map
\[\tilde\ph\colon C_S(\tup a) \to C_S(\reptup a)\]
such that $\tilde\ph|_{\gen{\tup a}} = \ph$. We define
\[\change_{\tup a}^{\reptup a} := [C_S(\tup a), \tilde\ph]_{C_S(\tup a)}^{C_S(\reptup a)} \cmp \omega_{C_\cF(\reptup a)} = [C_S(\tup a), \tilde\ph]_{C_S(\tup a)}^{C_\cF(\reptup a)}.\]
Given any other choice of extension $\psi \colon C_S(\tup a) \to C_S(\reptup a)$ that maps $\tup a$ to $\reptup a$, the composite $\psi\circ (\tilde \ph)^{-1}$ defines a map in $\cF$ from $\tilde\ph(C_S(\tup a))$ to $\psi(C_S(\tup a))$. The composite $\psi\circ (\tilde \ph)^{-1}$ maps the fully centralized tuple $\reptup a$ to itself by the identity, so $\rho=\psi\circ (\tilde \ph)^{-1}$ defines a map in $C_\cF(\reptup a)$ from the subgroup $\tilde\ph(C_S(\tup a))$ to $\psi(C_S(\tup a))$. Since $\psi=\rho \circ \tilde\ph$ and $\rho\in C_\cF(\reptup a)$, the two maps $\psi$ and $\tilde\ph$ give rise to the same virtual biset
\[[C_S(\tup a), \psi]_{C_S(\tup a)}^{C_\cF(\reptup a)} = [C_S(\tup a), \rho\circ \tilde\ph]_{C_S(\tup a)}^{C_\cF(\reptup a)} = [C_S(\tup a), \tilde\ph]_{C_S(\tup a)}^{C_\cF(\reptup a)}.\]
The virtual biset $\change_{\tup a}^{\reptup a}$ is therefore independent of the choice of extension $\tilde\ph$ of the map $\ph$.

We now claim that $\tilde\ph\colon C_S(\tup a)\to C_S(\reptup a)$ is in fact fusion preserving from $C_\cF(\tup a)$ to $C_\cF(\reptup a)$. Given $\zeta\in C_\cF(\tup a)(Q,P)$, it extends to $\tilde\zeta\colon \gen{\tup a}Q \to \gen{\tup a}P$ sending the elements of $\tup a$ to themselves. The composite $\tilde\ph\circ \zeta \circ (\tilde\ph)^{-1}\colon \tilde\ph(Q)\to \tilde\ph(P)$ is a map in $\cF$ that sends $\reptup a$ to itself, hence this composite lies in $C_\cF(\reptup a)$ as required, so $\tilde \ph$ is fusion preserving. Because $\tilde\ph$ is fusion preserving, it follows from \cite{RSS_p-completion}*{Lemma 4.6} that
\[\change_{\tup a}^{\reptup a} = [C_S(\tup a), \tilde\ph]_{C_S(\tup a)}^{C_S(\reptup a)} \cmp \omega_{C_\cF(\reptup a)}\]
is left $C_\cF(\tup a)$-stable.

Given a further map $\theta$ in $\cF$  from $\reptup a$ to $\reptup a'$, we have
\[\change_{\tup a}^{\reptup a} \cmp \change_{\reptup a}^{\reptup a'} = [C_S(\tup a), \tilde\ph]_{C_S(\tup a)}^{C_\cF(\reptup a)} \cmp [C_S(\reptup a), \tilde\theta]_{C_\cF(\reptup a)}^{C_\cF(\reptup a')}  = [C_S(\tup a), \tilde \theta\circ \tilde\ph]_{C_S(\tup a)}^{C_\cF(\reptup a')},\]
by Proposition \ref{propFusionDoubleCosetFormula}, since $\tilde\theta$ is fusion preserving. The composite \[\tilde\theta\circ\tilde\ph\colon C_S(\tup a) \to C_S(\reptup a')\]
sends $\tup a$ to $\reptup a'$ and is therefore a valid choice of extension $\widetilde{\theta\circ\ph}:= \tilde\theta\circ \tilde\ph$. Using this choice for $\widetilde{\theta\circ\ph}$, we then get
\[\change_{\tup a}^{\reptup a'} = [C_S(\tup a), \tilde \theta\circ \tilde\ph]_{C_S(\tup a)}^{C_\cF(\reptup a')} = \change_{\tup a}^{\reptup a} \cmp \change_{\reptup a}^{\reptup a'}.\qedhere\]
\end{proof}

\begin{lemma}\label{lemmaChangeOfRepresentatives}
Suppose $\tup a_1,\dotsc,\tup a_r$ and $\tup b_1,\dotsc,\tup b_r$ are two choices of fully centralized representatives for $\cntuples n\cF$, and suppose the labelling is such that $\tup a_i$ is $\cF$-conjugate to $\tup b_i$ for $1\leq i\leq r$. There is then a canonical isomorphism of formal unions
$\displaystyle\coprod_{i} C_\cF(\tup a_i) \xrightarrow{\cong} \coprod_{i} C_\cF(\tup b_i)$ in $\AF$ via the diagonal matrix whose $i$th entry is $\change_{\tup a_i}^{\tup b_i}\in \AF(C_\cF(\tup a_i),C_\cF(\tup b_i))$.

Given further choices of representatives, the isomorphisms are compatible with respect to composition.
\end{lemma}

\begin{proof}
Let $\iso_a^b$ be the diagonal matrix with entries $\change_{\tup a_i}^{\tup b_i}$. Given a third set of fully centralized representatives $\tup c_1,\dotsc,\tup c_r$ for $\cntuples n\cF$, the fact that the isomorphisms are compatible,
\[\iso_a^c=\iso_a^b\cmp \iso_b^c,\]
is immediate from the composition of diagonal entries $\change_{\tup a_i}^{\tup c_i} = \change_{\tup a_i}^{\tup b_i} \cmp \change_{\tup b_i}^{\tup c_i}$, which follows from Lemma \ref{lemmaCentralizerInclusion}.

The inverse matrix to $\iso_a^b$ is just the diagonal matrix $\iso_b^a$ with diagonal entries $\change_{\tup b_i}^{\tup a_i}$. The fact that that $\iso_a^b$ and $\iso_b^a$ are inverses follows from the equality $\change_{\tup a_i}^{\tup b_i} \cmp \change_{\tup b_i}^{\tup a_i}=\change_{\tup a_i}^{\tup a_i}$. Here $\change_{\tup a_i}^{\tup a_i}$ is the identity element in $A(C_\cF(\tup a_i),C_\cF(\tup a_i))$ since $\change_{\tup a_i}^{\tup a_i}$ is induced by the identity map on $C_S(\tup a_i)$ in $\cF$. The analogous statement is true for the $\tup b$'s.
\end{proof}

\begin{convention}\label{conventionTupleReps}
As with \cite[Convention 3.14]{RSS_Bold1}, we will now suppose that a choice of preferred representatives for $\cntuples n\cF$ has been made for all $n\geq 0$. We require that each representative $n$-tuple $\tup a$ is fully $\cF$-centralized. 
Lemma \ref{lemmaIteratedRepr} furthermore enables us to chose representatives such that each chosen representative $n$-tuple $\tup a=(a_1,\dotsc,a_n)$ satisfies that $\ptup a=(a_1,\dotsc,a_{n-1})$ is one of the previously chosen representative $(n-1)$-tuples (and $a_n$ is fully $C_\cF(\ptup a)$-centralized).
\end{convention}

Since $B\cF$, as constructed by \cite{BLO2} and \cite{Chermak}, is the $p$-completion of a finite category, the usual $n$-fold free loop space $\Map(B(\Z^n),B\cF)$ is equivalent to the colimit over cyclic $p$-groups: 
\[\freeO n B\cF \simeq \colim_{e\to \infty} \Map(B(\ZZ/p^e)^n, B\cF)\]
for any union $\cF$ of saturated fusion systems. In the following we shall replace $S^1$ with the classifying space $B(\Z/p^e)$ for sufficiently large $e$. We will follow \cite{RSS_Bold1}*{Convention 3.21} and suppose $e$ is large enough to work for all the finitely many fusion systems in each calculation.

As in \cite[Definition 3.11]{RSS_Bold1} we introduce an algebraic model for the $n$-fold free loop space of $B\cF$ as a union of centralizer fusion systems:

\begin{definition}\label{defFusionLoop}\label{defFusionAlgLoop}
Let $\cF$ be a saturated fusion system on $S$. We define $\freeO n \cF$ to be the saturated fusoid
\[\freeO n \cF := \coprod_{[\tup a]\in \cntuples n\cF} C_\cF(\tup a),\]
where the chosen representatives are fully centralized (according to Convention \ref{conventionTupleReps}).
\end{definition}

By Lemma \ref{lemmaChangeOfRepresentatives}, different choices of representatives for the conjugacy classes result in isomorphic fusoids $\freeO n \cF$.

Mapping spaces into $B\cF$ are described in detail by Broto-Levi-Oliver in their paper \cite{BLO2}, and applying their results when mapping out of $B(\ZZ/p^e)^n$ essentially gives a proof that $\freeO n\cF$ models the $n$-fold free loop space of $B\cF$.

\begin{prop}[\cite{BLO2}]\label{propFreeLoopSpaceDef}
Let $\cF$ be a saturated fusoid or fusion system. The fusoid $\freeO n\cF$ is an algebraic model for the free loop space $\displaystyle\freeO n B\cF \simeq \colim_{e\to\infty} \Map(B(\ZZ/p^e)^n , B\cF)$, that is we have homotopy equivalences
\[\freeO n B\cF \simeq B(\freeO n \cF) = \coprod_{[\tup a]\in \cntuples n\cF} BC_\cF(\tup a),\]
where the chosen representatives $\tup a$ are fully centralized.
\end{prop}

Before we go through the proof, we make the following observation based on the construction preceding Theorem 6.3 of \cite{BLO2}.

\begin{remark}\label{remarkFreeLoopSpaceDef}
Suppose $\cF$ is a saturated fusion system. The homotopy equivalence in Proposition \ref{propFreeLoopSpaceDef} is given in terms of maps
\[f_{\tup a}\colon B(\ZZ/p^e)^n\x BC_\cF(\tup a)  \to B\cF,\]
for each representing $n$-tuple $\tup a$. Each map $f_{\tup a}$, when restricted to the underlying $p$-groups, is $B(\ev_{\tup a})$, where $\ev_{\tup a}$ is the evaluation map of \cite[Lemma 3.20]{RSS_Bold1}
\[\ev_{\tup a}\colon (\ZZ/p^e)^n\x C_S(\tup a)  \to S\]
given by $\ev_{\tup a}(t_1,\dotsc,t_n,z) := (a_1)^{t_1}\dotsm(a_n)^{t_n}z$. Because $\ev_{\tup a}$ is a homomorphism of $p$-groups, when we pass to the Burnside category $\AF$, the map $f_{\tup a}$ is just
\[ [(\ZZ/p^e)^n\x C_S(\tup a) , \ev_{\tup a}]_{(\ZZ/p^e)^n\x C_\cF(\tup a) }^\cF. \]

We will see below in Lemma \ref{lemmaEvalFusionPreserving} that each $\ev_{\tup a}$ is in fact fusion preserving.
\end{remark}

\begin{proof}[Proof of Proposition \ref{propFreeLoopSpaceDef}]
It is enough to consider the case when $\cF$ is a saturated fusion system over a finite $p$-group $S$.
Let $e\gg0$, in fact $p^e \geq \abs S$ is enough. Corollary 4.5 of \cite{BLO2} tells us that
$[B(\ZZ/p^e)^n, B\cF]$ is in bijection with the classes of commuting $n$-tuples $\cntuples n\cF$.
For each class in $\cntuples n\cF$ choose a fully centralized representative $\tup a$. Then Theorem 6.3 of \cite{BLO2} states that the connected component of $\Map(B(\ZZ/p^e)^n , B\cF)$ corresponding to $\tup a$ is homotopy equivalent to $BC_\cF(\tup a)$. In total this shows that
\[\Map(B(\ZZ/p^e)^n , B\cF) \simeq \coprod_{[\tup a]\in \cntuples n\cF} BC_\cF(\tup a)\]
for sufficiently large $e$.
\end{proof}

\begin{lemma}\label{lemmaEvalFusionPreserving}
Let $\cF$ be a saturated fusion system, and let $\tup a$ be a representative $n$-tuple. The homomorphism $\ev_{\tup a}\colon (\Z/p^e)^n\x C_S(\tup a) \to S$ given by 
\[
\ev_{\tup a}(t_1,\dotsc,t_n,z) := (a_1)^{t_1}\dotsm(a_n)^{t_n}z
\]
is a fusion preserving with respect to the fusion systems $(\Z/p^e)^n\x C_\cF(\tup a)$ and $\cF$.
\end{lemma}

\begin{proof}
Any morphism in the product fusion system $(\Z/p^e)^n\x C_\cF(\tup a)$ has the form $\id\x \ph\colon D\to E$ for subgroups $D,E\leq (\Z/p^e)^n\x C_S(\tup a)$ and where $\ph$ is a morphism in the centralizer fusion system $C_\cF(\tup a)$.
Denote the projections of $D$ by $D_1\leq (\Z/p^e)^n$ and $D_2\leq C_S(\tup a)$ respectively, so $D\leq D_1\x D_2$, and similarly $E\leq E_1\x E_2$. Then $\ph\in \Hom_{C_\cF(\tup a)}(D_2, E_2)$ by definition of a product fusion system.

By definition of the centralizer fusion system, $\ph$ extends to $\tilde\ph\colon \gen{\tup a, D_2}\to \gen{\tup a, E_2}$ with $\tilde\ph(\tup a) = \tup a$ and $\tilde \ph\in \cF$. Now the diagram
\[
\begin{tikzpicture}
\node [matrix of math nodes] (M) {
D &[1cm] D_1\x D_2 &[2cm] \gen{\tup a, D_2} \\[2cm]
E & E_1\x E_2 & \gen{\tup a, E_2} \\
};
\path 
    (M-1-1) -- node{$\leq$} (M-1-2)
    (M-2-1) -- node{$\leq$} (M-2-2)
;
\path [auto, ->, arrow]
    (M-1-1) edge node{$\id\x \ph$} (M-2-1)
    (M-1-2) edge node{$\id\x \ph$} (M-2-2)
            edge node{$\ev_{\tup a}$} (M-1-3)
    (M-1-3) edge node{$\tilde\ph$} (M-2-3)
    (M-2-2) edge node{$\ev_{\tup a}$} (M-2-3)
;
\end{tikzpicture}
\]
commutes, so $\tilde \ph$ satisfies $\ev_{\tup a}\circ (\id \x \ph) = \tilde\ph \circ \ev_{\tup a}$ as homomorphisms $D\to \ev_{\tup a}(E)$. Hence $\ev_{\tup a}$ is fusion preserving.
\end{proof}

\begin{remark}
As homotopy classes of stable maps the isomorphism of Lemma \ref{lemmaChangeOfRepresentatives} commutes with the equivalence
\[\freeO n B\cF \xrightarrow{\simeq} \coprod_{i} BC_\cF(\tup a_i).\]
This can be seen the following way since we know from Lemma \ref{lemmaEvalFusionPreserving} that the evaluation maps are fusion preserving:

As mentioned in Remark \ref{remarkFreeLoopSpaceDef}, the equivalence $\freeO nB\cF \simeq \coprod_{i} BC_\cF(\tup a_i)$ is adjoint to the maps $f_i\colon B(\ZZ/p^e)^n \x BC_\cF(\tup a_i)  \to B\cF_i$, where $\cF_i$ is the component of $\cF$ containing $\tup a_i$. The stable homotopy class of $f_i$ is represented by the element
\[ [ (\ZZ/p^e)^n\x C_S(\tup a_i), \ev_{\tup a_i}]_{ (\ZZ/p^e)^n\x C_\cF(\tup a_i) }^{\cF_i} \]
in $\AF$.
By construction $\iso_a^b$ has entries
\[\change_{\tup a_i}^{\tup b_i}=[C_S(\tup a_i), \tilde\ph]_{C_\cF(\tup a_i)}^{C_\cF(\tup b_i)},\]
where $\tilde\ph$ is an $\cF$-isomorphism $C_S(\tup a_i)\xrightarrow{\cong} C_S(\tup b_i)$ that sends $\tup a_i$ to $\tup b_i$.

As group homomorphisms we have $\ev_{\tup b_i} \circ (\id\x \tilde\ph) = \tilde\ph \circ \ev_{\tup a_i}$ since both composites send $(t_1,\dotsc,t_n,z)\in (\ZZ/p^e)^n\x C_S(\tup a_i)$ to the same element \[((\tup b_i)_1)^{t_1}\dotsm ((\tup b_i)_n)^{t_n}\cdot \tilde\ph(z) = \tilde\ph\bigl(((\tup a_i)_1)^{t_1}\dotsm ((\tup a_i)_n)^{t_n}\cdot z \bigr).\]
If we compose $\change_{\tup a_i}^{\tup b_i}$ with the evaluation map for $\tup b_i$, we can use that $\ev_{\tup b_i}$ is fusion preserving and apply the special case of Proposition \ref{propFusionDoubleCosetFormula}. This gives us
\begin{align*}
&((\ZZ/p^e)^n\x \change_{\tup a_i}^{\tup b_i}) \cmp [ (\ZZ/p^e)^n\x C_S(\tup b_i), \ev_{\tup b_i}]_{(\ZZ/p^e)^n\x C_\cF(\tup b_i) }^{\cF_i}
\\={}&
[ (\ZZ/p^e)^n\x C_S(\tup a_i), \ev_{\tup b_i}\circ (\id\x \tilde\ph)]_{ (\ZZ/p^e)^n\x C_\cF(\tup a_i)}^{\cF_i}
\\ ={}& [  (\ZZ/p^e)^n\x C_S(\tup a_i), \tilde\ph \circ \ev_{\tup a_i}]_{(\ZZ/p^e)^n\x C_\cF(\tup a_i)}^{\cF_i}
\\ ={}& [(\ZZ/p^e)^n\x C_S(\tup a_i), \ev_{\tup a_i}]_{(\ZZ/p^e)^n\x C_\cF(\tup a_i)}^{\cF_i} \text{ since $\tilde\ph\in \cF_i$.}
\end{align*}
Taking adjoints it follows that $\change_{\tup a_i}^{\tup b_i}$ commutes with the maps from $BC_\cF(\tup a_i)$ and $BC_\cF(\tup b_i)$ to $\freeO n B\cF$ as homotopy classes of stable maps.
\end{remark}

\section{The functor $\freeO n$ for the category of fusion systems}\label{secFusionFreeLoopFunctor}
Given a saturated fusoid $\cF$ over a formal union of $p$-groups $S$, we can apply $\freeO n$ to the characteristic idempotent $\omega_\cF\in\AF(S,S)$. The result is an idempotent endomorphism $\freeO n(\omega_F)$ from $\freeO n S$ to itself.
\begin{definition}
Let $\tel_{\freeO n(\omega_\cF)}$ denote the mapping telescope
\[\tel_{\freeO n(\omega_\cF)} = \colim(\Sinfpp B\freeO n S\xto{\freeO n(\omega_\cF)} \Sinfpp B\freeO n S \xto{\freeO n(\omega_\cF)} \dotsb).\]
Then $\tel_{\freeO n(\omega_\cF)}$ is the retract of $B\freeO n S$ with respect to the idempotent $\freeO n(\omega_\cF)\in \AF(S,S)$.
\end{definition}
Given any $(\cE,\cF)$-stable virtual biset $X\in \AF(\cE,\cF)$, if we apply $\freeO n$ to the relation $X= \omega_\cE\cmp X\cmp \omega_\cF$, we get
\[\freeO n(X) = \freeO n(\omega_\cE) \cmp \freeO n(X) \cmp \freeO n(\omega_\cF).\]
This implies that $\freeO n(X)$ descends to a map $\freeO n(X)\colon \tel_{\freeO n(\omega_\cE)} \to \tel_{\freeO n(\omega_\cF)}$.

The next step for us is to compare $\tel_{\freeO n(\omega_\cF)}$ with the algebraic model for $\freeO n \cF$ in Definition \ref{defFusionAlgLoop}.

\begin{definition}\label{defLnFRetractOfLnS}
Let $\cF$ be a saturated fusoid with underlying union of $p$-groups $S$. We define a matrix $I_\cF\in\AF(\freeO nS,\freeO n \cF)$ as follows:
\[
(I_\cF)_{\tup a,\tup b} = \begin{cases}
0 & \text{if $\tup a$ is not $\cF$-conjugate to $\tup b$,}
\\ \change_{\tup a}^{\tup b} & \text{if $\tup b$ is the representative for the $\cF$-conjugacy class of $\tup a$,}
\end{cases}
\]
where $(I_\cF)_{\tup a,\tup b}$ is an element of $\AF(C_S(\tup a),C_\cF(\tup b))$.

Next we define a matrix $T_\cF\in \AF(\freeO n\cF, \freeO nS)$ by the formula
\[
(T_\cF)_{\tup a,\tup b} = \lc{\tup a}(\omega_\cF)^{\tup b}\in \AF(C_\cF(\tup a),C_S(\tup b)).
\]
It follows that $(T_\cF)_{\tup a,\tup b}$ is zero unless $\tup b$ is in the $\cF$-conjugacy class represented by $\tup a$. Furthermore, $(T_\cF)_{\tup a,\tup b} = \freeO n((\omega_\cF)_S^S)_{\tup a,\tup b}$ by \cite[Proposition 3.15]{RSS_Bold1}, so we can obtain $T_\cF$ from $\freeO n((\omega_\cF)_S^S)$ by deleting all rows not belonging to the chosen representatives $\tup a$ for the $\cF$-conjugacy classes of commuting $n$-tuples in $S$.
\end{definition}

\begin{remark}\label{remarkTransferMapWellDefined}
In order for $T_\cF$ to be well-defined, we need $\lc{\tup a}(\omega_\cF)^{\tup b}$ to be left $C_\cF(\tup a)$-stable. This is a consequence of the fact that $\omega_\cF$ is $\cF$-stable and can be seen as follows:

Suppose $P\leq C_S(\tup a)$ and $\ph\in C_\cF(\tup a)(P,C_S(\tup a))$. Let $\gen{\tup a}P$ be the subgroup of $S$ generated by $\tup a$ and $P$. By definition of the centralizer fusion system, we have an extension of $\ph$, $\tilde\ph\colon \gen{\tup a}P \to C_S(\tup a)$, inside $C_\cF(\tup a)$ with the property that $\tilde\ph(\tup a)=\tup a$.
If we restrict $\lc{\tup a}(\omega_\cF)^{\tup b}$ along $\ph$, we find that
\begin{align*}
    [P,\ph]_P^{C_S(\tup a)}\cmp \lc{\tup a}(\omega_\cF)^{\tup b} &= [P,\incl]_P^{\gen{\tup a}P}\cmp [\gen{\tup a}P,\tilde\ph]_{\gen{\tup a}P}^{C_S(\tup a)}\cmp \lc{\tup a}(\omega_\cF)^{\tup b}
\\ &= [P,\incl]_P^{\gen{\tup a}P}\cmp \lc{\tup a}( [\gen{\tup a}P,\tilde\ph]_{\gen{\tup a}P}^{S}\cmp \omega_\cF)^{\tup b}
\\ &= [P,\incl]_P^{\gen{\tup a}P}\cmp \lc{\tup a}( [\gen{\tup a}P,\incl]_{\gen{\tup a}P}^{S}\cmp \omega_\cF)^{\tup b}
\\ &= [P,\incl]_P^{\gen{\tup a}P}\cmp [\gen{\tup a}P,\incl]_{\gen{\tup a}P}^{C_S(\tup a)}\cmp \lc{\tup a}(\omega_\cF)^{\tup b}
\\ &= [P,\incl]_P^{C_S(\tup a)}\cmp \lc{\tup a}(\omega_\cF)^{\tup b}.
\end{align*}
We conclude that $\lc{\tup a}(\omega_\cF)^{\tup b}$ is left $C_\cF(\tup a)$-stable and an element of $\AF(C_\cF(\tup a),C_S(\tup b))$ as required.
\end{remark}

We claim that $I_\cF$ and $T_\cF$ express $\freeO n \cF$ as a retract of $\freeO n S$ and that $I_\cF\cmp T_\cF=\freeO n((\omega_\cF)_S^S)$. The first claim is proved as Proposition \ref{propFreeLoopRetractIdentity} below. The second claim is an easy consequence of $\cF$-stability for $\omega_\cF$, hence we shall prove the second claim first.
\begin{lemma}\label{lemmaIdempotentOfRetract}
The matrices $I_\cF$ and $T_\cF$ satisfy
\[I_\cF\cmp T_\cF = \freeO n((\omega_\cF)_S^S)\in \AF(\freeO n S,\freeO nS).\]
\end{lemma}

\begin{proof}
Recall from \cite[Proposition 3.15]{RSS_Bold1} that $\freeO n((\omega_\cF)_S^S)_{\tup a,\tup b} = \lc{\tup a}(\omega_\cF)^{\tup b}$. Since $\omega_\cF$ is $\cF$-generated (Definition \ref{defFCharacteristic}), we have $\freeO n((\omega_\cF)_S^S)_{\tup a,\tup b}=0$ unless $\tup a$ and $\tup b$ are $\cF$-conjugate.

By construction, we also have $(I_\cF\cmp T_\cF)_{\tup a,\tup b}=0$ unless $\tup a$ and $\tup b$ are $\cF$-conjugate, in which case
\[
(I_\cF\cmp T_\cF)_{\tup a,\tup b} = (I_{\cF})_{\tup a,\tup c} \cmp (T_\cF)_{\tup c,\tup b} = \change_{\tup a}^{\tup c} \cmp \lc{\tup c}(\omega_\cF)^{\tup b},
\]
where $\tup c$ is the chosen representative for the $\cF$-conjugacy class of $\tup a$ and $\tup b$.
By Lemma \ref{lemmaCentralizerInclusion}, $\change_{\tup a}^{\tup c}=[C_S(\tup a),\ph]_{C_S(\tup a)}^{C_\cF(\tup c)}$, for any $\ph\colon C_S(\tup a)\to C_S(\tup c)$ in $\cF$ such that $\ph(\tup a)=\tup c$. Precomposing with $[C_S(\tup a),\ph]$ is the same as the restriction $\Res_\ph$ along $\ph$. It follows that
\[
\change_{\tup a}^{\tup c} \cmp \lc{\tup c}(\omega_\cF)^{\tup b} = \lc{\tup a}(\Res_\ph \omega_\cF)^{\tup b} = \lc{\tup a}(\Res_{C_S(\tup a)}^S \omega_\cF)^{\tup b} = \lc{\tup a}(\omega_\cF)^{\tup b},
\]
since the restrictions $\Res_\ph \omega_\cF$ and $\Res_{C_S(\tup a)}^S \omega_\cF$ are equal in $\AF(C_S(\tup a),S)$ by left $\cF$-stability of $\omega_\cF$.
We conclude that $(I_\cF\cmp T_\cF)_{\tup a,\tup b}= \freeO n((\omega_\cF)_S^S)_{\tup a,\tup b}$ for all commuting $n$-tuples $\tup a$ and $\tup b$.
\end{proof}

\begin{prop}\label{propRetractOfSimpleLoopFunctor}
Let $\cE$ and $\cF$ be saturated fusoids with underlying unions of $p$-groups $R$ and $S$ respectively. Suppose $X\in \AF(\cE,\cF)$, we can apply $\freeO n$ to $X_R^S\in \AF(R,S)$ and get a matrix of virtual bisets $\freeO n (X_R^S)\in \AF(\freeO n R,\freeO n S)$. Precomposing with $T_\cE$ and postcomposing with $I_\cF$ then gives us a matrix in $\AF(\freeO n \cE,\freeO n\cF)$ with entries
\[
(T_\cE\cmp \freeO n(X_R^S) \cmp I_\cF)_{\tup a,\tup b} =\sum_{\substack{[\tup b']\in \cntuples nS \\  \tup b' \sim_\cF \tup b}} (\freeO n(X_R^S))_{\tup a,\tup b'} \cmp \change_{\tup b'}^{\tup b} = \sum_{\substack{[\tup b']\in \cntuples nS \\ \tup b' \sim_\cF \tup b}} \lc{\tup a}X^{\tup b'}\cmp \change_{\tup b'}^{\tup b}.
\]
\end{prop}

\begin{proof}
We first calculate the entries of $T_\cE\cmp \freeO n(X_R^S)$. Suppose $\tup a$ represents an $\cF$-conjugacy class and $\tup c$ represents an $S$-conjugacy class of $n$-tuples in $S$.
We have
\[
(T_\cE\cmp \freeO n(X_R^S))_{\tup a,\tup c} = \sum_{[\tup d]\in\cntuples n S} (T_\cE)_{\tup a,\tup d}\cmp \freeO n(X_R^S)_{\tup d,\tup c}.
\]
By Definition \ref{defLnFRetractOfLnS}, $(T_\cE)_{\tup a,\tup d} = \freeO n((\omega_\cE)_R^R)_{\tup a,\tup d}$. We plug this in above to get
\begin{multline*}
(T_\cE\cmp \freeO n(X_R^S))_{\tup a,\tup c} = \sum_{[\tup d]\in\cntuples n S} \freeO n((\omega_\cE)_R^R)_{\tup a,\tup d}\cmp \freeO n(X_R^S)_{\tup d,\tup c} 
\\= \freeO n((\omega_\cE\cmp X)_R^S)_{\tup a,\tup c} = (\freeO n(X_R^S))_{\tup a,\tup c},
\end{multline*}
making use of the fact that $X$ is left $\cE$-stable. From here we can easily calculate $(T_\cE\cmp \freeO n(X_R^S) \cmp I_\cF)_{\tup a,\tup b}$, if we recall that $(\freeO n(X_R^S))_{\tup a,\tup c} = \lc{\tup a} X^{\tup c}$ by \cite[Proposition 3.15]{RSS_Bold1}. We have
\begin{multline*}
(T_\cE\cmp \freeO n(X_R^S) \cmp I_\cF)_{\tup a,\tup b} = \sum_{[\tup d]\in\cntuples n S} (\freeO n(X_R^S))_{\tup a,\tup c} \cmp (I_\cF)_{\tup c,\tup b} 
\\ =\sum_{\substack{[\tup b']\in \cntuples nS \\  \tup b' \sim_\cF \tup b}} (\freeO n(X_R^S))_{\tup a,\tup b'} \cmp \change_{\tup b'}^{\tup b} = \sum_{\substack{[\tup b']\in \cntuples nS \\ \tup b' \sim_\cF \tup b}} \lc{\tup a}X^{\tup b'} \cmp \change_{\tup b'}^{\tup b}.\qedhere
\end{multline*}
\end{proof}

\begin{prop}\label{propFreeLoopRetractIdentity}
Let $\cF$ be a saturated fusoid. Then $T_\cF\cmp I_\cF$ is the identity in $\AF(\freeO n\cF,\freeO n\cF)$, i.e. the diagonal matrix with diagonal entries $\omega_{C_\cF(\tup a)}$ for $[\tup a]\in \cntuples n\cF$. 
\end{prop}

\begin{proof}
First note that
\[T_\cF = T_\cF \cmp \freeO n((\omega_\cF)_S^S).\]
As in the proof of Proposition \ref{propRetractOfSimpleLoopFunctor}, this is an easy consequence of Definition \ref{defLnFRetractOfLnS}:
\begin{multline*}
\bigl(T_\cF\cmp \freeO n((\omega_\cF)_S^S)\bigr)_{\tup a,\tup b} =  \sum_{[c]\in\cntuples n S} \freeO n((\omega_\cF)_S^S)_{\tup a,\tup c}\cmp \freeO n((\omega_\cF)_S^S)_{\tup c,\tup b}
\\= \freeO n((\omega_\cF\cmp \omega_\cF)_S^S)_{\tup a,\tup b} 
= \freeO n((\omega_\cF)_S^S)_{\tup a,\tup b} = (T_\cF)_{\tup a,\tup b}.
\end{multline*}
We now apply Proposition \ref{propRetractOfSimpleLoopFunctor} with $X=\omega_\cF\in\AF(S,S)$:
\[(T_\cF \cmp I_\cF)_{\tup a,\tup b} = (T_\cF \cmp \freeO n((\omega_\cF)_S^S)\cmp I_\cF)_{\tup a,\tup b} = \sum_{\substack{[\tup b']\in \cntuples nS \\ \tup b' \sim_\cF \tup b}} \lc{\tup a}(\omega_\cF)^{\tup b'}\cmp \change_{\tup b'}^{\tup b}.
\]
Since $\omega_\cF$ is $\cF$-generated, $\lc{\tup a}(\omega_\cF)^{\tup b'} = 0$ unless $\tup a$ is $\cF$-conjugate to $\tup b'$. Hence $\tup a$ is $\cF$-conjugate to $\tup b$, and since both $\tup a$ and $\tup b$ are representatives for their $\cF$-conjugacy class, we conclude $\tup a=\tup b$ unless $\lc{\tup a}(\omega_\cF)^{\tup b'}=0$.

Thus $T_\cF \cmp I_\cF\in \AF(\freeO n\cF,\freeO n\cF)$ is a diagonal matrix. It remains to show that each diagonal entry $(T_\cF \cmp I_\cF)_{\tup a,\tup a}$ is the characteristic idempotent for $C_\cF(\tup a)$.
We will prove that $(T_\cF \cmp I_\cF)_{\tup a,\tup a}$ has all of the properties of Definition \ref{defFCharacteristic} and is idempotent.

A direct application of Lemma \ref{lemmaIdempotentOfRetract} gives
\[T_\cF \cmp I_\cF \cmp T_\cF \cmp I_\cF = T_\cF \cmp \freeO n((\omega_\cF)_S^S)\cmp I_\cF = T_\cF \cmp I_\cF\]
so $(T_\cF \cmp I_\cF)_{\tup a,\tup a}$ is idempotent.

To see that $(T_\cF \cmp I_\cF)_{\tup a,\tup a}$ is $C_\cF(\tup a)$-generated consider the orbit decomposition of $\lc{\tup a}(\omega_\cF)^{\tup a'}$ as a virtual $(C_S(\tup a),C_S(\tup a'))$-biset for each $[\tup a']\in\cntuples n S$ with $\tup a'\sim_\cF \tup a$. For an orbit $[P,\ph]$, with $P\leq C_S(\tup a)$ and $\ph\colon P\to C_S(\tup a')$, to be summand of $\lc{\tup a}(\omega_\cF)^{\tup a'}$, we must have $\tup a\in P$ and $\ph(\tup a)=\tup a'$.
Suppose $\change_{\tup a'}^{\tup a} = [C_S(\tup a'),\rho]_{C_S(\tup a')}^{C_S(\tup a)}\cmp \omega_{C_\cF(\tup a)}$, where $\rho\colon C_S(\tup a')\to C_S(\tup a)$ is such that $\rho(\tup a')=\tup a$. Then we have
\[
[P,\ph]_{C_S(\tup a)}^{C_S(\tup a')}\cmp \change_{\tup a'}^{\tup a} = \Bigl( [P,\ph]_{C_S(\tup a)}^{C_S(\tup a')} \cmp [C_S(\tup a'),\rho]_{C_S(\tup a')}^{C_S(\tup a)} \Bigr) \cmp \omega_{C_\cF(\tup a)} = [P,\rho\circ \ph]_{C_S(\tup a)}^{C_S(\tup a)} \cmp \omega_{C_\cF(\tup a)}.
\]
Now $\rho(\ph(\tup a))=\tup a$, so $\rho\circ \ph$ is a morphism in $C_\cF(\tup a)$. Hence $[P,\rho\circ \ph]$ is $C_\cF(\tup a)$-generated, and $\omega_{C_\cF(\tup a)}$ is $C_\cF(\tup a)$-generated by definition of being the characteristic idempotent. Consequently the composite $[P,\rho\circ \ph] \cmp \omega_{C_\cF(\tup a)}$ is $C_\cF(\tup a)$-generated as well. The diagonal entry $(T_\cF \cmp I_\cF)_{\tup a,\tup a}$ is thus a linear combination of $C_\cF(\tup a)$-generated elements and therefore $C_\cF(\tup a)$-generated.

Next, $(T_\cF \cmp I_\cF)_{\tup a,\tup a}$ is right $C_\cF(\tup a)$-stable because each $\change_{\tup a'}^{\tup a}\in\AF(C_S(\tup a'),C_\cF(\tup a))$ is right $C_\cF(\tup a)$-stable. Similarly, $(T_\cF \cmp I_\cF)_{\tup a,\tup a}$ is left $C_\cF(\tup a)$-stable because each $\lc{\tup a}(\omega_\cF)^{\tup a'}$ is left $C_\cF(\tup a)$-stable by Remark \ref{remarkTransferMapWellDefined}.

Finally, we need to show that $\abs{(T_\cF \cmp I_\cF)_{\tup a,\tup a}}/\abs{C_S(\tup a)}$ is invertible in $\Z_p$, i.e. is not divisible by $p$.
According to \cite{ReehIdempotent}*{Theorem B}, we have
\[\abs{\lc{\tup a}(\omega_\cF)^{\tup a'}} = \frac{\abs S}{\abs{\cF(\gen{\tup a},S)}}\in \Z_p,\]
and $\abs{\cF(\gen{\tup a},S)}$ is simply the total number of $n$-tuples that are $\cF$-conjugate to $\tup a$.

Since $\omega_{C_\cF(\tup a)}$ is idempotent, $\abs{\omega_{C_\cF(\tup a)}}/\abs{C_S(\tup a)}=1$ (alternatively just apply \cite{ReehIdempotent}*{Theorem B} to the fixed points for $\omega_{C_\cF(\tup a)}$ with respect to  the trivial subgroup of $C_S(\tup a)$).
Hence we see that
\[
\frac{\abs{\change_{\tup a'}^{\tup a}}}{\abs{C_S(\tup a)}} = \frac{\abs{[C_S(\tup a'),\rho]_{C_S(\tup a')}^{C_S(\tup a)}}}{\abs{C_S(\tup a)}}\cdot \frac{\abs{\omega_{C_\cF(\tup a)}}}{\abs{C_S(\tup a)}} 
 = 1 \cdot 1 = 1.
\]
Putting the pieces together, we find that 
\begin{align*}
    \frac{\abs{(T_\cF \cmp I_\cF)_{\tup a,\tup a}}}{\abs{C_S(\tup a)}} &= \sum_{\substack{[\tup a']\in \cntuples nS \\ \tup a' \sim_\cF \tup a}} \frac{\abs{\lc{\tup a}(\omega_\cF)^{\tup a'}}}{\abs{C_S(\tup a')}}\cdot\frac{\abs{ \change_{\tup a'}^{\tup a}}}{\abs{C_S(\tup a)}}
\\ &= \sum_{\substack{[\tup a']\in \cntuples nS \\ \tup a' \sim_\cF \tup a}} \frac{\abs{S}}{\abs{C_S(\tup a')}\cdot\abs{\cF(\gen{\tup a},S)}}\cdot 1
\\ &= \frac 1{\abs{\cF(\gen{\tup a},S)}} \sum_{\substack{[\tup a']\in \cntuples nS \\ \tup a' \sim_\cF \tup a}} \frac{\abs{S}}{\abs{C_S(\tup a')}}
\\ &= \frac 1{\abs{\cF(\gen{\tup a},S)}} \sum_{\substack{[\tup a']\in \cntuples nS \\ \tup a' \sim_\cF \tup a}} \abs{[\tup a']}
\\ &= \frac 1{\abs{\cF(\gen{\tup a},S)}} \cdot \abs{\{\text{$n$-tuples $\tup a'\in S$}\mid \tup a'\sim_{\cF}\tup a\}}
\\ &= 1.
\end{align*}
This completes the proof that $(T_\cF \cmp I_\cF)_{\tup a,\tup a}$ is in fact $C_\cF(\tup a)$-characteristic, hence by uniqueness of characteristic idempotents $(T_\cF \cmp I_\cF)_{\tup a,\tup a}=\omega_{C_\cF(\tup a)}$ as required.
\end{proof}

\begin{cor}\label{corFreeLoopTelescopeEquiv}
Let $\cF$ be a saturated fusoid over a union of $p$-groups $S$. Then $I_\cF$ and $T_\cF$ induce inverse equivalences 
\[\tel_{\freeO n((\omega_\cF)_S^S)} \simeq \Sinfpp B\freeO n \cF\]
in $\Ho(\Spp)$.
\end{cor}

\begin{proof}
By Lemma \ref{lemmaIdempotentOfRetract} and Proposition \ref{propFreeLoopRetractIdentity} the matrices $I_\cF$ and $T_\cF$ induce maps between the towers
\[\freeO n S \xto{\freeO n(\omega_\cF)} \freeO n S  \xto{\freeO n(\omega_\cF)} \dotsb\]
and
\[\freeO n\cF \xto{\id} \freeO n\cF\xto{\id} \dotsb.\]
The composite $T_\cF\cmp I_\cF$ is simply the identity on the constant tower $\freeO n \cF$. The composite $I_\cF\cmp T_\cF$ applies $\freeO n (\omega_\cF)$ levelwise to the tower $\freeO n S \xto{\freeO n(\omega_\cF)} \freeO n S  \xto{\freeO n(\omega_\cF)} \dotsb$ which is homotopic to the identity on the colimit $\tel_{\freeO n((\omega_\cF)_S^S)}$.
\end{proof}

\begin{prop}\label{propSimpleFusionLoop}
The functor $\freeO n$ on unions of $p$-groups extends to a functor $\freeO n\colon \AF\to \AF$ given on objects by $\cF\mapsto \freeO n \cF$ and on morphisms $X\in \AF(\cE,\cF)$ by the matrix with entries
\[\freeO n(X_\cE^\cF)_{\tup a,\tup b} =\sum_{\substack{[\tup b']\in \cntuples nS \\  \tup b' \sim_\cF \tup b}} (\freeO n(X_R^S))_{\tup a,\tup b'} \cmp \change_{\tup b'}^{\tup b} = \sum_{\substack{[\tup b']\in \cntuples nS \\ \tup b' \sim_\cF \tup b}} \lc{\tup a}X^{\tup b'}\cmp \change_{\tup b'}^{\tup b} \qquad\in \AF(C_\cE(\tup a),C_\cF(\tup b)),\]
where $R$ and $S$ are the underlying unions of $p$-groups for $\cE$ and $\cF$ respectively.
\end{prop}

\begin{remark} \label{rem:idempotents}
By Proposition \ref{propRetractOfSimpleLoopFunctor}, we then have $\freeO n(X_\cE^\cF) = T_\cE\cmp \freeO n(X_R^S) \cmp I_\cF$. In particular, we have $I_\cF=\freeO n((\omega_\cF)_S^\cF)$ and $T_\cF=\freeO n((\omega_\cF)_\cF^S)$ for any saturated fusoid $\cF$ over $S$.
\end{remark}

\begin{proof}
Proposition \ref{propFreeLoopRetractIdentity} states that $\freeO n$ takes the identity $\omega_\cF$ on $\cF$ to the identity on $\freeO n\cF$.

Let $\cE$, $\cF$, and $\cG$ be saturated fusoids over $R$, $S$, and $T$, respectively. Suppose $X\in \AF(\cE,\cF)$ and $Y\in \AF(\cF,\cG)$. As in the Remark \ref{rem:idempotents}, we have 
\[\freeO n(X_\cE^\cF) = T_\cE\cmp \freeO n(X_R^S) \cmp I_\cF.\]
It follows that
\begin{align*}
    \freeO n(X_\cE^\cF)\cmp \freeO n(Y_\cF^\cG) &= T_\cE\cmp \freeO n(X_R^S)\cmp I_\cF\cmp T_\cF \cmp \freeO n(Y_S^T) \cmp I_\cG
    \\ &= T_\cE\cmp \freeO n(X_R^S)\cmp \freeO n((\omega_\cF)_S^S) \cmp \freeO n(Y_S^T) \cmp I_\cG
    \\ &= T_\cE\cmp \freeO n(X_R^S\cmp (\omega_\cF)_S^S \cmp Y_S^T) \cmp I_\cG
    \\ &= T_\cE\cmp \freeO n((X\cmp Y)_S^T) \cmp I_\cG
    \\ &= \freeO n((X\cmp Y)_\cE^\cG).
\end{align*}
The characteristic idempotent $\omega_\cF$ disappears from the middle because both $X$ and $Y$ are $\cF$-stable (in fact either of these would be enough). Thus $\freeO n$ preserves composition of virtual bisets between fusoids.
\end{proof}

Similarly to \cite[Corollary 3.18]{RSS_Bold1}, there is also a formula for $\freeO n(X_\cE^\cF)_{\tup a,\tup b}$ in terms of the restriction of $X\in \AF(\cE,\cF)$ to the centralizer $C_\cE(\tup a)$. We state the formula for fusion systems and the formula can easily be applied component-wise in the general case.
\begin{cor}\label{corSimpleFusionLoopOrbits}
Let $\cE$ and $\cF$ be saturated fusion systems over $p$-groups $R$ and $S$ respectively, and suppose $X\in \AF(\cE,\cF)$ is a virtual biset. Furthermore, let $\tup a$ in $\cE$ and  $\tup b$ in $\cF$ be chosen representatives for conjugacy classes of commuting $n$-tuples (according to Convention \ref{conventionTupleReps}). Consider the restriction of $X$ to the centralizer fusion system $C_\cE(\tup a)$, and write $M_{C_\cE(\tup a)}^\cF$ as a linear combination of basis elements (recalling Convention \ref{conventionFusionOrbitDecomposition}):
\[
X_{C_\cE(\tup a)}^\cF =\sum_{(P,\ph)}c_{P,\ph} \cdot [P,\ph]_{C_\cE(\tup a)}^\cF,
\]
where $P\leq C_R(\tup a)$ and $\ph\colon P\to S$.

The matrix entry $\freeO n(X)_{\tup a,\tup b}$ then satisfies the formula
\[
\freeO n(X)_{\tup a,\tup b} = \sum_{\substack{(P,\ph)\text{ s.t. $\tup a\in P$ and}\\\text{$\ph(\tup a)$ is $\cF$-conjugate to $\tup b$}}}c_{P,\ph}\cdot [P,\ph]_{C_\cE(\tup a)}^{C_S(\ph(\tup a))}\cmp \change_{\ph(\tup a)}^{\tup b},
\]
with $P$, $\ph$, and $c_{P,\ph}$ as in the linear combination above.
\end{cor}

\begin{proof}
Instead of writing $X_{C_\cE(\tup a)}^\cF$ as a linear combination of basis elements in $\AF(C_\cE(\tup a),\cF)$, consider $X$ as an element of $\AF(C_R(\tup a),S)$:
\[
X_{C_R(\tup a)}^S = \sum_{(P,\ph)}u_{P,\ph} \cdot [P,\ph]_{C_R(\tup a)}^S,
\]
with a (possibly) different collection of coefficients $u_{P,\ph}$.
If we precompose with $\omega_{C_\cE(\tup a)}$ and postcompose with $\omega_\cF$, the linear combination above becomes
\[
X_{C_\cE(\tup a)}^\cF = \omega_{C_\cE(\tup a)}\cmp X_{C_R(\tup a)}^S \cmp \omega_{\cF} = \sum_{(P,\ph)}u_{P,\ph} \cdot [P,\ph]_{C_\cE(\tup a)}^\cF.
\]
Thus the coefficients $u_{P,\ph}$ and $c_{P,\ph}$ provide different choices for the decomposition in $\AF(C_\cE(\tup a),\cF)$, but only $u_{P,\ph}$ provides a decomposition in $\AF(C_R(\tup a),S)$. We shall start by proving that the formula for $\freeO n(X)_{\tup a,\tup b}$ in the statement of the corollary is independent of the choice of linear combination for $X_{C_\cE(\tup a)}^\cF$.

Two basis elements $[P,\ph]_{C_\cE(\tup a)}^\cF$ and $[Q,\psi]_{C_\cE(\tup a)}^\cF$ are equal if and only if $Q$ is $C_\cE(\tup a)$-isomorphic to $P$ and $\psi$ arises from $\ph$ by precomposing with a map in $C_\cE(\tup a)$ and postcomposing with a map in $\cF$.
The total coefficient in front of a particular basis element $[P,\ph]_{C_\cE(\tup a)}^\cF$ in the linear combinations for $X_{C_\cE(\tup a)}^\cF$ is therefore given by the two sums
\[
\sum_{\text{$(P',\ph')$, $(C_\cE(\tup a),\cF)$-conj. to $(P,\ph)$}} c_{P',\ph'} \qquad\text{and}\qquad \sum_{\text{$(P',\ph')$, $(C_\cE(\tup a),\cF)$-conj. to $(P,\ph)$}} u_{P',\ph'}.
\]
Hence these two coefficient sums must be equal. 

Suppose $(P',\ph')$ is $(C_\cE(\tup a),\cF)$-conjugate to a particular $(P,\ph)$, and suppose further that $\tup a\in P$ and $\ph(\tup a)$ is $\cF$-conjugate to $\tup b$. Since $P'$ is isomorphic to $P$ in $C_\cE(\tup a)$, we also have $\tup a\in P'$. Let $\ph'=\gamma\circ \ph\circ \alpha$ with $\alpha\in C_\cE(\tup a)$ and $\gamma\in \cF$, then $\alpha(\tup a)=\tup a$ and $\ph'(\tup a)=\gamma(\ph(\tup a))$ is $\cF$-conjugate to $\ph(\tup a)$ and hence to $\tup b$.

Let $\change_{\ph'(\tup a)}^{\tup b}=[C_S(\ph'(\tup a)),\rho]_{C_S(\tup a)}^{C_\cF(\tup b)}$, then $\rho\circ\gamma\colon \ph(C_R(\tup a)) \to C_S(\tup b)$ extends to a morphism $\tilde{\rho\circ\gamma}\colon C_S(\ph(\tup a))\to C_S(\tup b)$ and $[C_S(\ph(\tup a)), \tilde{\rho\circ\gamma}]_{C_S(\tup a)}^{C_\cF(\tup b)}=\change_{\ph(\tup a)}^{\tup b}$. We next have
\begin{multline*}
[P,\ph']_{C_\cE(\tup a)}^{C_S(\ph'(\tup a)}\cmp \change_{\ph'(\tup a)}^{\tup b} = [P,\gamma\circ \ph\circ \alpha]_{C_\cE(\tup a)}^{C_S(\ph'(\tup a)}\cmp \change_{\ph'(\tup a)}^{\tup b} \\= [P,\ph]_{C_\cE(\tup a)}^{C_S(\ph'(\tup a)}\cmp [C_S(\ph(\tup a)), \tilde{\rho\circ\gamma}]_{C_S(\tup a)}^{C_\cF(\tup b)} = [P,\ph]_{C_\cE(\tup a)}^{C_S(\ph'(\tup a)}\cmp \change_{\ph(\tup a)}^{\tup b}.
\end{multline*}
The pairs $(P',\ph')$ and $(P,\ph)$ therefore give the same contribution to the formula in the statement of the corollary whenever the pairs are $(C_\cE(\tup a),\cF)$-conjugate. We conclude that
\begin{multline}\label{eqChangeOfCoefficientsForLinearComb}
    \sum_{\substack{(P,\ph)\text{ s.t. $\tup a\in P$ and}\\\text{$\ph(\tup a)$ is $\cF$-conjugate to $\tup b$}}}c_{P,\ph}\cdot [P,\ph]_{C_\cE(\tup a)}^{C_S(\ph(\tup a))}\cmp \change_{\ph(\tup a)}^{\tup b} \\= \sum_{\substack{(P,\ph)\text{ s.t. $\tup a\in P$ and}\\\text{$\ph(\tup a)$ is $\cF$-conjugate to $\tup b$}}}u_{P,\ph}\cdot [P,\ph]_{C_\cE(\tup a)}^{C_S(\ph(\tup a))}\cmp \change_{\ph(\tup a)}^{\tup b}.
\end{multline}
Hence it is sufficient to prove the corollary with the coefficients $u_{P,\ph}$.

By Proposition \ref{propSimpleFusionLoop}, we can write $\freeO n(X_\cE^\cF)_{\tup a,\tup b}$ as
\begin{align*}
    \freeO n(X_\cE^\cF)_{\tup a,\tup b} &=  \sum_{\substack{[\tup b']\in \cntuples n S\\\tup b'\sim_\cF \tup b}} \freeO n(X_R^S)_{\tup a,\tup b'}\cmp \change_{\tup b'}^{\tup b}.
\end{align*}
We can then apply \cite[Corollary 3.18]{RSS_Bold1} with the linear combination for $X_{C_R(\tup a)}^S$ given by the coefficients $u_{P,\ph}$:
\begin{align*}
    \freeO n(X_\cE^\cF)_{\tup a,\tup b} &= \sum_{\substack{[\tup b']\in \cntuples n S\\\tup b'\sim_\cF \tup b}} \freeO n(X_R^S)_{\tup a,\tup b'}\cmp \change_{\tup b'}^{\tup b}
    \\ &= \sum_{\substack{[\tup b']\in \cntuples n S\\\tup b'\sim_\cF \tup b}}   \sum_{\substack{(P,\ph)\text{ s.t. $\tup a\in P$ and}\\\text{$\ph(\tup a)$ is $S$-conjugate to $\tup b'$}}}u_{P,\ph}\cdot [P,\ph]_{C_R(\tup a)}^{C_S(\ph(\tup a))}\cmp \change_{\ph(\tup a)}^{\tup b'}\cmp \change_{\tup b'}^{\tup b}
    \\&= \sum_{\substack{(P,\ph)\text{ s.t. $\tup a\in P$ and}\\\text{$\ph(\tup a)$ is $\cF$-conjugate to $\tup b$}}}u_{P,\ph}\cdot [P,\ph]_{C_R(\tup a)}^{C_S(\ph(\tup a))}\cmp \change_{\ph(\tup a)}^{\tup b}.
\end{align*}
Finally, note that $\freeO n(X_\cE^\cF)_{\tup a,\tup b}\in \AF(C_\cE(\tup a),C_\cF(\tup b))$ is left $C_\cE(\tup a)$-stable and as such does not change if we precompose with $\omega_{C_\cE(\tup a)}$:
\begin{align*}
    \freeO n(X_\cE^\cF)_{\tup a,\tup b} &= \omega_{C_\cE(\tup a)}\cmp \freeO n(X_\cE^\cF)_{\tup a,\tup b} 
    \\&= \sum_{\substack{(P,\ph)\text{ s.t. $\tup a\in P$ and}\\\text{$\ph(\tup a)$ is $\cF$-conjugate to $\tup b$}}}u_{P,\ph}\cdot \omega_{C_\cE(\tup a)}\cmp [P,\ph]_{C_R(\tup a)}^{C_S(\ph(\tup a))}\cmp \change_{\ph(\tup a)}^{\tup b}
    \\&= \sum_{\substack{(P,\ph)\text{ s.t. $\tup a\in P$ and}\\\text{$\ph(\tup a)$ is $\cF$-conjugate to $\tup b$}}}u_{P,\ph}\cdot [P,\ph]_{C_\cE(\tup a)}^{C_S(\ph(\tup a))}\cmp \change_{\ph(\tup a)}^{\tup b}.
\end{align*}
Finally, by \eqref{eqChangeOfCoefficientsForLinearComb} we can replace the coefficients $u_{P,\ph}$ with the coefficients $c_{P,\ph}$ to get the formula in the corollary:
\[\freeO n(X_\cE^\cF)_{\tup a,\tup b}=\sum_{\substack{(P,\ph)\text{ s.t. $\tup a\in P$ and}\\\text{$\ph(\tup a)$ is $\cF$-conjugate to $\tup b$}}}c_{P,\ph}\cdot [P,\ph]_{C_\cE(\tup a)}^{C_S(\ph(\tup a))}\cmp \change_{\ph(\tup a)}^{\tup b}.\qedhere\]

\end{proof}

\section{Evaluation maps and $\twistO n$ for fusion systems}\label{secFusionEvaluation}
Given a saturated fusoid $\cF$ over a union of $p$-groups $S$, we can apply $\twistO n$ to the characteristic idempotent $\omega_\cF\in \AF(S,S)$ to get an idempotent endomorphism of $(\Z/p^e)^n\x \freeO n S$ in $\AF$.
As in Section \ref{secFusionFreeLoopFunctor}, we let $\tel_{\twistO n(\omega_\cF)}$ denote the mapping telescope 
\[\tel_{\twistO n(\omega_\cF)} = \colim(\Sinfpp  B((\Z/p^e)^n\x \freeO n S)\xto{\twistO n(\omega_\cF)} \Sinfpp  B((\Z/p^e)^n\x \freeO n S) \xto{\twistO n(\omega_\cF)} \dotsb).\]
This is a retract of $\Sinfpp B((\Z/p^e)^n\x \freeO n S)$.

Our next step is to prove that $\tel_{\twistO n(\omega_\cF)}$ is equivalent to $\Sinfpp B((\Z/p^e)^n)\wedge \tel_{\freeO n(\omega_\cF)}$ and hence to $\Sinfpp B((\Z/p^e)^n\x \freeO n\cF)$ via Corollary \ref{corFreeLoopTelescopeEquiv}.
In order to produce this equivalence, we first need a technical lemma about the characteristic idempotent $\omega_\cF$, followed by a proposition relating the idempotents $\twistO n(\omega_\cF)$ and $\freeO n(\omega_\cF)$.

\begin{lemma}\label{lemmaIdempotentCoefficients}
Let $\cF$ be a saturated fusion system over $S$, and let $\tup a$ be a commuting $n$-tuple in $S$. Write the characteristic idempotent $\omega_\cF$, restricted to $C_S(\tup a)$, as a linear combination according to \cite[Convention 3.2]{RSS_Bold1}:
\[(\omega_\cF)_{C_S(\tup a)}^S = \sum_{\substack{R\leq C_S(\tup a)\\ \ph\in\cF( R, S)}} c_{R,\ph}\cdot [R,\ph]_{C_S(\tup a)}^S.\]
Then the coefficients $c_{R,\ph}$ satisfy the following relation for each subgroup $R\leq C_S(\tup a)$
\[
\sum_{\substack{R'\sim_{C_S(\tup a)} R\\ \ph\in \cF(R',S)}} c_{R',\ph} = \begin{cases}1&\text{if $R=C_S(\tup a)$,}\\ 0&\text{otherwise.}\end{cases}
\]
\end{lemma}

\begin{proof}
The restriction $(\omega_\cF)_{C_S(\tup a)}^S$ equals $[C_S(\tup a),\incl]_{C_S(\tup a)}^S\cmp \omega_\cF$.
According to \cite{ReehIdempotent}*{Corollary 5.13} the map $\pi\colon \AF(C_S(\tup a),S)\to \AF(C_S(\tup a),\cF)$ given by $X\mapsto X\cmp \omega_\cF$ coincides with the map of Burnside rings $\pi'\colon A(C_S(\tup a)\x S)^\wedge_p \to A(C_S(\tup a)\x \cF)^\wedge_p$ given in \cite{ReehIdempotent}*{Theorem A} for the product fusion system $C_S(\tup a)\x \cF$. We are particularly interested in 
\[(\omega_\cF)_{C_S(\tup a)}^S = [C_S(\tup a),\incl]_{C_S(\tup a)}^S\cmp \omega_\cF = \pi([C_S(\tup a),\incl]_{C_S(\tup a)}^S).\]

The result is now a consequence of \cite{ReehIdempotent}*{Remark 4.7}, which describes how the coefficients of an element $X\in A(C_S(\tup a)\x S)^\wedge_p$ relates to the coefficients of $\pi'(X)\in A(C_S(\tup a)\x \cF)^\wedge_p$. Since $\pi'$ for the product fusion system $C_S(\tup a)\x \cF$ coincides with $\pi$ for bisets, we have the same relation between coefficients of $X\in \AF(C_S(\tup a),S)$ and $\pi(X)\in \AF(C_S(\tup a),\cF)$. 

Let $c_{R',\ph}(X)$ denote the coefficient of the orbit $[R',\ph]_{C_S(\tup a)}^S$ in the decomposition of a general element $X\in \AF(C_S(\tup a),S)$.
\cite{ReehIdempotent}*{Remark 4.7} states that, since $(\omega_\cF)_{C_S(\tup a)}^S= \pi([C_S(\tup a),\incl]_{C_S(\tup a)}^S)$ is the image of a transitive biset, we have
\begin{multline*}
\sum_{\substack{R'\sim_{C_S(\tup a)} R\\ \ph\in \cF(R',S)}} c_{R',\ph}(\pi([C_S(\tup a),\incl]_{C_S(\tup a)}^S)) \\= \sum_{\substack{R'\sim_{C_S(\tup a)} R\\ \ph\in \cF(R',S)}} c_{R',\ph}([C_S(\tup a),\incl]_{C_S(\tup a)}^S) = \begin{cases}1&\text{if $R=C_S(\tup a)$,}\\ 0&\text{otherwise.}\end{cases}\qedhere
\end{multline*}
\end{proof}

\begin{prop}\label{propComposingLoopIdempotents}
Let $\cF$ be a saturated fusion system over $S$, and let $n\geq 0$. The characteristic idempotent gives us two idempotent endomorphisms of $(\Z/p^e)^n\x \freeO nS$ coming from the functors $(\Z/p^e)^n\x \freeO n(-)$ and $\twistO n (-)$. The two resulting idempotents satisfy the following relations when composed with each other:
\begin{enumerate}
\renewcommand{\theenumi}{$(\roman{enumi})$}\renewcommand{\labelenumi}{\theenumi}
    \item\label{itemIdempotentCompTwist} $\displaystyle \bigl( (\Z/p^e)^n\x \freeO n(\omega_\cF)\bigr)\cmp \twistO n(\omega_\cF) = \twistO n(\omega_\cF)$
    \item\label{itemIdempotentCompFree} $\displaystyle  \twistO n(\omega_\cF)\cmp \bigl((\Z/p^e)^n\x \freeO n(\omega_\cF) \bigr)= (\Z/p^e)^n\x \freeO n(\omega_\cF)$.
\end{enumerate}
Both of these composites are taken in $\AF((\Z/p^e)^n\x \freeO nS, (\Z/p^e)^n\x \freeO nS)$.
\end{prop}

\begin{proof} 
Let $\tup a$ and $\tup b$ be representatives for the conjugacy classes of commuting $n$-tuples in $S$. Restrict $\omega_\cF$ to $C_S(\tup a)$ on the left and write
\[(\omega_\cF)_{C_S(\tup a)}^S = \sum_{\substack{R\leq C_S(\tup a)\\ \ph\in\cF( R, S)}} c_{R,\ph}\cdot [R,\ph]_{C_S(\tup a)}^S\]
as in Lemma \ref{lemmaIdempotentCoefficients}.

Let us first consider part \ref{itemIdempotentCompTwist} of the proposition. According to \cite[Corollary 3.18]{RSS_Bold1} the matrix $\freeO n(\omega_\cF)\in \AF(\freeO n S,\freeO nS)$ has entries
\[
\freeO n(M)_{\tup a,\tup c} = \sum_{\substack{R\leq C_S(\tup a)\\\ph\in\cF(R,S)\\\text{ s.t. $\tup a\in R$ and}\\\text{$\ph(\tup a)$ is $S$-conjugate to $\tup c$}}}c_{R,\ph}\cdot [R,\ph]_{C_S(\tup a)}^{C_S(\ph(\tup a))}\cmp \change_{\ph(\tup a)}^{\tup c}
\]
for any representative $n$-tuple $\tup c$. Let $\psi_{\ph(\tup a)}\colon C_S(\ph(\tup a))\to C_S(\tup c)$ be any conjugation map in $S$ taking $\ph(\tup a)$ to the representative of its conjugacy class. By the description of $\change_{\ph(\tup a)}^{\tup c}$ in \cite[Lemma 3.12]{RSS_Bold1}, we then have
\[
\freeO n(M)_{\tup a,\tup c} = \sum_{\substack{R\leq C_S(\tup a)\\\ph\in\cF(R,S)\\\text{ s.t. $\tup a\in R$ and}\\\text{$\ph(\tup a)$ is $S$-conjugate to $\tup c$}}}c_{R,\ph}\cdot [R,\psi_{\ph(\tup a)}\circ \ph]_{C_S(\tup a)}^{C_S(\tup c)}.
\]
When we compose $(\Z/p^e)^n\x \freeO n(\omega_\cF)$ with $\twistO n(\omega_\cF)$, we take a sum over all conjugacy classes of commuting $n$-tuples:
\begin{align*}
&\Bigl(\bigl( (\Z/p^e)^n\x \freeO n(\omega_\cF)\bigr)\cmp \twistO n(\omega_\cF)\Bigr)_{\tup a,\tup b} 
\\={}& \sum_{[\tup c]\in \cntuples n S} \bigl( (\Z/p^e)^n\x \freeO n(\omega_\cF)_{\tup a,\tup c}\bigr)\cmp \twistO n(\omega_\cF)_{\tup c,\tup b}
\\={}& \sum_{[\tup c]\in \cntuples n S} \sum_{\substack{R\leq C_S(\tup a)\\\ph\in\cF(R,S)\\\text{ s.t. $\tup a\in R$ and}\\\text{$\ph(\tup a)$ is $S$-conjugate to $\tup c$}}}c_{R,\ph}\cdot \bigl((\Z/p^e)^n\x [R,\psi_{\ph(\tup a)}\circ \ph]_{C_S(\tup a)}^{C_S(\tup c)}\bigr)\cmp \twistO n(\omega_\cF)_{\tup c,\tup b}
\\={}& \sum_{\substack{R\leq C_S(\tup a)\\\ph\in\cF(R,S)\\\text{ s.t. $\tup a\in R$}}}c_{R,\ph}\cdot \bigl((\Z/p^e)^n\x [R,\psi_{\ph(\tup a)}\circ \ph]_{C_S(\tup a)}^{C_S(\psi_{\ph(\tup a)}(\ph(\tup a)))}\bigr)\cmp \twistO n(\omega_\cF)_{\psi_{\ph(\tup a)}(\ph(\tup a)),\tup b}
\\={}& \sum_{\substack{R\leq C_S(\tup a)\\\ph\in\cF(R,S)\\\text{ s.t. $\tup a\in R$}}}c_{R,\ph}\cdot \bigl((\Z/p^e)^n\x [R,\id]_{C_S(\tup a)}^R\bigr)
\\* &\hspace{2.5cm}\cmp \bigl((\Z/p^e)^n\x [R,\psi_{\ph(\tup a)}\circ \ph]_R^{C_S(\psi_{\ph(\tup a)}(\ph(\tup a)))}\bigr)\cmp \twistO n(\omega_\cF)_{\psi_{\ph(\tup a)}(\ph(\tup a)),\tup b}.
\end{align*}
In the last line we have split $[R,\psi_{\ph(\tup a)}\circ \ph]_{C_S(\tup a)}^{C_S(\psi_{\ph(\tup a)}(\ph(\tup a)))}$ into its transfer and homomorphism parts. This enables us to now apply \cite[Theorem 3.33.(iv)]{RSS_Bold1} to note that 
\begin{multline*}
    (\Z/p^e)^n\x [R,\psi_{\ph(\tup a)}\circ \ph]_R^{C_S(\psi_{\ph(\tup a)}(\ph(\tup a)))} \\= (\Z/p^e)^n\x \freeO n([R,\psi_{\ph(\tup a)}\circ \ph]_R^S)_{\tup a,\psi_{\ph(\tup a)}(\ph(\tup a))} = \twistO n([R,\psi_{\ph(\tup a)}\circ \ph]_R^S)_{\tup a,\psi_{\ph(\tup a)}(\ph(\tup a))}.
\end{multline*}
Functoriality of $\twistO n$ then gives us 
\begin{align*}
&\Bigl(\bigl( (\Z/p^e)^n\x \freeO n(\omega_\cF)\bigr)\cmp \twistO n(\omega_\cF)\Bigr)_{\tup a,\tup b} 
\\={}&  \sum_{\substack{R\leq C_S(\tup a)\\\ph\in\cF(R,S)\\\text{ s.t. $\tup a\in R$}}}c_{R,\ph}\cdot \bigl((\Z/p^e)^n\x [R,\id]_{C_S(\tup a)}^R\bigr)
\\* &\hspace{2.5cm}\cmp \twistO n([R,\psi_{\ph(\tup a)}\circ \ph]_R^S)_{\tup a,\psi_{\ph(\tup a)}(\ph(\tup a))}\cmp \twistO n(\omega_\cF)_{\psi_{\ph(\tup a)}(\ph(\tup a)),\tup b}
\\={}&  \sum_{\substack{R\leq C_S(\tup a)\\\ph\in\cF(R,S)\\\text{ s.t. $\tup a\in R$}}}c_{R,\ph}\cdot \bigl((\Z/p^e)^n\x [R,\id]_{C_S(\tup a)}^R\bigr)
\cmp \twistO n([R,\psi_{\ph(\tup a)}\circ \ph]_R^S\cmp \omega_\cF)_{\tup a,\tup b}.
\end{align*}
The characteristic idempotent $\omega_\cF\in \AF(S,S)$ is $\cF$-stable, so restricting $\omega_\cF$ along the homomorphism $(\psi_{\ph(\tup a)}\circ \ph)\in \cF(R,S)$ on the left is equivalent to restricting along the inclusion $\incl\colon R\to S$. This means that $[R,\psi_{\ph(\tup a)}\circ \ph]_R^S\cmp \omega_\cF = [R,\incl]_R^S\cmp \omega_\cF$
and we can run our intermediary calculations in reverse:
\begin{align*}
&\Bigl(\bigl( (\Z/p^e)^n\x \freeO n(\omega_\cF)\bigr)\cmp \twistO n(\omega_\cF)\Bigr)_{\tup a,\tup b} 
\\={}&  \sum_{\substack{R\leq C_S(\tup a)\\\ph\in\cF(R,S)\\\text{ s.t. $\tup a\in R$}}}c_{R,\ph}\cdot \bigl((\Z/p^e)^n\x [R,\id]_{C_S(\tup a)}^R\bigr)
\cmp \twistO n([R,\incl]_R^S\cmp \omega_\cF)_{\tup a,\tup b}
\\={}&  \sum_{\substack{R\leq C_S(\tup a)\\\ph\in\cF(R,S)\\\text{ s.t. $\tup a\in R$}}}c_{R,\ph}\cdot \bigl((\Z/p^e)^n\x [R,\id]_{C_S(\tup a)}^R\bigr)
\cmp \twistO n([R,\incl]_R^S)_{\tup a,\tup a}\cmp \twistO n(\omega_\cF)_{\tup a,\tup b}
\\={}&  \sum_{\substack{R\leq C_S(\tup a)\\\ph\in\cF(R,S)\\\text{ s.t. $\tup a\in R$}}}c_{R,\ph}\cdot \bigl((\Z/p^e)^n\x [R,\incl]_{C_S(\tup a)}^{C_S(\tup a)}\bigr)
\cmp \twistO n(\omega_\cF)_{\tup a,\tup b}
\\={}&  \sum_{\substack{R\leq C_S(\tup a)\text{ up to $C_S(\tup a)$-conj.}\\\text{s.t. $a\in R$}}} \biggl(\sum_{\substack{R'\sim_{C_S(\tup a)} R\\ \ph\in\cF(R',S)}} c_{R',\ph} \biggr)\cdot \bigl((\Z/p^e)^n\x [R,\incl]_{C_S(\tup a)}^{C_S(\tup a)}\bigr)
\cmp \twistO n(\omega_\cF)_{\tup a,\tup b}.
\end{align*}
Lemma \ref{lemmaIdempotentCoefficients} implies that the sum of coefficients is $0$, when $R< C_S(\tup a)$, and the sum of coefficients equals $1$, when $R=C_S(\tup a)$. We complete our calculation with 
\begin{align*}
&\Bigl(\bigl( (\Z/p^e)^n\x \freeO n(\omega_\cF)\bigr)\cmp \twistO n(\omega_\cF)\Bigr)_{\tup a,\tup b} 
\\={}&  \sum_{\substack{R\leq C_S(\tup a)\text{ up to $C_S(\tup a)$-conj.}\\\text{s.t. $a\in R$}}} \biggl(\sum_{\substack{R'\sim_{C_S(\tup a)} R\\ \ph\in\cF(R',S)}} c_{R',\ph} \biggr)\cdot \bigl((\Z/p^e)^n\x [R,\incl]_{C_S(\tup a)}^{C_S(\tup a)}\bigr)
\cmp \twistO n(\omega_\cF)_{\tup a,\tup b}
\\={}& \bigl((\Z/p^e)^n\x [C_S(\tup a),\id]_{C_S(\tup a)}^{C_S(\tup a)}\bigr)
\cmp \twistO n(\omega_\cF)_{\tup a,\tup b}
\\={}& \twistO n(\omega_\cF)_{\tup a,\tup b}.
\end{align*}
This completes the proof of part \ref{itemIdempotentCompTwist} of the proposition.

\medskip
The proof of part \ref{itemIdempotentCompFree} makes use of the exact same tricks as the proof of part \ref{itemIdempotentCompTwist}. By \cite[Proposition 3.31]{RSS_Bold1}, the matrix $\twistO n(\omega_\cF)\in \AF((\Z/p^e)^n\x \freeO nS, (\Z/p^e)^n\x \freeO nS)$ has entries
\begin{align*}
&\twistO n(\omega_\cF)_{\tup a,\tup c} 
\\={}& \sum_{\substack{R\leq C_S(\tup a)\\\ph\in \cF(R,S)\\ \text{ s.t. } \ph(\tup a^{k(\tup a,R)})\\\text{is $S$-conj. to $\tup c$}}} c_{R,\ph}\cdot 
\Bigl( [\ev_{\tup a}^{-1}(R),(\id_{(\Z/p^e)^n}\x\ph)\circ \wind(\tup a, R)]_{(\Z/p^e)^n\x C_S(\tup a)}^{(\Z/p^e)^n\x C_S(\ph(\tup a^{k(\tup a,R)}))} 
\\* & \hspace{2.5cm}\cmp ((\Z/p^e)^n\x \change_{\ph(\tup a^{k(\tup a,R)})}^{\tup c}) \Bigr)
\\={}& \sum_{\substack{R\leq C_S(\tup a)\\\ph\in \cF(R,S)\\ \text{ s.t. } \ph(\tup a^{k(\tup a,R)})\\\text{is $S$-conj. to $\tup c$}}} c_{R,\ph}\cdot 
 [\ev_{\tup a}^{-1}(R),(\id_{(\Z/p^e)^n}\x(\psi_{\ph(\tup a^{k(\tup a,R)})}\circ \ph)\circ \wind(\tup a, R)]_{(\Z/p^e)^n\x C_S(\tup a)}^{(\Z/p^e)^n\x C_S(\tup c)}, 
\end{align*}
where $\psi_{\ph(\tup a^{k(\tup a,R)})}\colon C_S(\ph(\tup a^{k(\tup a,R)})) \to C_S(\tup c)$ is an $S$-conjugation map taking $\ph(\tup a^{k(\tup a,R)})$ to the representative of its conjugacy class in $S$.
We will name this representative $\tup z$ to ease the notation in the following calculations. Note that $\tup z= \psi_{\ph(\tup a^{k(\tup a,R)})}(\ph(\tup a^{k(\tup a,R)}))$ represents the conjugacy class of $\ph(\tup a^{k(\tup a,R)})$ and as such depends on $R$, $\ph$, and $\tup a$.

We decompose $[\ev_{\tup a}^{-1}(R),(\id_{(\Z/p^e)^n}\x(\psi_{\ph(\tup a^{k(\tup a,R)})}\circ \ph)\circ \wind(\tup a, R)]$ into 
\[
    [\ev_{\tup a}^{-1}(R),\wind(\tup a,R)]_{(\Z/p^e)^n\x C_S(\tup a)}^{(\Z/p^e)^n\x R}\cmp \Bigl((\Z/p^e)^n\x[R,\psi_{\ph(\tup a^{k(\tup a,R)})}\circ \ph]_{R}^{C_S(\tup z)}\Bigr).
\]
Using the functoriality of $\freeO n$ and the $\cF$-stability of $\omega_\cF$, we then proceed as in part \ref{itemIdempotentCompTwist} to replace $\psi_{\ph(\tup a^{k(\tup a,R)})}\circ \ph$ with the inclusion $\incl\colon R\to S$. The main steps are the following:
\begin{align*}
    &\Bigl( \twistO n(\omega_\cF)\cmp \bigl( (\Z/p^e)^n\x \freeO n(\omega_\cF)\bigr)\Bigr)_{\tup a,\tup b}
\\={}&  \sum_{\substack{R\leq C_S(\tup a)\\\ph\in\cF(R,S)}}c_{R,\ph}\cdot [\ev_{\tup a}^{-1}(R),\wind(\tup a,R)]_{(\Z/p^e)^n\x C_S(\tup a)}^{(\Z/p^e)^n\x R}
\\* &\hspace{2.5cm}\cmp \Bigl((\Z/p^e)^n\x [R,\psi_{\ph(\tup a^{k(\tup a,R)})}\circ \ph]_R^{C_S(\tup z)}\Bigr)
\\* &\hspace{2.5cm}\cmp \bigl((\Z/p^e)^n\x \freeO n(\omega_\cF)_{\tup z,\tup b}\bigr)
\\={}& \sum_{\substack{R\leq C_S(\tup a)\\\ph\in\cF(R,S)}}c_{R,\ph}\cdot [\ev_{\tup a}^{-1}(R),\wind(\tup a,R)]_{(\Z/p^e)^n\x C_S(\tup a)}^{(\Z/p^e)^n\x R}
\\* &\hspace{2.5cm}\cmp \bigl((\Z/p^e)^n\x \freeO n([R,\psi_{\ph(\tup a^{k(\tup a,R)})}\circ \ph]_R^S\cmp \omega_\cF)_{\tup a^{k(\tup a,R)},\tup b}\bigr)
\\={}& \sum_{\substack{R\leq C_S(\tup a)\\\ph\in\cF(R,S)}}c_{R,\ph}\cdot [\ev_{\tup a}^{-1}(R),\wind(\tup a,R)]_{(\Z/p^e)^n\x C_S(\tup a)}^{(\Z/p^e)^n\x R}
\\* &\hspace{2.5cm}\cmp \bigl((\Z/p^e)^n\x \freeO n([R,\incl]_R^S\cmp \omega_\cF)_{\tup a^{k(\tup a,R)},\tup b}\bigr)
\\={}& \sum_{\substack{R\leq C_S(\tup a)\\\ph\in\cF(R,S)}}c_{R,\ph}\cdot [\ev_{\tup a}^{-1}(R),(\id_{(\Z/p^e)^n}\x \psi_{\tup a^{k(\tup a,R)}})\circ \wind(\tup a,R)]_{(\Z/p^e)^n\x C_S(\tup a)}^{(\Z/p^e)^n\x C_S(\psi_{\tup a^{k(\tup a,R)}}(\tup a^{k(\tup a,R)}))}
\\* &\hspace{2.5cm}\cmp \bigl((\Z/p^e)^n\x \freeO n(\omega_\cF)_{\psi_{\tup a^{k(\tup a,R)}}(\tup a^{k(\tup a,R)}),\tup b}\bigr)
\\={}&  \sum_{\substack{R\leq C_S(\tup a)\\\text{ up to $C_S(\tup a)$-conj.}}} \biggl(\sum_{\substack{R'\sim_{C_S(\tup a)} R\\ \ph\in\cF(R',S)}} c_{R',\ph} \biggr)\cdot
\\* &\hspace{2.5cm} [\ev_{\tup a}^{-1}(R),(\id_{(\Z/p^e)^n}\x \psi_{\tup a^{k(\tup a,R)}})\circ \wind(\tup a,R)]_{(\Z/p^e)^n\x C_S(\tup a)}^{(\Z/p^e)^n\x C_S(\psi_{\tup a^{k(\tup a,R)}}(\tup a^{k(\tup a,R)}))}
\\* &\hspace{2.5cm}\cmp \bigl((\Z/p^e)^n\x \freeO n(\omega_\cF)_{\psi_{\tup a^{k(\tup a,R)}}(\tup a^{k(\tup a,R)}),\tup b}\bigr).
\end{align*}
We again apply Lemma \ref{lemmaIdempotentCoefficients} to remove all summands except $R=C_S(\tup a)$. For $R=C_S(\tup a)$, we have $\tup a^{k(\tup a,R)} = \tup a$ and $\wind(\tup a,C_S(\tup a)) = \id_{(\Z/p^e)^n\x C_S(\tup a)}$. We finish the calculation with 
\begin{align*}
    &\Bigl( \twistO n(\omega_\cF)\cmp \bigl( (\Z/p^e)^n\x \freeO n(\omega_\cF)\bigr)\Bigr)_{\tup a,\tup b}
\\={}&  [\ev_{\tup a}^{-1}(C_S(\tup a)),(\id_{(\Z/p^e)^n}\x \psi_{\tup a})\circ \wind(\tup a,C_S(\tup a))]_{(\Z/p^e)^n\x C_S(\tup a)}^{(\Z/p^e)^n\x C_S(\psi_{\tup a}(\tup a))}
\\* &\hspace{2.5cm}\cmp \bigl((\Z/p^e)^n\x \freeO n(\omega_\cF)_{\tup a,\tup b}\bigr)
\\={}&  [(\Z/p^e)^n\x C_S(\tup a),\id]_{(\Z/p^e)^n\x C_S(\tup a)}^{(\Z/p^e)^n\x C_S(\tup a)}
\\* &\hspace{2.5cm}\cmp \bigl((\Z/p^e)^n\x \freeO n(\omega_\cF)_{\tup a,\tup b}\bigr)
\\={}&  (\Z/p^e)^n\x \freeO n(\omega_\cF)_{\tup a,\tup b}.
\end{align*}
This completes the proof of part \ref{itemIdempotentCompFree}.
\end{proof}

\begin{cor}\label{corEquivalentTelescopes}
The idempotents $\twistO n ((\omega_\cF)_S^S)$ and $(\Z/p^e)^n\x \freeO n((\omega_\cF)_S^S)$ induce inverse equivalences
\[
\begin{tikzpicture}
\node (M) [matrix of math nodes] {
\tel_{\twistO n (\omega_\cF)} &[3cm] \Sinfpp B((\Z/p^e)^n)\wedge \tel_{\freeO n(\omega_\cF)}. \\
};
\path[->,arrow,auto]
    (M-1-1.north east) edge[bend left=30] node{$(\Z/p^e)^n\x \freeO n((\omega_\cF)_S^S)$} (M-1-2.north west)
    (M-1-2.south west) edge[bend left=30] node{$\twistO n((\omega_\cF)_S^S)$} (M-1-1.south east)
;
\end{tikzpicture}
\]
\end{cor}

\begin{proof}
By Proposition \ref{propComposingLoopIdempotents}, we can apply  $(\Z/p^e)^n\x \freeO n((\omega_\cF)_S^S)$ and $\twistO n((\omega_\cF)_S^S)$ level-wise to get maps back and forth between the two towers
\[
\Sinfpp  B((\Z/p^e)^n\x \freeO n S)\xto{(\Z/p^e)^n\x \freeO n(\omega_\cF)} \Sinfpp  B((\Z/p^e)^n\x \freeO n S) \xto{(\Z/p^e)^n\x \freeO n(\omega_\cF)} \dotsb
\]
and
\[
\Sinfpp  B((\Z/p^e)^n\x \freeO n S)\xto{\twistO n(\omega_\cF)} \Sinfpp  B((\Z/p^e)^n\x \freeO n S) \xto{\twistO n(\omega_\cF)} \dotsb.
\]
Again by Proposition \ref{propComposingLoopIdempotents}, the composite $\twistO n((\omega_\cF)_S^S)\cmp ((\Z/p^e)^n\x \freeO n((\omega_\cF)_S^S)) = (\Z/p^e)^n\x \freeO n((\omega_\cF)_S^S)$ is the identity on the telescope $\Sinfpp B((\Z/p^e)^n)\wedge \tel_{\freeO n(\omega_\cF)}$ of the first tower, and the composite $((\Z/p^e)^n\x \freeO n((\omega_\cF)_S^S))\cmp \twistO n((\omega_\cF)_S^S) = \twistO n((\omega_\cF)_S^S)$ is the identity on the telescope $\tel_{\twistO n(\omega_\cF)}$ of the second tower.
\end{proof}

\begin{cor}\label{corTwistedLoopTelescopeEquiv}
Combining Corollaries \ref{corFreeLoopTelescopeEquiv} and \ref{corEquivalentTelescopes}, we have inverse equivalences
\[
\begin{tikzpicture}
\node (M) [matrix of math nodes] {
\tel_{\twistO n (\omega_\cF)} &[3cm] \Sinfpp B((\Z/p^e)^n\x \freeO n\cF), \\
};
\path[->,arrow,auto]
    (M-1-1.north east) edge[bend left=30] node{$(\Z/p^e)^n\x I_\cF$} (M-1-2.north west)
    (M-1-2.south west) edge[bend left=30] node{$\twistT\cF$} (M-1-1.south east)
;
\end{tikzpicture}
\]
with matrices $(\Z/p^e)^n\x I_\cF\in \AF((\Z/p^e)\x \freeO n S, (\Z/p^e)^n\x \freeO n \cF)$ and\linebreak $\twistT\cF \in \AF((\Z/p^e)\x \freeO n \cF, (\Z/p^e)^n\x \freeO n S)$. These matrices have the following entries:
\[
((\Z/p^e)^n\x I_\cF)_{\tup a,\tup b} = \begin{cases}
0 & \text{if $\tup a$ is not $\cF$-conjugate to $\tup b$,}
\\ (\Z/p^e)^n\x \change_{\tup a}^{\tup b} & \text{if $\tup b$ is the representative for the $\cF$-conjugacy class of $\tup a$,}
\end{cases}
\]
with $((\Z/p^e)^n\x I_\cF)_{\tup a,\tup b}\in \AF((\Z/p^e)^n\x C_S(\tup a), (\Z/p^e)^n\x C_\cF(\tup b)) $,
and
\[
(\twistT\cF)_{\tup a,\tup b} = \twistO n((\omega_\cF)_S^S)_{\tup a,\tup b},
\]
with $(\twistT\cF)_{\tup a,\tup b}\in \AF((\Z/p^e)^n\x C_\cF(\tup a),(\Z/p^e)^n\x C_S(\tup b))$.
\end{cor}

\begin{proof}
We compose the maps of Corollaries \ref{corFreeLoopTelescopeEquiv} and \ref{corEquivalentTelescopes}. In one direction we have the composite 
\[
\tel_{\twistO n(\omega_\cF)} \xto{(\Z/p^e)^n\x \freeO n((\omega_\cF)_S^S)} \Sinfpp B(\Z/p^e)^n \wedge \tel_{\freeO n(\omega_\cF)} \xto {(\Z/p^e)^n\x I_\cF} \Sinfpp B((\Z/p^e)^n\x \freeO n\cF).
\]
Due to Lemma \ref{lemmaIdempotentOfRetract}, we have
\[\freeO n((\omega_\cF)_S^S) \cmp I_\cF = I_\cF\cmp T_\cF\cmp I_\cF = I_\cF.\]
Hence the composed equivalence above is simply 
\[
\tel_{\twistO n(\omega_\cF)} \xto {(\Z/p^e)^n\x I_\cF} \Sinfpp B((\Z/p^e)^n\x \freeO n\cF)
\]
induced by $(\Z/p^e)^n\x I_\cF\in \AF((\Z/p^e)^n\x \freeO n S, (\Z/p^e)^n\x \freeO n \cF)$, which is invariant with respect to composition with $\twistO n((\omega_\cF)_S^S)$ on the left.

In the other direction, we have the composite
\[
\Sinfpp B((\Z/p^e)^n\x \freeO n\cF) \xto{(\Z/p^e)^n\x T_\cF} \Sinfpp B(\Z/p^e)^n \wedge \tel_{\freeO n(\omega_\cF)} \xto {\twistO n((\omega_\cF)_S^S)} \tel_{\twistO n(\omega_\cF)}.
\]
By Definition \ref{defLnFRetractOfLnS}, the matrix $T_\cF\in \AF(\freeO n \cF,\freeO n S)$ has entries $(T_\cF)_{\tup a,\tup b}=\freeO n((\omega_\cF)_S^S)_{\tup a,\tup b}$, whenever $\tup a$ represents the class $[\tup a]\in \cntuples n \cF$ and $\tup b$ represents $[\tup b]\in \cntuples n S$. Applying Proposition \ref{propComposingLoopIdempotents}\ref{itemIdempotentCompTwist}, we then calculate the entries of the composite equivalence:
\begin{align*}
    \bigl(((\Z/p^e)^n\x T_\cF)\cmp \twistO n((\omega_\cF)_S^S)\bigr)_{\tup a,\tup b} &= 
    \sum_{[\tup c]\in \cntuples nS} ((\Z/p^e)^n\x T_\cF)_{\tup a,\tup c} \cmp \twistO n((\omega_\cF)_S^S)_{\tup c,\tup b} 
    \\ &= \sum_{[\tup c]\in \cntuples nS} ((\Z/p^e)^n\x \freeO n((\omega_\cF)_S^S)_{\tup a,\tup c} \cmp \twistO n((\omega_\cF)_S^S)_{\tup c,\tup b} 
    \\ &= \bigl(((\Z/p^e)^n\x \freeO n((\omega_\cF)_S^S))\cmp \twistO n((\omega_\cF)_S^S)\bigr)_{\tup a,\tup b}
    \\ &= \twistO n((\omega_\cF)_S^S)_{\tup a,\tup b}.
\end{align*}
This is the matrix $\twistT \cF\in \AF((\Z/p^e)^n\x \freeO n\cF, (\Z/p^e)^n\x \freeO nS)$ described in the statement of the corollary. Furthermore $\twistT\cF$ is invariant with respect to $\twistO n((\omega_\cF)_S^S)$ on the right.
\end{proof}

We define the functor $\twistO n$ for saturated fusoids similarly to Proposition \ref{propSimpleFusionLoop}:
\begin{prop}\label{propTwistedFusionLoop}
The functor $\twistO n$ on unions of $p$-groups extends to a functor $\twistO n\colon \AF\to \AF$ given on objects by $\cF\mapsto (\Z/p^e)^n\x \freeO n \cF$ and on morphisms $X\in \AF(\cE,\cF)$ by the matrix with entries
\[\twistO n(X_\cE^\cF)_{\tup a,\tup b} =\sum_{\substack{[\tup b']\in \cntuples nS \\  \tup b' \sim_\cF \tup b}} (\twistO n(X_R^S))_{\tup a,\tup b'} \cmp ((\Z/p^e)^n\x \change_{b'}^b),\]
with $\twistO n(X_\cE^\cF)_{\tup a,\tup b}\in \AF((\Z/p^e)^n\x C_\cE(\tup a),(\Z/p^e)^n\x C_\cF(\tup b))$, and where $R$ and $S$ are the underlying unions of $p$-groups for $\cE$ and $\cF$ respectively.
\end{prop}

\begin{remark} \label{remarkFusionTwistedLoopFunctor}
In the proof we shall see that $\twistO n$ can alternatively be given as the composite $\twistO n(X_\cE^\cF) = \twistT\cE\cmp \twistO n(X_R^S) \cmp ((\Z/p^e)^n\x I_\cF)$.
\end{remark}

\begin{proof}
We first note that $\twistO n(X_\cE^\cF) = \twistT\cE\cmp \twistO n(X_R^S) \cmp ((\Z/p^e)^n\x I_\cF)$. We check this by calculating the entries of the right hand side:
\begin{align*}
    & \bigl( \twistT\cE\cmp \twistO n(X_R^S) \cmp ((\Z/p^e)^n\x I_\cF) \bigr)_{\tup a,\tup b}
\\ ={}& \sum_{[\tup c]\in \cntuples nR} (\twistT\cE)_{\tup a,\tup c} \cmp \bigl( \twistO n(X_R^S) \cmp ((\Z/p^e)^n\x I_\cF) \bigr)_{\tup c,\tup b}
\\ ={}& \sum_{[\tup c]\in \cntuples nR} (\twistO n((\omega_\cE)_R^R))_{\tup a,\tup c} \cmp \bigl( \twistO n(X_R^S) \cmp ((\Z/p^e)^n\x I_\cF) \bigr)_{\tup c,\tup b}
\\ ={}& \bigl( \twistO n((\omega_\cE)_R^R)\cmp \twistO n(X_R^S) \cmp ((\Z/p^e)^n\x I_\cF) \bigr)_{\tup a,\tup b}
\\ ={}& \bigl( \twistO n((\omega_\cE)_R^R\cmp X_R^S) \cmp ((\Z/p^e)^n\x I_\cF) \bigr)_{\tup a,\tup b}
\\ ={}& \bigl( \twistO n( X_R^S) \cmp ((\Z/p^e)^n\x I_\cF) \bigr)_{\tup a,\tup b}
\\ ={}& \sum_{[\tup d]\in \cntuples nS} \twistO n( X_R^S)_{\tup a,\tup d}\cmp ((\Z/p^e)^n\x I_\cF)_{\tup d,\tup b}
\\ ={}& \sum_{\substack{[\tup b']\in \cntuples nS \\  \tup b' \sim_\cF \tup b}} (\twistO n(X_R^S))_{\tup a,\tup b'} \cmp ((\Z/p^e)^n\x \change_{b'}^b)
\\ ={}& \twistO n(X_\cE^\cF)_{\tup a,\tup b}.
\end{align*}

It is now straightforward to check that $\twistO n$ is a functor. First note that $\twistO n$ takes identity maps to identity maps since 
\begin{multline*}
\twistO n((\omega_\cF)_\cF^\cF) =\twistT\cF\cmp \twistO n((\omega_\cF)_S^S) \cmp ((\Z/p^e)^n\x I_\cF) \\= \twistT\cF \cmp ((\Z/p^e)^n\x I_\cF) = \id_{(\Z/p^e)^n\x\freeO n\cF}.
\end{multline*}
We have used the fact that $\twistT\cF$ is the inverse to $((\Z/p^e)^n\x I_\cF)$ by Corollary \ref{corTwistedLoopTelescopeEquiv}.

Let $\cE$, $\cF$, and $\cG$ be saturated fusoids over $R$, $S$, and $T$ respectively. Suppose $X\in \AF(\cE,\cF)$ and $Y\in \AF(\cF,\cG)$. We check that $\twistO n$ preserves composition:
\begin{align*}
    \twistO n(X_\cE^\cF)\cmp \twistO n(Y_\cF^\cG) &= \twistT \cE\cmp \twistO n(X_R^S)\cmp ((\Z/p^e)^n\x I_\cF)\cmp \twistT \cF \cmp \twistO n(Y_S^T) \cmp ((\Z/p^e)^n\x I_\cG)
    \\ &= \twistT \cE\cmp \twistO n(X_R^S)\cmp \twistO n((\omega_\cF)_S^S) \cmp \twistO n(Y_S^T) \cmp ((\Z/p^e)^n\x I_\cG)
    \\ &= \twistT \cE\cmp \twistO n(X_R^S\cmp (\omega_\cF)_S^S \cmp Y_S^T) \cmp ((\Z/p^e)^n\x I_\cG)
    \\ &= \twistT \cE\cmp \twistO n((X\cmp Y)_S^T) \cmp ((\Z/p^e)^n\x I_\cG)
    \\ &= \twistO n((X\cmp Y)_\cE^\cG).
\end{align*}
The characteristic idempotent $\omega_\cF$ disappears from the middle because  $X$ and $Y$ are $\cF$-stable. 
\end{proof}

We next give a formula for $\twistO n$ in terms of the decomposition of $X\in \AF(\cE,\cF)$ into basis elements. The formula is analogous to \cite[Proposition 3.31]{RSS_Bold1} and the proof follows the same lines as the proof of Corollary \ref{corSimpleFusionLoopOrbits}.

\begin{prop}\label{propFusionTwistedLoopOrbits}
Let $\cE$ and $\cF$ be saturated fusion systems over $p$-groups $R$ and $S$ respectively, and suppose $X\in \AF(\cE,\cF)$ is a virtual biset. Furthermore, let $\tup a$ in $\cE$ and  $\tup b$ in $\cF$ be chosen representatives for conjugacy classes of commuting $n$-tuples (according to Convention \ref{conventionTupleReps}). Consider the restriction of $X$ to the centralizer fusion system $C_\cE(\tup a)$, and write $M_{C_\cE(\tup a)}^\cF$ as a linear combination of basis elements (recalling Convention \ref{conventionFusionOrbitDecomposition}):
\[
X_{C_\cE(\tup a)}^\cF =\sum_{(P,\ph)}c_{P,\ph} \cdot [P,\ph]_{C_\cE(\tup a)}^\cF,
\]
where $P\leq C_R(\tup a)$ and $\ph\colon P\to S$.

The matrix entry $\twistO n(X)_{\tup a,\tup b}$ then satisfies the formula
\begin{multline*}
\twistO n(X)_{\tup a,\tup b} = \hspace{-.4cm} \sum_{\substack{(P,\ph) \text{ s.t. } \ph(\tup a^{k(\tup a,P)})\\\text{is $\cF$-conj. to $\tup b$}}}\hspace{-.4cm} c_{P,\ph}\cdot 
\Bigl( [\ev_{\tup a}^{-1}(P),(\id_{(\Z/p^e)^n}\x\ph)\circ \wind(\tup a, P)]_{(\Z/p^e)^n\x C_\cE(\tup a)}^{(\Z/p^e)^n\x C_S(\ph(\tup a^{k(\tup a,P)}))} \\ \cmp ((\Z/p^e)^n\x \change_{\ph(\tup a^{k(\tup a,P)})}^{\tup b}) \Bigr),
\end{multline*}
with $P$, $\ph$, and $c_{P,\ph}$ as in the linear combination above, and with $k(\tup a,P)$ and $\wind(\tup a,P)$ as given in \cite[Definition 3.25]{RSS_Bold1} and \cite[Lemma 3.29]{RSS_Bold1}, respectively.
\end{prop}

\begin{proof}
As in the proof of Corollary \ref{corSimpleFusionLoopOrbits}, we first wish to replace the coefficients $c_{P,\ph}$ with an alternative set of coefficients that play nicely with the underlying $p$-groups. View $X$ as an element of $\AF(C_R(\tup a),S)$ and write $X_{C_R(\tup a)}^S$ as a linear combination:
\[
X_{C_R(\tup a)}^S = \sum_{(P,\ph)}u_{P,\ph} \cdot [P,\ph]_{C_R(\tup a)}^S,
\]
with a (possibly) different collection of coefficients $u_{P,\ph}$.
If we precompose with $\omega_{C_\cE(\tup a)}$ and postcompose with $\omega_\cF$, the linear combination above becomes
\[
X_{C_\cE(\tup a)}^\cF = \omega_{C_\cE(\tup a)}\cmp X_{C_R(\tup a)}^S \cmp \omega_{\cF} = \sum_{(P,\ph)}u_{P,\ph} \cdot [P,\ph]_{C_\cE(\tup a)}^\cF.
\]
We shall then prove that the formula for $\freeO n(X)_{\tup a,\tup b}$ in the statement of the proposition is independent of the choice of linear combination for $X_{C_\cE(\tup a)}^\cF$, so that we can use the coefficients $u_{P,\ph}$ instead of $c_{P,\ph}$.

As in the proof of Corollary \ref{corSimpleFusionLoopOrbits}, it suffices to prove that whenever two basis elements $[P,\ph]_{C_\cE(\tup a)}^\cF$ and $[P',\ph']_{C_\cE(\tup a)}^\cF$ are equal, then the composite
\[
[\ev_{\tup a}^{-1}(P),(\id_{(\Z/p^e)^n}\x\ph)\circ \wind(\tup a, P)]_{(\Z/p^e)^n\x C_\cE(\tup a)}^{(\Z/p^e)^n\x C_S(\ph(\tup a^{k(\tup a,P)}))} \cmp ((\Z/p^e)^n\x \change_{\ph(\tup a^{k(\tup a,P)})}^{\tup b})
\]
and the same composite for the pair $(P',\ph')$ are also equal.

Suppose $[P,\ph]_{C_\cE(\tup a)}^\cF=[P',\ph']_{C_\cE(\tup a)}^\cF$ so that $(P',\ph')$ is $(C_\cE(\tup a),\cF)$-conjugate to $(P,\ph)$. Suppose further that $\ph(\tup a^{k(\tup a,P)})$ is $\cF$-conjugate to $\tup b$. Let $\ph'=\gamma\circ \ph\circ \alpha$ with $\alpha\in C_\cE(\tup a)$ and $\gamma\in \cF$, and where $\alpha\colon \gen{\tup a}P'\xto\cong \gen{\tup a}P$ satisfies $\alpha(\tup a)=\tup a$.

Via the isomorphism $\alpha$, it is clear that powers of elements in $\tup a$ lie in $P$ if and only if they lie in $P'$. Hence $k(\tup a,P)=k(\tup a,P')$, and 
\[\wind(\tup a,P)\circ (\id_{(\Z/p^e)^n}\x \alpha) = (\id_{(\Z/p^e)^n}\x \alpha)\circ \wind(\tup a,P')\] 
as maps $\ev^{-1}_{\tup a}(P') \xto\cong (\Z/p^e)^n\x P$.
In $S$ we have $\ph(\tup a^{k(\tup a,P')})=\gamma(\ph'(\tup a^{k(\tup a, P')}))$, and since $\gamma\in \cF$, we conclude that $\ph'(\tup a^{k(\tup a,P')})$ is also $\cF$-conjugate to $\tup b$. Furthermore, we have
\[[(\Z/p^e)^n\x P,\id\x \gamma]\cmp ((\Z/p^e)^n\x\change_{\ph'(\tup a^{k(\tup a,P')})}^{\tup b}) =  (\Z/p^e)^n\x\change_{\ph(\tup a^{k(\tup a,P)})}^{\tup b}.\]
Combining these observations, we get 
\begin{align*}
    & [\ev_{\tup a}^{-1}(P'),(\id_{(\Z/p^e)^n}\x\ph')\circ \wind(\tup a, P')]_{(\Z/p^e)^n\x C_\cE(\tup a)}^{(\Z/p^e)^n\x C_S(\ph'(\tup a^{k(\tup a,P')}))} \cmp ((\Z/p^e)^n\x \change_{\ph'(\tup a^{k(\tup a,P')})}^{\tup b})
    \\ &{}= [\ev_{\tup a}^{-1}(P'),(\id_{(\Z/p^e)^n}\x(\gamma\circ\ph\circ\alpha))\circ \wind(\tup a, P')]_{(\Z/p^e)^n\x C_\cE(\tup a)}^{(\Z/p^e)^n\x C_S(\ph'(\tup a^{k(\tup a,P')}))} \cmp ((\Z/p^e)^n\x \change_{\ph'(\tup a^{k(\tup a,P')})}^{\tup b})
    \\ &{}= [\ev_{\tup a}^{-1}(P'),(\id_{(\Z/p^e)^n}\x\ph)\circ \wind(\tup a, P)\circ (\id_{(\Z/p^e)^n}\x\alpha)]_{(\Z/p^e)^n\x C_\cE(\tup a)}^{(\Z/p^e)^n\x C_S(\ph(\tup a^{k(\tup a,P)}))} \cmp ((\Z/p^e)^n\x \change_{\ph(\tup a^{k(\tup a,P)})}^{\tup b})
    \\ &{}= [\ev_{\tup a}^{-1}(P),(\id_{(\Z/p^e)^n}\x\ph)\circ \wind(\tup a, P)]_{(\Z/p^e)^n\x C_\cE(\tup a)}^{(\Z/p^e)^n\x C_S(\ph(\tup a^{k(\tup a,P)}))} \cmp ((\Z/p^e)^n\x \change_{\ph(\tup a^{k(\tup a,P)})}^{\tup b}),
\end{align*}
where the last equality follows from the fact that $\alpha\in C_\cE(\tup a)$.

As in the proof of Corollary \ref{corSimpleFusionLoopOrbits}, the equality 
\[ \sum_{(P,\ph)}c_{P,\ph} \cdot [P,\ph]_{C_\cE(\tup a)}^\cF = X_{C_\cE(\tup a)}^\cF = \sum_{(P,\ph)}u_{P,\ph} \cdot [P,\ph]_{C_\cE(\tup a)}^\cF\]
implies that we can replace $c_{P,\ph}$ with $u_{P,\ph}$ in the formula
\begin{multline*}
\hspace{-.4cm} \sum_{\substack{(P,\ph) \text{ s.t. } \ph(\tup a^{k(\tup a,P)})\\\text{is $\cF$-conj. to $\tup b$}}}\hspace{-.4cm} c_{P,\ph}\cdot 
\Bigl( [\ev_{\tup a}^{-1}(P),(\id_{(\Z/p^e)^n}\x\ph)\circ \wind(\tup a, P)]_{(\Z/p^e)^n\x C_\cE(\tup a)}^{(\Z/p^e)^n\x C_S(\ph(\tup a^{k(\tup a,P)}))} \\ \cmp ((\Z/p^e)^n\x \change_{\ph(\tup a^{k(\tup a,P)})}^{\tup b}) \Bigr)
\end{multline*}
and get the same sum.

By Proposition \ref{propTwistedFusionLoop}, we can write $\twistO n(X_\cE^\cF)_{\tup a,\tup b}$ as
\begin{align*}
    \twistO n(X_\cE^\cF)_{\tup a,\tup b} &=  \sum_{\substack{[\tup b']\in \cntuples n S\\\tup b'\sim_\cF \tup b}} \twistO n(X_R^S)_{\tup a,\tup b'}\cmp ((\Z/p^e)^n\x \change_{\tup b'}^{\tup b}).
\end{align*}
We then apply Proposition \ref{propFusionTwistedLoopOrbits} with the linear combination for $X_{C_R(\tup a)}^S$ given by the coefficients $u_{P,\ph}$:
\begin{align*}
    &\twistO n(X_\cE^\cF)_{\tup a,\tup b} 
    \\={}& \sum_{\substack{[\tup b']\in \cntuples n S\\\tup b'\sim_\cF \tup b}} \twistO n(X_R^S)_{\tup a,\tup b'}\cmp ((\Z/p^e)^n\x \change_{\tup b'}^{\tup b})
    \\ ={}& \sum_{\substack{[\tup b']\in \cntuples n S\\\tup b'\sim_\cF \tup b}}    \sum_{\substack{(P,\ph) \text{ s.t. } \ph(\tup a^{k(\tup a,P)})\\\text{is $S$-conj. to $\tup b'$}}}\hspace{-.4cm} u_{P,\ph}\cdot \Bigl( [\ev_{\tup a}^{-1}(P),(\id_{(\Z/p^e)^n}\x\ph)\circ \wind(\tup a, P)]_{(\Z/p^e)^n\x C_R(\tup a)}^{(\Z/p^e)^n\x C_S(\ph(\tup a^{k(\tup a,P)}))} 
    \\*& \cmp ((\Z/p^e)^n\x \change_{\ph(\tup a^{k(\tup a,P)})}^{\tup b'}) \Bigr)\cmp ((\Z/p^e)^n\x \change_{\tup b'}^{\tup b})
    \\={}& \hspace{-.4cm} \sum_{\substack{(P,\ph) \text{ s.t. } \ph(\tup a^{k(\tup a,P)})\\\text{is $\cF$-conj. to $\tup b$}}}\hspace{-.4cm} u_{P,\ph}\cdot \Bigl( [\ev_{\tup a}^{-1}(P),(\id_{(\Z/p^e)^n}\x\ph)\circ \wind(\tup a, P)]_{(\Z/p^e)^n\x C_R(\tup a)}^{(\Z/p^e)^n\x C_S(\ph(\tup a^{k(\tup a,P)}))} 
    \\*& \cmp ((\Z/p^e)^n\x \change_{\ph(\tup a^{k(\tup a,P)})}^{\tup b}) \Bigr).
\end{align*}
Finally, note that $\twistO n(X_\cE^\cF)_{\tup a,\tup b}\in \AF((\Z/p^e)^n\x C_\cE(\tup a),(\Z/p^e)^n\x C_\cF(\tup b))$ is left\linebreak $((\Z/p^e)^n\x C_\cE(\tup a))$-stable and as such does not change if we precompose with the idempotent $(\Z/p^e)^n\x \omega_{C_\cE(\tup a)}$:
\begin{align*}
    & \twistO n(X_\cE^\cF)_{\tup a,\tup b} 
    \\ ={}& ((\Z/p^e)^n\x \omega_{C_\cE(\tup a)})\cmp \twistO n(X_\cE^\cF)_{\tup a,\tup b} 
    \\= {}& \hspace{-.4cm} \sum_{\substack{(P,\ph) \text{ s.t. } \ph(\tup a^{k(\tup a,P)})\\\text{is $\cF$-conj. to $\tup b$}}}\hspace{-.4cm} u_{P,\ph}\cdot \omega_{C_\cE(\tup a)}
    \\*& \cmp \Bigl( [\ev_{\tup a}^{-1}(P),(\id_{(\Z/p^e)^n}\x\ph)\circ \wind(\tup a, P)]_{(\Z/p^e)^n\x C_R(\tup a)}^{(\Z/p^e)^n\x C_S(\ph(\tup a^{k(\tup a,P)}))} 
     \cmp ((\Z/p^e)^n\x \change_{\ph(\tup a^{k(\tup a,P)})}^{\tup b}) \Bigr)
    \\={}& \hspace{-.4cm} \sum_{\substack{(P,\ph) \text{ s.t. } \ph(\tup a^{k(\tup a,P)})\\\text{is $\cF$-conj. to $\tup b$}}}\hspace{-.4cm} u_{P,\ph}\cdot \Bigl( [\ev_{\tup a}^{-1}(P),(\id_{(\Z/p^e)^n}\x\ph)\circ \wind(\tup a, P)]_{(\Z/p^e)^n\x C_\cE(\tup a)}^{(\Z/p^e)^n\x C_S(\ph(\tup a^{k(\tup a,P)}))} 
    \\*& \cmp ((\Z/p^e)^n\x \change_{\ph(\tup a^{k(\tup a,P)})}^{\tup b}) \Bigr)
    \\={}&\hspace{-.4cm} \sum_{\substack{(P,\ph) \text{ s.t. } \ph(\tup a^{k(\tup a,P)})\\\text{is $\cF$-conj. to $\tup b$}}}\hspace{-.4cm} c_{P,\ph}\cdot \Bigl( [\ev_{\tup a}^{-1}(P),(\id_{(\Z/p^e)^n}\x\ph)\circ \wind(\tup a, P)]_{(\Z/p^e)^n\x C_\cE(\tup a)}^{(\Z/p^e)^n\x C_S(\ph(\tup a^{k(\tup a,P)}))} 
    \\*& \cmp ((\Z/p^e)^n\x \change_{\ph(\tup a^{k(\tup a,P)})}^{\tup b}) \Bigr).\qedhere
\end{align*}
\end{proof}

\begin{remark}\label{remarkFusionUntwisted}
As in \cite[Remark 3.32]{RSS_Bold1}, we can form a functor $\untwistO n$ by taking the formula in Proposition \ref{propFusionTwistedLoopOrbits} and leaving out all summands indexed by $(P,\ph)$, where $P$ does not contain $\tup a$. By \cite[Lemma 3.27]{RSS_Bold1}, the summands for which $\tup a\in P$ are precisely those summands where $\ev^{-1}_{\tup a}(P)=(\Z/p^e)^n\x P$, and these are also the summands for which $k(\tup a,P)_i=1$, for all $1\leq i\leq n$ .

For saturated fusion systems  $\cE$ and $\cF$ over $p$-groups $R$ and $S$, and for $X\in \AF(\cE,\cF)$, we thus define $\untwistO n(X)\in \AF((\Z/p^e)^n\x \cE, (\Z/p^e)^n\x \cF)$ to be the matrix with entries as in Proposition \ref{propFusionTwistedLoopOrbits} except leaving out all summands where $\tup a$ is not in $P$:
\begin{align*}
\untwistO n(X)_{\tup a,\tup b} ={}& \hspace{-.4cm} \sum_{\substack{(P,\ph) \text{ s.t. }\tup a\in P\text{ and} \\ \ph(\tup a^{k(\tup a,P)})\text{ is $\cF$-conj. to $\tup b$}}}\hspace{-.4cm} c_{P,\ph}\cdot 
\Bigl( [\ev_{\tup a}^{-1}(P),(\id_{(\Z/p^e)^n}\x\ph)\circ \wind(\tup a, P)]_{(\Z/p^e)^n\x C_\cE(\tup a)}^{(\Z/p^e)^n\x C_S(\ph(\tup a^{k(\tup a,P)}))} 
\\* &\cmp ((\Z/p^e)^n\x \change_{\ph(\tup a^{k(\tup a,P)})}^{\tup b}) \Bigr)
\\ ={}& \hspace{-.4cm} \sum_{\substack{(P,\ph) \text{ s.t. } \tup a\in P\text{ and}\\ \ph(\tup a)\text{ is $\cF$-conj. to $\tup b$}}}\hspace{-.4cm} c_{P,\ph}\cdot 
\Bigl( [(\Z/p^e)^n\x P,(\id_{(\Z/p^e)^n}\x\ph)]_{(\Z/p^e)^n\x C_\cE(\tup a)}^{(\Z/p^e)^n\x C_S(\ph(\tup a))} 
\\* &\cmp ((\Z/p^e)^n\x \change_{\ph(\tup a)}^{\tup b}) \Bigr).
\end{align*}
Comparing with Corollary \ref{corSimpleFusionLoopOrbits}, we see that $\untwistO n(X)$ coincides with $(\Z/p^e)^n\x \freeO n(X)$. As mentioned, we observed this for groups in \cite[Remark 3.32]{RSS_Bold1} and for the category $\Cov$ in \cite[Remark 2.11]{RSS_Bold1}.
\end{remark}

\section{Properties of $\twistO n$}\label{secFusionMainTheorem}

Before we state the fusion system version of \cite[Theorems 2.13 and 3.33]{RSS_Bold1}, we need to discuss the auxiliary maps describing evaluation and partial evaluation for $\freeO n\cF$, the action of $\Sigma_n$ on $\freeO n\cF$, and the embedding of $(\Z/p^e)^{n+m}\x \freeO {n+m}\cF$ into $(\Z/p^e)^n\x \freeO n((\Z/p^e)^m\x \freeO m\cF)$.

We shall handle all four auxiliary maps at once according to the following common framework: First the auxiliary map on $p$-groups gives rise to a map between towers corresponding to suitable idempotents. The map of towers induces a map between the mapping telescopes in spectra, and the telescopes are each equivalent to classifying spectra of fusoids. We calculate the induced map between the classifying spectra of fusoids and confirm that this map takes the form that we would expect for the auxiliary map in question.

Constructing the auxiliary maps in terms of maps between towers of idempotents instead of just defining the auxiliary maps for fusion systems directly has one significant advantage: It allows us to formally, using standard methods, turn natural transformations for $p$-groups into natural transformations for fusion systems (see Proposition \ref{propExtendNaturalTransformation}). This allows us to use \cite[Theorem 3.33]{RSS_Bold1} for $p$-groups to prove a significant part of Theorem \ref{thmFusionMain} for fusion systems.

\begin{lemma}\label{lemmaInducedTelescopeEquivalence}
Let $F\colon \AF\to \AF$ be a functor, and let $\cF$ be a saturated fusoid with underlying union of $p$-groups $S$. Consider the idempotent endomorphism $F((\omega_\cF)_S^S)\in \AF(F(S),F(S))$ and the associated mapping telescope \[\tel_{F(\omega_\cF)} = \colim(\Sinfpp BF(S) \xto{F((\omega_\cF)_S^S)} \Sinfpp BF(S) \xto{F((\omega_\cF)_S^S)} \dotsb).\]
We then have a pair of inverse equivalences
\[
\begin{tikzpicture}
\node (M) [matrix of math nodes] {
\tel_{F (\omega_\cF)} &[3cm] \Sinfpp BF(\cF). \\
};
\path[->,arrow,auto]
    (M-1-1.north east) edge[bend left=30] node{$F((\omega_\cF)_S^\cF)$} (M-1-2.north west)
    (M-1-2.south west) edge[bend left=30] node{$F((\omega_\cF)_\cF^S)$} (M-1-1.south east)
;
\end{tikzpicture}
\]
\end{lemma}

\begin{proof}
Consider the two towers
\[S\xto{(\omega_\cF)_S^S} S\xto{(\omega_\cF)_S^S}\dotsb\]
as well as 
\[\cF\xto{\id_\cF} \cF \xto{\id_\cF} \dotsb.\]
We can apply $(\omega_\cF)_\cF^S$ and $(\omega_\cF)_S^\cF$ level-wise to get maps back and forth between the towers, where we recall that $\id_\cF$ is simply $(\omega_\cF)_\cF^\cF\in \AF(\cF,\cF)$. The induced maps between telescopes $\tel_{\omega_\cF}$ and $\Sinfpp B\cF$ are equivalences (or even the identity if we used $\tel_{\omega_\cF}$ as the construction of $\Sinfpp B\cF$).

If we apply the functor $F\colon \AF\to \AF$ to the elements $(\omega_\cF)_\cF^S$ and $(\omega_\cF)_S^\cF$, we can apply the resulting maps $F((\omega_\cF)_\cF^S)$ and $F((\omega_\cF)_S^\cF)$ level-wise to the towers
\[F(S)\xto{F((\omega_\cF)_S^S)} F(S)\xto{F(\omega_\cF)_S^S)}\dotsb\]
and
\[F(\cF)\xto{\id_{F(\cF)}} F(\cF) \xto{\id_{F(\cF)}} \dotsb.\]
The composites $F((\omega_\cF)_S^\cF)\cmp F((\omega_\cF)_\cF^S) = F((\omega_\cF)_S^S)$ and $F((\omega_\cF)_\cF^S)\cmp F((\omega_\cF)_S^\cF) = F((\omega_\cF)_\cF^\cF)= \id_{F(\cF)}$ recover the idempotents of the towers, so the induced maps on telescopes 
\[
\begin{tikzpicture}
\node (M) [matrix of math nodes] {
\tel_{F (\omega_\cF)} &[3cm] \Sinfpp BF(\cF) \\
};
\path[->,arrow,auto]
    (M-1-1.north east) edge[bend left=30] node{$F((\omega_\cF)_S^\cF)$} (M-1-2.north west)
    (M-1-2.south west) edge[bend left=30] node{$F((\omega_\cF)_\cF^S)$} (M-1-1.south east)
;
\end{tikzpicture}
\]
are inverse to each other in $\Ho(\Spp)$.
\end{proof}

\begin{remark}
If we apply Lemma \ref{lemmaInducedTelescopeEquivalence} to the functors $\freeO n$ and $\twistO n$, we recover Corollaries \ref{corFreeLoopTelescopeEquiv} and \ref{corTwistedLoopTelescopeEquiv} with the same equivalences. However, we need those corollaries in order to construct $\freeO n$ and $\twistO n$ as functors $\AF\to \AF$ in the first place.
\end{remark}

\begin{definition}\label{defInducedNaturalTransformation}
Suppose $F,G\colon \AF\to \AF$ are functors defined on all saturated fusoids, and suppose we have a natural transformation $\eta\colon F|_{\textup{$p$-groups}} \Rightarrow G|_{\textup{$p$-groups}}$ defined only on formal unions of $p$-groups. For each saturated fusoid $\cF$ over $S$, we then define a map $\eta_\cF\colon F(\cF) \to G(\cF)$ as the composite
\[\eta_\cF \colon F(\cF) \xto{F((\omega_\cF)_\cF^S)} F(S) \xto{\eta_S} G(S) \xto{G((\omega_\cF)_S^\cF)} G(\cF).\]
When $\cF$ is the trivial fusion system on $S$, this just recovers $\eta_S$, so there is no ambiguity of notation.
\end{definition}

Since $\eta$ is a natural transformation on $p$-groups, $\eta_S$ fits in a diagram of towers:
\[
\begin{tikzpicture}
\node (M) [matrix of math nodes] {
    F(\cF) &[2.5cm] F(\cF) &[2.5cm] \dotsb \\[1.5cm]
    F(S) & F(S) & \dotsb\\[1.5cm]
    G(S) & G(S) & \dotsb\\[1.5cm]
    G(\cF) & G(\cF) & \dotsb.\\
};
\path[->,arrow,auto]
    (M-1-1) edge node{$\id_{F(\cF)}$} (M-1-2)
            edge node{$F((\omega_\cF)_\cF^S)$} (M-2-1)
    (M-1-2) edge node{$\id_{F(\cF)}$} (M-1-3)
            edge node{$F((\omega_\cF)_\cF^S)$} (M-2-2)
    (M-2-1) edge node{$F((\omega_\cF)_S^S)$} (M-2-2)
            edge node{$\eta_S$} (M-3-1)
    (M-2-2) edge node{$F((\omega_\cF)_S^S)$} (M-2-3)
            edge node{$\eta_S$} (M-3-2)
    (M-3-1) edge node{$G((\omega_\cF)_S^S)$} (M-3-2)
            edge node{$G((\omega_\cF)_S^\cF)$} (M-4-1)
    (M-3-2) edge node{$G((\omega_\cF)_S^S)$} (M-3-3)
            edge node{$G((\omega_\cF)_S^\cF)$} (M-4-2)
    (M-4-1) edge node{$\id_{G(\cF)}$} (M-4-2)
    (M-4-2) edge node{$\id_{G(\cF)}$} (M-4-3)
;
\end{tikzpicture}
\]
The map $\eta_S$ induces a map between the telescopes $\eta_S\colon \tel_{F(\omega_\cF)}\to \tel_{G(\omega_\cF)}$, and as such $\eta_\cF$ is simply the composite
\[\eta_\cF\colon \Sinfpp BF(\cF) \xto[\simeq]{F((\omega_\cF)_\cF^S)} \tel_{F(\omega_\cF)} \xto{\eta_S} \tel_{G(\omega_\cF)} \xto[\simeq]{G((\omega_\cF)_S^\cF)} \Sinfpp BG(\cF)\]
in $\Ho(\Spp)$. It is now easy to prove that the extension of $\eta$ to fusoids defines a natural transformation $F\Rightarrow G$ on all of $\AF$.
\begin{prop}\label{propExtendNaturalTransformation}
Suppose $F,G\colon \AF\to \AF$ are functors defined on all saturated fusoids, and suppose we have a natural transformation $\eta\colon F|_{\textup{$p$-groups}} \Rightarrow G|_{\textup{$p$-groups}}$ defined only on formal unions of $p$-groups.
If we extend $\eta$ to all saturated fusoids by Definition \ref{defInducedNaturalTransformation}, then the extension defines a natural transformation $\eta\colon F\Rightarrow G$ on all of $\AF$.
\end{prop}

\begin{proof}
Let $\cE$ and $\cF$ be saturated fusoids over $R$ and $S$ respectively, and let $X\in \AF(\cE,\cF)$ be any matrix of virtual bisets. We have to prove that $F(X_\cE^\cF) \cmp \eta_\cF = \eta_\cE \cmp G(X_\cE^\cF)$ in $\AF(F(\cE),G(\cF))$.
Since $\AF\to \Ho(\Spp)$ is fully faithful, it suffices to prove this as homotopy classes of maps $\Sinfpp BF(\cE) \to \Sinfpp BG(\cF)$. 

By the naturality of $\eta$ on $p$-groups as well as the definitions of $\eta_\cE$ and $\eta_\cF$, we have the following commutative diagram in $\Ho(\Spp)$:
\[
\begin{tikzpicture}
\node (M) [matrix of math nodes] {
    \Sinfpp BF(\cE) &[2cm] &[2cm] &[2cm] \Sinfpp BG(\cE) \\[2cm]
    & \tel_{F(\omega_\cE)}  & \tel_{G(\omega_\cE)} & \\[2cm]
    & \tel_{F(\omega_\cF)}  & \tel_{G(\omega_\cF)} & \\[2cm]
    \Sinfpp BF(\cF) & & &  \Sinfpp BG(\cF).\\
};
\path[->,arrow,auto]
    (M-1-1) edge node{$\eta_\cE$} (M-1-4)
            edge node{$F((\omega_\cE)_\cE^R)$} node[swap]{$\simeq$} (M-2-2)
            edge node{$F(X_\cE^\cF)$} (M-4-1)
    (M-2-2) edge node{$\eta_R$} (M-2-3)
            edge node{$F(X_R^S)$} (M-3-2)
    (M-2-3) edge node{$G((\omega_\cE)_R^\cE)$} node[swap]{$\simeq$} (M-1-4)
            edge node{$G(X_R^S)$} (M-3-3)
    (M-1-4) edge node{$G(X_\cE^\cF)$} (M-4-4)
    (M-4-1) edge node[swap]{$F((\omega_\cF)_\cF^S$} node{$\simeq$} (M-3-2)
            edge node{$\eta_\cF$} (M-4-4)
    (M-3-2) edge node{$\eta_S$} (M-3-3)
    (M-3-3) edge node[swap]{$G((\omega_\cF)_S^\cF)$} node{$\simeq$} (M-4-4)
;
\end{tikzpicture}
\]
Since all the smaller squares commute, the outer square commutes as well.
\end{proof}

As the first of the four auxiliary maps needed for Theorem \ref{thmFusionMain}, let us describe the evaluation map from $(\Z/p^e)^n\x \freeO n\cF$ to $\cF$. 
Consider the two endofunctors $\twistO n, \Id_{\AF}\colon \AF\to \AF$. By Theorem \cite[Theorem 3.33.(v)]{RSS_Bold1}, the evaluation maps for unions of $p$-groups define a natural transformation $\ev \colon \twistO n|_{\textup{$p$-groups}}\Rightarrow \Id_{\AF}|_{\textup{$p$-groups}}$.
We can thus extend $\ev$ by Proposition \ref{propExtendNaturalTransformation} to a natural transformation $\ev\colon \twistO n\Rightarrow \Id_{\AF}$ on all of $\AF$. Let us determine a formula for the evaluation map $\ev_\cF\colon (\Z/p^e)^n\x \freeO n\cF\to \cF$ in this extension.

By Definition \ref{defInducedNaturalTransformation}, the biset matrix $\ev_\cF$ is given as the composite
\[
\ev_\cF = \twistO n((\omega_\cF)_\cF^S) \cmp \ev_S \cmp (\omega_\cF)_S^\cF
\]
inside $\AF((\Z/p^e)^n\x \freeO n \cF, \cF)$. In order to calculate the matrix entries of $\ev_\cF$, recall that $\twistO n((\omega_\cF)_\cF^S) = \twistT \cF$ by Remark \ref{remarkFusionTwistedLoopFunctor}, and, by Corollary \ref{corEquivalentTelescopes}, we have $(\twistT \cF)_{\tup a,\tup b}=\twistO n((\omega_\cF)_S^S)_{\tup a,\tup b}$, whenever $\tup a$ represents a conjugacy class of tuples in $\cF$ and $\tup b$ represents a class in $S$.

Calculating the entries of $\ev_\cF\in \AF((\Z/p^e)^n\x \freeO n \cF, \cF)$, we then get
\begin{align*}
    (\ev_\cF)_{\tup a} &= \sum_{[\tup b]\in \cntuples n S} \twistO n((\omega_\cF)_\cF^S)_{\tup a,\tup b} \cmp [(\Z/p^e)^n\x C_S(\tup b),\ev_{\tup b}]_{(\Z/p^e)^n\x C_S(\tup b)}^S \cmp (\omega_\cF)_S^\cF
    \\ &= \sum_{[\tup b]\in \cntuples n S} \twistO n((\omega_\cF)_S^S)_{\tup a,\tup b} \cmp [(\Z/p^e)^n\x C_S(\tup b),\ev_{\tup b}]_{(\Z/p^e)^n\x C_S(\tup b)}^S \cmp (\omega_\cF)_S^\cF
    \\ &= (\twistO n((\omega_\cF)_S^S) \cmp \ev_S \cmp (\omega_\cF)_S^\cF)_{\tup a}
    \\ &= (\ev_S \cmp (\omega_\cF)_S^S\cmp (\omega_\cF)_S^\cF)_{\tup a}
    \\ &= (\ev_S \cmp (\omega_\cF)_S^\cF)_{\tup a}
    \\ &= [(\Z/p^e)^n\x C_S(\tup a),\ev_{\tup a}]_{(\Z/p^e)^n\x C_S(\tup a)}^S \cmp (\omega_\cF)_S^\cF
    \\ &= [(\Z/p^e)^n\x C_S(\tup a),\ev_{\tup a}]_{(\Z/p^e)^n\x C_S(\tup a)}^\cF
\end{align*}
as an element of $\AF((\Z/p^e)^n\x C_\cF(\tup a), \cF)$. 

By Lemma \ref{lemmaEvalFusionPreserving}, the homomorphism $\ev_{\tup a}$ is fusion preserving from $(\Z/p^e)^n\x C_\cF(\tup a)$ to $\cF$.
Lemma \ref{lemmaFusionPreserving} now implies that we can also write
\[(\ev_\cF)_{\tup a} = [(\Z/p^e)^n\x C_S(\tup a),\ev_{\tup a}]_{(\Z/p^e)^n\x C_\cF(\tup a)}^\cF.\]

Thus, the evaluation map $\ev_\cF$ simply applies the fusion preserving map $\ev_{\tup a}\colon (\Z/p^e)^n\x C_\cF(\tup a)\to \cF$ to each component of $(\Z/p^e)^n\x \freeO n\cF$. Let us record this fact for future reference:
\begin{lemma}\label{lemmaFusionEvaluation}
Let $\cF$ be a saturated fusion system (or fusoid) over a finite $p$-group $S$ (or formal union of such). The evaluation map $\ev_\cF\colon (\Z/p^e)^n\x\freeO n\cF \to \cF$ takes each component $(\Z/p^e)^n\x C_\cF(\tup a)$ to $\cF$ via the fusion preserving map $\ev_{\tup a}$. As such, the entries of $\ev_\cF$ are given by
\[(\ev_\cF)_{\tup a} = [(\Z/p^e)^n\x C_S(\tup a),\ev_{\tup a}]_{(\Z/p^e)^n\x C_\cF(\tup a)}^\cF.\]
Furthermore, the maps $\ev_\cF$ assemble into a natural transformation $\ev\colon \twistO n\Rightarrow \Id_{\AF}$ between endofunctors on $\AF$.
\end{lemma}
We could have easily defined $\ev_\cF$ directly for fusion systems in terms of the fusion preserving maps $(\Z/p^e)^n\x C_\cF(\tup a)\to \cF$. However, by constructing $\ev_\cF$ via Proposition \ref{propExtendNaturalTransformation}, we know that $\ev_\cF$ gives rise to a natural transformation which will greatly simplify the proof of Theorem \ref{thmFusionMain}.

The next auxiliary map we shall work on is the partial evaluation map $\pev_\cF\colon \Z/p^e\x \freeO {n+1}\cF \to \freeO n \cF$. Note that for $p$-groups the partial evaluation map $\pev_S\colon \Z/p^e\x \freeO {n+1}S \to \freeO n S$ of \cite[Theorem 3.33.(vi)]{RSS_Bold1}, given by
\[\pev_S(t,z) =  (a_{n+1})^t\cdot z\in C_S(a_1,\dotsc,a_n), \quad \text{for } t\in \Z/p^e,  z\in \coprod_{\tup a\in \cntuples{n+1}S} C_S(\tup a),\]
coincides with the 1-fold evaluation map for $\freeO n S$, i.e. we have $\pev_S = \ev_{\freeO n S}\colon \Z/p^e\x \freeO 1(\freeO n S) \to \freeO n S$. As such, $\pev_S$ provides a natural transformation $\pev\colon \twistO 1(\freeO n(-))\Rightarrow \freeO n(-)$ on $p$-groups. By Lemma \ref{lemmaFusionEvaluation}, the maps $\pev_\cF = \ev_{\freeO n \cF} \colon \Z/p^e\x \freeO {n+1}\cF \to \freeO n \cF$ provide a natural transformation $\pev\colon \twistO 1(\freeO n(-))\Rightarrow \freeO n(-)$ on all of $\AF$. The entries of the biset matrix, $(\pev_\cF)_{\tup a,\tup b}\in \AF(\Z/p^e \x C_\cF(\tup a),C_\cF(\tup b))$ have the following form for each representative $(n+1)$-tuple $\tup a = (a_1,\dotsc, a_n,a_{n+1})$:
\[
(\pev_\cF)_{\tup a,\tup b} =\begin{cases} 
[\Z/p^e \x C_S(\tup a), \pev_{\tup a}]_{\Z/p^e \x C_\cF(\tup a)}^{C_\cF(a_1,\dotsc,a_n)} &\text{if }\tup b=(a_1,\dotsc,a_n), 
\\ 0 &\text{otherwise.}
\end{cases}
\]
Note that $(a_1,\dotsc,a_n)$ is itself a chosen representative $n$-tuple for $\cF$ by Convention \ref{conventionTupleReps}.

According to \cite[Theorem 3.33.(vi)]{RSS_Bold1}, the partial evaluation maps give a natural transformation $(\Z/p^e)^n\x \pev\colon \twistO {n+1} \Rightarrow \twistO n$ for all unions of $p$-groups. By Proposition \ref{propExtendNaturalTransformation}, we can extend this natural transformation to all saturated fusoids to get $\eta\colon \twistO {n+1} \Rightarrow \twistO n$ on $\AF$ with $\eta_\cF$ given by
\[\eta_\cF = \twistO{n+1}((\omega_\cF)_\cF^S) \cmp ((\Z/p^e)^n\x \pev_S) \cmp \twistO n((\omega_\cF)_S^\cF).\]
We claim that $\eta_\cF$ simply recovers $(\Z/p^e)^n\x \pev_\cF$ for fusoids as well. The calculation is completely analogous to the calculation preceding Lemma \ref{lemmaFusionEvaluation}.

We note that $\twistO{n+1}((\omega_\cF)_\cF^S)_{\tup a,\tup b} = \twistO{n+1}((\omega_\cF)_S^S)_{\tup a,\tup b}$ for all $(n+1)$-tuple representatives $\tup a$ for $\cF$ and $\tup b$ for $S$. Additionally, Remark \ref{remarkFusionTwistedLoopFunctor} implies that $\twistO n((\omega_\cF)_S^\cF)= (\Z/p^e)^n\x I_\cF$ has entries
\[
\twistO n((\omega_\cF)_S^\cF)_{\tup a,\tup b} =
\begin{cases}
(\Z/p^e)^n\x \change_{\tup a}^{\tup b} &\text{if $\tup a$ is $\cF$-conj. to $\tup b$,}
\\ 0 &\text{otherwise.}
\end{cases}
\]
The calculation of $\eta_\cF$ the proceeds accordingly:
\begin{align*}
    &(\eta_\cF)_{\tup a,\tup b} 
\\ ={}& \Bigl(\twistO{n+1}((\omega_\cF)_\cF^S) \cmp ((\Z/p^e)^n\x \pev_S) \cmp \twistO n((\omega_\cF)_S^\cF)\Bigr)_{\tup a,\tup b}
\\ ={}& \Bigl(\twistO{n+1}((\omega_\cF)_S^S) \cmp ((\Z/p^e)^n\x \pev_S) \cmp \twistO n((\omega_\cF)_S^\cF)\Bigr)_{\tup a,\tup b}
\\ ={}& \Bigl(((\Z/p^e)^n\x \pev_S) \cmp \twistO{n}((\omega_\cF)_S^S) \cmp \twistO n((\omega_\cF)_S^\cF)\Bigr)_{\tup a,\tup b}
\\ ={}& \Bigl(((\Z/p^e)^n\x \pev_S) \cmp \twistO n((\omega_\cF)_S^\cF)\Bigr)_{\tup a,\tup b}
\\ ={}& \begin{cases} 
\!\begin{aligned}
 [(\Z/p^e)^{n+1} \x C_S(\tup a), (\Z/p^e)^n\x \pev_{\tup a}&]_{(\Z/p^e)^{n+1} \x C_S(\tup a)}^{(\Z/p^e)^n \x C_S(a_1,\dotsc,a_n)} 
\\*& \cmp ((\Z/p^e)^n\x \change_{(a_1,\dotsc,a_n)}^{\tup b})
\end{aligned}&\text{if $(a_1,\dotsc,a_n)\sim_\cF \tup b$}, 
\\ 0 &\text{otherwise.}
\end{cases}
\end{align*}
However if $(a_1,\dotsc,a_n)$ is $\cF$-conjugate to $\tup b$, then the two $n$-tuples are equal as they both represent the same $\cF$-conjugacy class, hence $\change_{(a_1,\dotsc,a_n)}^{\tup b}$ is the inclusion $C_S(a_1,\dotsc,a_n) \to C_\cF(a_1,\dotsc,a_n)$.

In total we have
\[
(\eta_\cF)_{\tup a,\tup b} = \begin{cases}
[(\Z/p^e)^{n+1} \x C_S(\tup a), (\Z/p^e)^n\x \pev_{\tup a}]_{(\Z/p^e)^{n+1} \x C_S(\tup a)}^{(\Z/p^e)^n \x C_\cF(a_1,\dotsc,a_n)} &\text{if }\tup b=(a_1,\dotsc,a_n), 
\\ 0 &\text{otherwise.}
\end{cases}
\]
By Lemma \ref{lemmaEvalFusionPreserving}, $\pev_{\tup a}\colon \Z/p^e \x C_\cF(\tup a) \to C_\cF(a_1,\dotsc,a_n)$ is fusion preserving, so 
\begin{multline*}
[(\Z/p^e)^{n+1} \x C_S(\tup a), (\Z/p^e)^n\x \pev_{\tup a}]_{(\Z/p^e)^{n+1} \x C_S(\tup a)}^{(\Z/p^e)^n \x C_\cF(a_1,\dotsc,a_n)} 
\\= [(\Z/p^e)^{n+1} \x C_S(\tup a), (\Z/p^e)^n\x \pev_{\tup a}]_{(\Z/p^e)^{n+1} \x C_\cF(\tup a)}^{(\Z/p^e)^n \x C_\cF(a_1,\dotsc,a_n)}
\end{multline*}
and $\eta_\cF$ recovers $(\Z/p^e)^n\x \pev_\cF$ as claimed.
Again we record this as a lemma:
\begin{lemma}\label{lemmaFusionPartialEvaluation}
Let $\cF$ be a saturated fusion system (or fusoid) over a finite $p$-group $S$ (or formal union of such). The partial evaluation map $\pev_\cF\colon \Z/p^e\x\freeO{n+1}\cF \to \freeO n\cF$ takes each component $\Z/p^e\x C_\cF(\tup a)$ to $C_\cF(a_1,\dotsc,a_n)$ via the fusion preserving map $\pev_{\tup a}$ given by
\[
\pev_{\tup a}(t,z) = (a_{n+1})^t\cdot z.
\]
The maps $(\Z/p^e)^n\x \pev_\cF$ define a natural transformation $(\Z/p^e)^n\x \pev\colon \twistO{n+1}\Rightarrow \twistO n$ between endofunctors on $\AF$.
\end{lemma}

Next we establish the action of $\Sigma_n$ on $\freeO n \cF$ in a similar fashion.
Given a permutation $\sigma\in \Sigma_n$, the action of $\sigma$ on $\freeO n S$ permutes the components (as described in \cite[Theorem 3.33.(iii)]{RSS_Bold1}) by sending $C_S(\tup a)=C_S(\sigma(\tup a))$ to $C_S(\tilde{\sigma(\tup a)})$ via the isomorphism $\change_{\sigma(\tup a)}^{\tilde{\sigma(\tup a)}}$, where $\tilde{\sigma(\tup a)}$ is the chosen representative for the $S$-conjugacy class of $\sigma(\tup a)$.

The action of $\sigma$ defines a natural transformation $\sigma\colon \freeO n\Rightarrow \freeO n$ for all unions of $p$-groups (note that \cite[Proposition 3.15]{RSS_Bold1} is invariant with respect to permuting the $n$-tuples by $\sigma$). Proposition \ref{propExtendNaturalTransformation} then provides an extension $\sigma\colon \freeO n\Rightarrow \freeO n$ to all of $\AF$.

As for the (partial) evaluation maps, we calculate that the induced map $\sigma_\cF = \freeO n((\omega_\cF^S))\cmp \sigma_S \cmp \freeO n((\omega_S^\cF))$ is a matrix in which the only non-zero entries are
\[(\sigma_\cF)_{\tup a,\tilde{\sigma(\tup a)}} = \change_{\sigma(\tup a)}^{\tilde{\sigma(\tup a)}} \in \AF(C_\cF(\tup a),C_\cF(\tilde{\sigma(\tup a)}),\]
and where $\tilde{\sigma(\tup a)}$ is the chosen representative for the $\cF$-conjugacy class of $\sigma(\tup a)$.

Similarly, by \cite[Theorem 3.33.(iii)]{RSS_Bold1} the diagonal action of $\sigma$ on $(\Z/p^e)^n\x \freeO n S$ provides a natural transformation $\sigma\colon \twistO n\Rightarrow \twistO n$ for unions of $p$-groups. 
Again we extend this natural transformation to a natural transformation $\eta\colon\twistO n\Rightarrow \twistO n$ on all of $\AF$, and we proceed to calculate that $\eta_\cF = \twistO n((\omega_\cF^S))\cmp \sigma_S \cmp \twistO n((\omega_S^\cF))$ is simply the diagonal action of $\sigma$ on $(\Z/p^e)^n\x \freeO n\cF$. We record this fact:
\begin{lemma}\label{lemmaFusionEquivariant}
Let $\cF$ be a saturated fusion system (or fusoid) over a finite $p$-group $S$ (or formal union of such).
The symmetric group $\Sigma_n$ acts on $\freeO n\cF$ by permuting the coordinates of the $n$-tuples. If $\sigma\in \Sigma_n$, then the action of $\sigma$ on $\freeO n\cF$ sends each component $C_\cF(\tup a)=C_\cF(\sigma(\tup a))$ to $C_\cF(\tilde{\sigma(\tup a)})$ via the isomorphism $\change_{\sigma(\tup a)}^{\tilde{\sigma(\tup a)}}$, where $\tilde{\sigma(\tup a)}$ is the chosen representative for the $\cF$-conjugacy class of $\sigma(\tup a)$.

The diagonal action of $\sigma$ on $(\Z/p^e)^n\x \freeO n\cF$ gives a natural transformation $\sigma\colon \twistO n\Rightarrow \twistO n$ between endofunctors on $\AF$.
\end{lemma}

The final auxiliary map needed for Theorem \ref{thmFusionMain} is the embedding of $(\Z/p^e)^{n+m}\x \freeO{n+m}\cF$ into $(\Z/p^e)^m\x \freeO m((\Z/p^e)^n\x \freeO n\cF)$. \cite[Theorem 3.33.(vii)]{RSS_Bold1} states that we have a natural transformation $\iota\colon \twistO{n+m}(-)\Rightarrow \twistO m(\twistO n(-))$ for all unions of $p$-groups, where $\iota_S\colon (\Z/p^e)^{n+m} \x \freeO{n+m}S \to (\Z/p^e)^m \x \freeO m((\Z/p^e)^n \x \freeO n S)$ takes each component
$(\Z/p^e)^{n+m} \x C_S(\tup x,\tup y)$ to the component $(\Z/p^e)^m\x C_{(\Z/p^e)^n\x C_S(\tup x)}(\tup 0 \x \tup y)$ via the map
\begin{multline*}
    \bigl((\tup s,\tup r), z\bigr)\in (\Z/p^e)^{n+m}\x C_S((\tup x,\tup y)) \mapsto \bigl(\tup r, (\tup s,z) \bigr) \in (\Z/p^e)^m\x C_{(\Z/p^e)^n\x C_S(\tup x)}(\tup 0\x \tup y),
\end{multline*}
for $\tup s\in (\Z/p^e)^n$, $\tup r\in (\Z/p^e)^m$, $\tup x\in \cntuples n S$, $\tup y\in \cntuples m {C_S(\tup x)}$, and $z\in C_S(\tup x,\tup y)$.

Suppose $(\tup x,\tup y)$ is a chosen representative $(n+m)$-tuple in $S$. Then by Convention \ref{conventionTupleReps} the commuting $n$-tuple $\tup x$ is a chosen representative in $S$ as well, and the commuting $m$-tuple $\tup y$ is a chosen representative in $C_S(\tup x)\subseteq \freeO n S$ -- note that $\tup y$ is \emph{not} necessarily a chosen representative in $S$.

The entries of the biset matrix $\iota_S\in \AF((\Z/p^e)^{n+m} \x \freeO{n+m}S, (\Z/p^e)^m \x \freeO m((\Z/p^e)^n \x \freeO n S))$ are given by
\begin{align*}
& (\iota_S)_{(\tup x,\tup y)\in \cntuples{n+m} S, (\tup t\x \tup w)\in \cntuples m {((\Z/p^e)^n\x \freeO n S)}}
\\ = {}& \begin{cases}
[(\Z/p^e)^{n+m} \x C_S(\tup x,\tup y), ((\tup s,\tup r),z)\mapsto (\tup r, (\tup s,z))] & \substack{\text{if $(\tup t\x \tup w) = (\tup 0\x \tup y)$}
\\ \text{in the component $(\Z/p^e)^n\x C_S(\tup x)$,}}
\\ 0 & \text{otherwise.}
\end{cases}
\end{align*}
We extend $\iota$ to a natural transformation $\twistO{n+m}(-)\Rightarrow \twistO m(\twistO n(-))$ on all of $\AF$ by Proposition \ref{propExtendNaturalTransformation}, and as such we have
\[\iota_\cF = \twistO{n+m}((\omega_\cF)_\cF^S) \cmp \iota_S\cmp \twistO m(\twistO n((\omega_\cF)_S^\cF)).\]
As for the previous auxiliary maps, we first note that for any representative $(n+m)$-tuple $(\tup x,\tup y)$ in $\cF$ and $m$-tuple $(\tup t\x \tup w)$ in $(\Z/p^e)^n\x \freeO n \cF$ we have
\begin{align*}
    (\iota_\cF)_{(\tup x,\tup y),(\tup t\x \tup w)} &= (\twistO{n+m}((\omega_\cF)_\cF^S) \cmp \iota_S\cmp \twistO m(\twistO n((\omega_\cF)_S^\cF)))_{(\tup x,\tup y),(\tup t\x \tup w)}
    \\ &= (\iota_S\cmp \twistO m(\twistO n((\omega_\cF)_S^\cF)))_{(\tup x,\tup y),(\tup t\x \tup w)}
    \\ &= \sum_{\tup t'\x \tup w'\in \cntuples m {((\Z/p^e)^n\x \freeO n S)}} (\iota_S)_{(\tup x,\tup y),(\tup t'\x \tup w')} \cmp \twistO m(\twistO n((\omega_\cF)_S^\cF))_{(\tup t'\x \tup w'),(\tup t\x \tup w)}.
\end{align*}
Now the matrix entries of $\iota_S$ are zero except when $\tup t'\x \tup w' = \tup 0\x \tup y$ as $m$-tuples in the component $(\Z/p^e)^n\x C_S(\tup x)$ of $(\Z/p^e)^n\x \freeO n S$. Hence the above equation becomes
\begin{align*}
& (\iota_\cF)_{(\tup x,\tup y),(\tup t\x \tup w)} 
\\={}& (\iota_S)_{(\tup x,\tup y),(\tup 0\x \tup y)} \cmp \twistO m(\twistO n((\omega_\cF)_S^\cF))_{(\tup 0\x \tup y),(\tup t\x \tup w)}
\\ ={}& [(\Z/p^e)^{n+m} \x C_S(\tup x,\tup y), ((\tup s,\tup r),z)\mapsto (\tup r, (\tup s,z))]_{(\Z/p^e)^{n+m}\x C_S(\tup x,\tup y)}^{(\Z/p^e)^m\x C_{(\Z/p^e)^n\x C_S(\tup x)}(\tup 0\x \tup y)} 
\\* &\cmp \twistO m(\twistO n((\omega_\cF)_S^\cF))_{(\tup 0\x \tup y),(\tup t\x \tup w)}.
\end{align*}
If we apply Proposition \ref{propFusionTwistedLoopOrbits} twice to the basis element $(\omega_\cF)_S^\cF = [S,\id]_S^\cF \in \AF(S,\cF)$, we see that
\begin{multline*}
    \twistO m(\twistO n([S,\id]_S^\cF))_{(\tup 0\x \tup y),(\tup t\x \tup w)}
    \\ = \begin{cases}
    (\Z/p^e)^m\x \change_{(\tup 0\x \tup y)}^{(\tup t\x \tup w)} & \text{if $(\tup t\x \tup w)$ represents the $(\Z/p^e)^n\x C_\cF(\tup x)$-conjugacy class of $(\tup 0\x \tup y)$,}
    \\ 0 &\text{otherwise.}
    \end{cases}
\end{multline*}
Because $(\tup x,\tup y)$ was assumed to be a representative $(n+m)$-tuple in $\cF$, the $m$-tuple $\tup y$ is a representative in $C_\cF(\tup x)$. At the same time, if the $m$-tuple $\tup t$ in $(\Z/p^e)^n$ is conjugate to $\tup 0$, then $\tup t=\tup 0$. Consequently $(\tup 0\x \tup y)$ is the chosen representative for its $(\Z/p^e)^n\x C_\cF(\tup x)$-conjugacy class. 

We conclude that $(\iota_\cF)_{(\tup x,\tup y),(\tup t\x \tup w)}$ is zero unless $\tup t\x \tup w=\tup 0\x \tup y$ in the component $(\Z/p^e)^n\x C_\cF(\tup x)$ of $(\Z/p^e)^n\x \freeO n S$. Furthermore,
\begin{align*}
& (\iota_\cF)_{(\tup x,\tup y),(\tup 0\x \tup y)} 
\\ ={}& [(\Z/p^e)^{n+m} \x C_S(\tup x,\tup y), ((\tup s,\tup r),z)\mapsto (\tup r, (\tup s,z))]_{(\Z/p^e)^{n+m}\x C_S(\tup x,\tup y)}^{(\Z/p^e)^m\x C_{(\Z/p^e)^n\x C_S(\tup x)}(\tup 0\x \tup y)} 
\\* &\cmp (\Z/p^e)^m\x \change_{(\tup 0\x \tup y)}^{(\tup 0\x \tup y)} 
\\ ={}& [(\Z/p^e)^{n+m} \x C_S(\tup x,\tup y), ((\tup s,\tup r),z)\mapsto (\tup r, (\tup s,z))]_{(\Z/p^e)^{n+m}\x C_S(\tup x,\tup y)}^{(\Z/p^e)^m\x C_{(\Z/p^e)^n\x C_S(\tup x)}(\tup 0\x \tup y)} 
\\* &\cmp (\Z/p^e)^m\x [C_{(\Z/p^e)^n\x C_S(\tup x)}(\tup 0\x \tup y),\id]_{C_{(\Z/p^e)^n\x C_S(\tup x)}(\tup 0\x \tup y)}^{C_{(\Z/p^e)^n\x C_\cF(\tup x)}(\tup 0\x \tup y)}
\\ ={}& [(\Z/p^e)^{n+m} \x C_S(\tup x,\tup y), ((\tup s,\tup r),z)\mapsto (\tup r, (\tup s,z))]_{(\Z/p^e)^{n+m}\x C_S(\tup x,\tup y)}^{(\Z/p^e)^m\x C_{(\Z/p^e)^n\x C_\cF(\tup x)}(\tup 0\x \tup y)}.
\end{align*}
As in Lemma \ref{lemmaEvalFusionPreserving}, it follows that the map $((\tup s,\tup r),z)\mapsto (\tup r, (\tup s,z))$ is fusion preserving from $(\Z/p^e)^{n+m}\x C_\cF(\tup x,\tup y)$ to $(\Z/p^e)^m\x C_{(\Z/p^e)^n\x C_\cF(\tup x)}(\tup 0\x \tup y)$ and additionally we can write 
\[
(\iota_\cF)_{(\tup x,\tup y),(\tup 0\x \tup y)} = [(\Z/p^e)^{n+m} \x C_S(\tup x,\tup y), ((\tup s,\tup r),z)\mapsto (\tup r, (\tup s,z))]_{(\Z/p^e)^{n+m}\x C_\cF(\tup x,\tup y)}^{(\Z/p^e)^m\x C_{(\Z/p^e)^n\x C_\cF(\tup x)}(\tup 0\x \tup y)}
\]
if we prefer. We record our final preliminary result about the auxiliary maps of the main theorem.

\begin{lemma}\label{lemmaFusionLnIteration}
Let $\cF$ be a saturated fusion system (or fusoid) over a finite $p$-group $S$ (or formal union of such). The embedding $\iota_\cF\colon (\Z/p^e)^{n+m}\x\freeO{n+m}\cF \to (\Z/p^e)^m \x \freeO m((\Z/p^e)^n\x \freeO n\cF)$ takes each component $(\Z/p^e)^{n+m}\x C_\cF(\tup x,\tup y)$ to the component $(\Z/p^e)^m \x C_{(\Z/p^e)^n\x C_\cF(\tup x)}(\tup 0\x \tup y)$ via the fusion preserving map given by
\[
((\tup s,\tup r),z) \mapsto (\tup r,(\tup s,z)).
\]
The maps $\iota_\cF$ define a natural transformation $\twistO{n+m}(-)\Rightarrow \twistO m(\twistO n(-))$ between endofunctors on $\AF$.
\end{lemma}

Now we are finally ready to prove that our extension of $\twistO n$ to fusion systems satisfies the same properties as for finite groups and finite sheeted covering maps. Recalling \cite[Convention 3.21]{RSS_Bold1}, we state and prove the following variant of \cite[Theorems 2.13 and 3.33]{RSS_Bold1}:
\begin{theorem}\label{thmFusionMain}
The endofunctors $\twistO n \colon \AF\to \AF$ for $n\geq 0$ of Proposition \ref{propTwistedFusionLoop} have the following properties:
\begin{enumerate}
\renewcommand{\theenumi}{$(\roman{enumi})$}\renewcommand{\labelenumi}{\theenumi}
\item[$(\emptyset )$] Let $L^{\dagger,\AG}_n\colon \AG \to \AG$ be the functor constructed in \cite[Section 3]{RSS_Bold1}. When restricted to the full subcategories of $\AG$ and $\AF$ spanned by formal unions of finite $p$-groups, the functor $\twistO n\colon \AF \to \AF$ is the $\Z_p$-linearization of $L^{\dagger,\AG}_n$.
\item\label{itemFusionLnZero} $\twistO 0$ is the identity functor on $\AF$.
\item\label{itemFusionLnObjects} On objects, $\twistO n$ takes a saturated fusoid $\cF$ to the saturated fusoid
\[\twistO n(\cF) = (\Z/p^e)^n\x \freeO n (\cF)=\coprod_{\tup a\in \cntuples n\cF} (\Z/p^e)^n\x C_\cF(\tup a).\]
\item\label{itemFusionEquivariant} The group $\Sigma_n$ acts on $\freeO n\cF=\coprod_{\tup a\in \cntuples n\cF} C_\cF(\tup a)$ by permuting the coordinates of the $n$-tuples $\tup a$. Explicitly, if $\sigma\in \Sigma_n$ and if $\widetilde{\sigma(\tup a)}$ is the representative for the $\cF$-conjugacy class of $\sigma(\tup a)$, then $\sigma\colon \freeO n\cF\to \freeO n\cF$ maps $C_\cF(\tup a)=C_\cF(\sigma(\tup a))$ to $C_\cF(\widetilde{\sigma(\tup a)})$ via the isomorphism $\change_{\sigma(\tup a)}^{\widetilde{\sigma(\tup a)}}\in \AF(C_\cF(\tup a),C_\cF(\widetilde{\sigma(\tup a)}))$. 

The functor $\twistO n$ is equivariant with respect to the $\Sigma_n$-action on $(\Z/p^e)^n\x \freeO n(-)$ that permutes the coordinates of both $(\Z/p^e)^n$ and $\freeO n(-)$, i.e. for every $\sigma\in\Sigma_n$ the diagonal action of $\sigma$ on $(\Z/p^e)^n\x \freeO n(-)$ induces a natural isomorphism $\sigma\colon \twistO n\overset\cong\Rightarrow \twistO n$.
\item\label{itemFusionLnForwardMaps} Let $\cE$ and $\cF$ be saturated fusion systems on $R$ and $S$ respectively. For forward maps, i.e. transitive bisets $[R,\ph]_\cE^\cF\in \AF(\cE,\cF)$ with $\ph\colon R\to S$ fusion preserving, the functor $\twistO n$ coincides with $(\Z/p^e)^n \x \freeO n(-)$ so that 
\[\twistO n([R,\ph]_\cE^\cF)=(\Z/p^e)^n\x \freeO n([R,\ph]_\cE^\cF).\]
In addition, $\freeO n([R,\ph]_\cE^\cF)$ is the biset matrix that takes a component $C_\cE(\tup a)$ of $\freeO n\cE$ to the component $C_\cF(\tup b)$ of $\freeO n\cF$ by the biset 
\[
[C_R(\tup a),\ph]_{C_\cE(\tup a)}^{C_S(\ph(\tup a))}\cmp \change_{\ph(\tup a)}^{\tup b}\in \AF( C_\cE(\tup a), C_\cF(\tup b)),
\]
where $\tup b$ represents the $\cF$-conjugacy class of $\ph(\tup a)$.
\item\label{itemFusionEvalSquare} For all $n \geq 0$, the functor $\twistO n$ commutes with evaluation maps, i.e. the evaluation maps $\ev_\cF\colon (\Z/p^e)^n\x \freeO n(\cF)\to \cF$ form a natural transformation $\ev\colon \twistO n \Rightarrow \Id_{\AF}$.
\item\label{itemFusionLnPartialEvaluation} For all $n \geq 0$, the partial evaluation maps $\pev_\cF\colon \Z/p^e\x \freeO {n+1}(\cF) \to \freeO n (\cF)$ given as fusion preserving maps $\pev_{\tup a}\colon \Z/p^e \x C_\cF(\tup a) \to C_\cF(a_1,\dotsc, a_n)$ in terms of the formula
\[
\pev_{\tup a}(t,z)= (a_{n+1})^t\cdot z\in C_S(a_1,\dotsc,a_n), \quad \text{for } t\in \Z/p^e,  z\in C_S(\tup a),
\] form natural transformations $(\Z/p^e)^n\x \pev\colon \twistO {n+1} \Rightarrow \twistO n$.
\item\label{itemFusionIterateLn} For all $n,m\geq 0$, and any saturated fusoid $\cF$ on $S$, the formal union $(\Z/p^e)^{n+m}\x \freeO {n+m} \cF$ embeds into $(\Z/p^e)^m\x \freeO m((\Z/p^e)^n\x \freeO n \cF)$ as the components corresponding to the commuting $m$-tuples in $(\Z/p^e)^n\x \freeO n \cF$ that are zero in the $(\Z/p^e)^n$-coordinate, i.e. the embedding takes each component $(\Z/p^e)^{n+m}\x C_\cF(\tup x,\tup y)$ to the component $(\Z/p^e)^m \x C_{(\Z/p^e)^n\x C_\cF(\tup x)}(\tup 0\x \tup y)$, for $\tup x\in \cntuples n \cF$ and $\tup y\in \cntuples m \cF$, via the fusion preserving map given by
\[
((\tup s,\tup r),z) \mapsto (\tup r,(\tup s,z)),
\]
for $\tup s\in (\Z/p^e)^n$, $\tup r\in (\Z/p^e)^m$, and $z\in C_S(\tup x,\tup y)$.

These embeddings $(\Z/p^e)^{n+m}\x \freeO {n+m} \cF\to (\Z/p^e)^m\x \freeO m((\Z/p^e)^n\x \freeO n \cF)$ then form a natural transformation $\twistO {n+m}(-)\Rightarrow \twistO m(\twistO n(-))$.
\end{enumerate}
\end{theorem}

\begin{proof}
\ref{itemFusionLnZero},\ref{itemFusionLnObjects}: Both follow immediately from the definition of $\twistO n$ in Proposition \ref{propTwistedFusionLoop}. In particular for $n=0$ and $X\in \AF(\cE,\cF)$ there is only a single $0$-tuple $()$ in each component and
\[\twistO 0(X_\cE^\cF) = \twistO 0(X_R^S) = X\]
in Proposition \ref{propTwistedFusionLoop}.

\ref{itemFusionEvalSquare},\ref{itemFusionLnPartialEvaluation},\ref{itemFusionEquivariant},\ref{itemFusionIterateLn}: These are Lemmas \ref{lemmaFusionEvaluation}-\ref{lemmaFusionLnIteration} respectively.

\ref{itemFusionLnForwardMaps}: With $\ph\colon R\to S$ a fusion preserving map from $\cE$ to $\cF$, it follows from Lemma \ref{lemmaFusionPreserving} that 
\[
    X := [R,\ph]_\cE^\cF = [R,\ph]_R^\cF.
\]
Let $\tup a$ be a representative commuting $n$-tuple in $\cE$, then restricting to $C_R(\tup a)$ we have
\[  X_{C_\cE(\tup a)}^\cF = X_{C_R(\tup a)}^\cF = [C_R(\tup a), \ph|_{C_R(\tup a)}]_{C_R(\tup a)}^\cF = [C_R(\tup a), \ph|_{C_R(\tup a)}]_{C_\cE(\tup a)}^\cF,
\]
where we note that the restriction $\ph|_{C_R(\tup a)}$ is fusion preserving from $C_\cE(\tup a)$ to $\cF$.

Now Proposition \ref{propFusionTwistedLoopOrbits} and Corollary \ref{corSimpleFusionLoopOrbits} together state that
\begin{align*}
&\twistO n(X)_{\tup a,\tup b} 
\\ ={}& \left.\begin{cases} (\Z/p^e)^n\x ([C_R(\tup a), \ph]_{C_\cE(\tup a)}^{C_S(\ph(\tup a))} \cmp  \zeta_{\ph(\tup a)}^{\tup b} )& \text{if $\tup b$ is $\cF$-conjugate to $\ph(\tup a)$} \\ 0 & \text{otherwise} \end{cases}\right\}
\\ ={}& (\Z/p^e)^n \x \freeO n(X)_{\tup a,\tup b}.
\end{align*}
This completes the proof of \ref{itemFusionLnForwardMaps} and the theorem.
\end{proof}

\section{$\twistO n$ commutes with $p$-completion for bisets of finite groups}\label{secFusionPCompletionCommutes}
In this final section we compare the functor $\twistO n\colon \AF \to \AF$ with the functor $\twistO {n,p}\colon \AG\to \AG$ for finite groups, where we restrict to centralizers of commuting $n$-tuples of $p$-power order elements (see \cite[Proposition 3.35]{RSS_Bold1}).
We shall see that these two functors are closely related via the $p$-completion functor $(-)^\wedge_p\colon \AG\to \AF$ described in \cite{RSS_p-completion}, such that we have 
\[
(-)^\wedge_p \circ \twistO {n,p} = \twistO n \circ (-)^\wedge_p
\]
as functors $\AG\to \AF$.

Suppose $S$ and $T$ are Sylow $p$-subgroups of $G$ and $H$ respectively. Let $\cF_G=\cF_S(G)$ and $\cF_H=\cF_T(H)$ be the associated fusion systems at the prime $p$. The $p$-completed classifying spectrum $\Sinfpp B\cF_G$ is equivalent to the $p$-completion of $\Sinfp BG$ via the composite
\[\Sinfpp B\cF_G \xto{\omega_{\cF_G}} \Sinfpp BS \xto{(B\incl_S^G)^\wedge_p} \Sinfpp BG.\]
Here we first include the summand $\Sinfpp B\cF_G$ into $\Sinfpp BS$, and then include $S$ into $G$. The fact that this map is an equivalence is essentially due to \cite{CartanEilenberg}*{XII.10.1} and \cite{BLO2}*{Proposition 5.5} (see \cite{RSS_p-completion}*{Proposition 3.3} for additional details).

Via the Segal conjectures for finite groups and fusion systems, we can interpret the $p$-completion functor $(-)^\wedge_p$ as a functor $(-)^\wedge_p\colon \AG \to \AF$, where we use the particular equivalence $\Sinfpp BG \simeq \Sinfpp B\cF_G$ described above:
\begin{definition}
The functor $(-)^\wedge_p\colon \AG \to \AF$ is defined on morphisms by 
\[(-)^\wedge_p\colon \AG(G,H) \to [\Sinfp BG, \Sinf BH] \xto{(-)^\wedge_p} [\Sinfpp B\cF_G, \Sinfpp B\cF_H] \cong \AF(\cF_G,\cF_H).\]
The first map is the Segal map for finite groups, and the last isomorphism is the Segal conjecture for saturated fusion systems.
\end{definition}

If we restrict the $(H,H)$-biset $H$ to a $(T,T)$-biset $H_T^T$ the resulting biset is stable with respect to $\cF_H$, hence $H_{\cF_H}^{\cF_H}\in \AF(\cF_H,\cF_H)$. The element $H_{\cF_H}^{\cF_H}$ is always invertible in $\AF(\cF_H,\cF_H)$ by \cite{RSS_p-completion}*{Lemma 3.6}, and using the inverse, we get the following algebraic formula for the $p$-completion functor $(-)^\wedge_p\colon \AG\to \AF$:
\begin{prop}[\cite{RSS_p-completion}*{Theorem 1.1}]\label{propPCompletion}
The $p$-completion functor $(-)^\wedge_p \colon \AG \to \AF$ satisfies the following formula for any virtual biset $X_G^H\in \AG(G,H)$:
\[(X_G^H)^\wedge_p = X_{\cF_G}^{\cF_H} \cmp (H_{\cF_H}^{\cF_H})^{-1} \in \AF(\cF_G,\cF_H).\]
\end{prop}

The free loop space $\freeO n G$ has components corresponding to all $G$-conjugacy classes of commuting $n$-tuples in $G$, while $\freeO n \cF_G$ only has components corresponding to $G$-conjugacy classes of commuting $n$-tuples in $S$, i.e. tuples of $p$-power order elements. Hence we cannot hope for some sort of equivalence between $\freeO n G$ and $\freeO n \cF_G$ on the nose.

However, recall from \cite[Proposition 3.35]{RSS_Bold1} that we define $\freeO n_p G$ to consist of the components in $\freeO n G$ corresponding to commuting tuples of $p$-power order elements up to $G$-conjugation. In this case, we do have a correspondence between the component of $\freeO n_p G$ and the components of $\freeO n \cF_G$.
\begin{convention}\label{conventionCommonReprs}
We use the same chosen representative $n$-tuples $\tup a$ in $\cF_G$ to represent both the components $C_{\cF_G}(\tup a)$ in $\freeO n\cF_G$ and the components $C_G(\tup a)$ in $\freeO n_p G$.
\end{convention}

Each component $C_{\cF_G}(\tup a)$ in $\freeO n \cF_G$ is the fusion system induced by $C_G(\tup a)$ on the Sylow $p$-subgroup $C_S(\tup a)$, and $C_G(\tup a)$ is the corresponding component of $\freeO n_p(G)$. Consequently the $p$-completion functor $(-)^\wedge_p\colon \AG\to \AF$ takes the group $C_G(\tup a)$ to $C_{\cF_G}(\tup a)$, and thus takes all of $\freeO n_p(G)$ to  $\freeO n \cF_G$.

On morphisms, according to Proposition \ref{propPCompletion}, $(-)^\wedge_p\colon \AG((\Z/p^e)^n \x\freeO n_p G, (\Z/p^e)^n\x \freeO n_p H) \to \AF((\Z/p^e)^n \x\freeO n \cF_G, (\Z/p^e)^n\x \freeO n \cF_H)$ is given on a biset matrix $X\in \AG((\Z/p^e)^n \x\freeO n_p G, (\Z/p^e)^n\x \freeO n_p H)$ by
\[
((X)^\wedge_p)_{\tup a,\tup b} = (X_{\tup a,\tup b})_{(\Z/p^e)^n\x C_{\cF_G}(\tup a)}^{(\Z/p^e)^n\x C_{\cF_H}(\tup b)} \cmp \Bigl(((\Z/p^e)^n\x C_H(\tup b))_{(\Z/p^e)^n\x C_{\cF_H}(\tup b)}^{(\Z/p^e)^n\x C_{\cF_H}(\tup b)}\Bigr)^{-1}
\]
for representative commuting $n$-tuples $\tup a$ in $\cF_G$ and $\tup b$ in $\cF_H$.

Now suppose we start with a virtual $(G,H)$-biset $X\in \AG(G,H)$. If we first apply $\twistO {n,p}$ and then apply $(-)^\wedge_p$, we get a matrix in $\AF((\Z/p^e)^n \x\freeO n \cF_G, (\Z/p^e)^n\x \freeO n \cF_H)$ consisting of virtual bisets
\begin{equation}\label{eqTwistThenComplete}
((\twistO {n,p} (X))^\wedge_p)_{\tup a,\tup b} = \Bigl(\twistO {n,p}(X)_{\tup a,\tup b}\Bigr)_{(\Z/p^e)^n\x C_{\cF_G}(\tup a)}^{(\Z/p^e)^n\x C_{\cF_H}(\tup b)} \cmp \Bigl(((\Z/p^e)^n\x C_H(\tup b))_{(\Z/p^e)^n\x C_{\cF_H}(\tup b)}^{(\Z/p^e)^n\x C_{\cF_H}(\tup b)}\Bigr)^{-1}.
\end{equation}
If, on the other hand, we first apply $(-)^\wedge_p$ and then apply $\twistO n\colon \AF\to \AF$, we can use the fact that $\twistO n$ is a functor to get a matrix in $\AF((\Z/p^e)^n \x\freeO n \cF_G, (\Z/p^e)^n\x \freeO n \cF_H)$ consisting of virtual bisets
\begin{equation}\label{eqCompleteThenTwist}
\begin{split}
\twistO n ((X)^\wedge_p)_{\tup a,\tup b} &= \Bigl(\twistO n(X_{\cF_G}^{\cF_H} \cmp (H_{\cF_H}^{\cF_H})^{-1})\Bigr)_{\tup a,\tup b} 
\\ &= \Bigl(\twistO n(X_{\cF_G}^{\cF_H}) \cmp \twistO n (H_{\cF_H}^{\cF_H})^{-1}\Bigr)_{\tup a,\tup b}.
\end{split}
\end{equation}
We claim that these two matrices are always equal (Theorem \ref{thmTwistedLoopsAndPCompletion} below), and in order to prove this we will make use of the following lemma:
\begin{lemma}\label{lemmaRestrictBifreeToSylow}
Suppose $H$ is a finite group with Sylow $p$-subgroup $T$ and let $S$ be a $p$-group. 
The restriction map $\AG(H,S)^\wedge_p \to \AG(T,S)^\wedge_p$ is injective on the subgroup of bifree virtual $(H,S)$-bisets.
\end{lemma}

\begin{proof}
Suppose we have a bifree virtual biset $Y\in \AG(H,S)^{\wedge}_p$, then $Y$ must be a $\Z_p$-linear combination of basis elements of the form $[P,\ph]_H^S$, where $P\leq H$ is a $p$-subgroup and $\ph\colon P\to S$ is injective. Such a linear combination is uniquely determined by its number of fixed points $\abs{Y^{(Q,\psi)}}\in \Z_p$, for all pairs $(Q,\psi)$ of a subgroup $Q\leq H$ and an injective group homomorphism $\psi\colon Q\to S$. Thus the virtual biset $Y$ is uniquely determined by $\abs{Y^{(Q,\psi)}}\in \Z_p$, where we only consider pairs $(Q,\psi)$ with $Q\leq H$ a $p$-subgroup.

Furthermore, since every $p$-subgroup of $H$ is $H$-conjugate to a subgroup of $T$, $Y$ is uniquely determined by the number of fixed points $\abs{Y^{(Q,\psi)}}\in \Z_p$ for $Q\leq T$ and $\psi\colon Q\to S$ injective. These fixed points happen to be preserved under restriction of the left action from $H$ to $T$, so it follows that any bifree $Y\in \AG(H,S)^{\wedge}_p$ is uniquely determined by its restriction $[T,\id]_T^H\cmp Y\in \AG(T,S)^{\wedge}_p$.
\end{proof}

We now prove Theorem \ref{thmIntroMainTheoremPCompletion} from the introduction.
\begin{theorem}\label{thmTwistedLoopsAndPCompletion}
We have $(-)^\wedge_p \circ \twistO {n,p} = \twistO n \circ (-)^\wedge_p$
as functors $\AG\to \AF$.
\end{theorem}

\begin{proof}
We aim to prove that the matrices described in \eqref{eqTwistThenComplete} and \eqref{eqCompleteThenTwist} are equal.

Let us for a moment consider the matrix $\twistO n (X_{\cF_G}^{\cF_H}).$ Since restricting the action of $G$ to $\cF_G$ is given by precomposition with $\omega_{\cF_G}\circ [S,\id]_S^G$, functoriality of $\twistO n$ gives us
\[
\twistO n (X_{\cF_G}^{\cF_H}) = \twistO n((\omega_{\cF_G})_{\cF_G}^S)\cmp \twistO n (X_S^T) \cmp \twistO n((\omega_{\cF_H})_T^{\cF_H}).
\]
Since every element in a $p$-group has $p$-power order, the functor $\twistO{n,p}$ coincides with $\twistO n$ when restricted to $p$-groups. Using functoriality of $\twistO{n,p}$, we can further decompose $\twistO n (X_{\cF_G}^{\cF_H})$ as
\begin{equation*}
\twistO n (X_{\cF_G}^{\cF_H}) = \twistO n((\omega_{\cF_G})_{\cF_G}^S) \cmp  \twistO{n,p}([S,\id]_S^G) \cmp \twistO {n,p} (X_G^H) \cmp \twistO{n,p}([T,\id]_H^T) \cmp \twistO n((\omega_{\cF_H})_T^{\cF_H}).
\end{equation*}
By Proposition \ref{propTwistedFusionLoop}, the entries of $\twistO n((\omega_{\cF_G})_{\cF_G}^S)$ satisfy
\[
\twistO n((\omega_{\cF_G})_{\cF_G}^S)_{\tup a,\tup s} = \twistO n((\omega_{\cF_G})_S^S)_{\tup a,\tup s},
\]
whenever $\tup a$ is a representative $n$-tuple in $\cF_G$, and $\tup s$ in $S$. For an $n$-tuple of $p$-power order elements $\tup g$ in $G$ and a representative $n$-tuple $\tup a$ in $\cF_G$, the entries of the composite $\twistO n((\omega_{\cF_G})_{\cF_G}^S) \cmp  \twistO{n,p}([S,\id]_S^G)$ satisfy
\begin{align*}
&(\twistO n((\omega_{\cF_G})_{\cF_G}^S) \cmp  \twistO{n,p}([S,\id]_S^G))_{\tup a, \tup g} 
\\ ={}& (\twistO n((\omega_{\cF_G})_S^S) \cmp  \twistO{n,p}([S,\id]_S^G))_{\tup a, \tup g} 
\\ ={}& (\twistO {n,p}((\omega_{\cF_G})_S^S) \cmp  \twistO{n,p}([S,\id]_S^G))_{\tup a, \tup g} 
\\ ={}& (\twistO {n,p}((\omega_{\cF_G})_S^S\cmp [S,\id]_S^G))_{\tup a, \tup g}
\\ ={}& (\twistO {n,p}([S,\id]_S^G))_{\tup a, \tup g} \hspace{2cm} \text{since $[S,\id]_S^G$ is left $\cF_G$-stable}.
\end{align*}
The biset $[S,\id]_S^G$ represents the forward map $\incl_S^G\colon S\to G$, hence by \cite[Theorems 3.33.(iv)]{RSS_Bold1} the entries $(\twistO {n,p}([S,\id]_S^G))_{\tup a, \tup g}$ are given by
\[
(\twistO {n,p}([S,\id]_S^G))_{\tup a, \tup g} = [(\Z/p^e)^n\x C_S(\tup a), \id]_{(\Z/p^e)^n\x C_S(\tup a)}^{(\Z/p^e)^n\x C_G(\tup g)},
\]
if $\tup g=\tup a$ (following Convention \ref{conventionCommonReprs}), and $0$ otherwise.

To sum up, the composite $\twistO n((\omega_{\cF_G})_{\cF_G}^S) \cmp  \twistO{n,p}([S,\id]_S^G)$ is a diagonal matrix with entries
\begin{equation}\label{eqPCompletionPreComposition}
(\twistO n((\omega_{\cF_G})_{\cF_G}^S) \cmp  \twistO{n,p}([S,\id]_S^G))_{\tup a, \tup a} = [(\Z/p^e)^n\x C_S(\tup a), \id]_{(\Z/p^e)^n\x C_S(\tup a)}^{(\Z/p^e)^n\x C_G(\tup a)},
\end{equation}
which are left $((\Z/p^e)^n\x C_{\cF_G}(\tup a))$-stable. This will be the first main ingredient in the proof of the theorem.

\medskip\noindent For the second ingredient, consider the matrix entry 
\begin{equation*} 
    \Bigl(\twistO{n,p}([T,\id]_H^T) \cmp \twistO n((\omega_{\cF_H})_T^{\cF_H})\cmp \twistO n(H_{\cF_H}^{\cF_H})^{-1}\Bigr)_{\tup h,\tup b}
\end{equation*}
inside the right $((\Z/p^e)^n\x C_{\cF_H}(\tup b))$-stable elements of $\AG((\Z/p^e)^n\x C_H(\tup h), (\Z/p^e)^n\x C_T(\tup b))^\wedge_p$.
We claim that this entry is $0$ unless $\tup h=\tup b$, in which case it equals the composite 
\[[(\Z/p^e)^n\x C_T(\tup b),\id]_{(\Z/p^e)^n\x C_H(\tup b)}^{(\Z/p^e)^n\x C_T(\tup b)}\cmp \Bigl(((\Z/p^e)^n\x C_H(\tup b))_{(\Z/p^e)^n\x C_{\cF_H}(\tup b)}^{(\Z/p^e)^n\x C_{\cF_H}(\tup b)}\Bigr)^{-1}.\] 
Each of the virtual bisets $\twistO{n,p}([T,\id]_H^T)$, $\twistO n((\omega_{\cF_H})_T^{\cF_H})$ and $\twistO n(H_{\cF_H}^{\cF_H})^{-1}$ are composites of bifree bisets and so are bifree as well -- note that $\twistO n$ preserves bifree actions, since the $\wind$-maps involved in the construction are all injective group homomorphisms.
We shall apply Lemma \ref{lemmaRestrictBifreeToSylow} to the bifree matrix entry 
\[\Bigl(\twistO{n,p}([T,\id]_H^T) \cmp \twistO n((\omega_{\cF_H})_T^{\cF_H})\cmp \twistO n(H_{\cF_H}^{\cF_H})^{-1}\Bigr)_{\tup h,\tup b}\]
in $\AG((\Z/p^e)^n\x C_H(\tup h), (\Z/p^e)^n\x C_T(\tup b))^\wedge_p$.

Recall \eqref{eqPCompletionPreComposition} from the first part of the proof, and replace $G$, $S$, and $\tup a$, with $H$, $T$, and $\tup h$, respectively. We then know that $\twistO n((\omega_{\cF_H})_{\cF_H}^T) \cmp  \twistO{n,p}([T,\id]_T^H)$ is a diagonal matrix with entries \[
(\twistO n((\omega_{\cF_H})_{\cF_H}^T) \cmp  \twistO{n,p}([T,\id]_T^H))_{\tup h, \tup h} = [(\Z/p^e)^n\x C_T(\tup h), \id]_{(\Z/p^e)^n\x C_T(\tup h)}^{(\Z/p^e)^n\x C_H(\tup h)}.
\]
From this it follows that restricting the virtual biset 
\[\Bigl(\twistO{n,p}([T,\id]_H^T) \cmp \twistO n((\omega_{\cF_H})_T^{\cF_H})\cmp \twistO n(H_{\cF_H}^{\cF_H})^{-1}\Bigr)_{\tup h,\tup b}\]
on the left from $(\Z/p^e)^n\x C_H(\tup h)$ to $(\Z/p^e)^n\x C_T(\tup h)$ can be achieved by precomposing the entire matrix with the diagonal matrix $\twistO n((\omega_{\cF_H})_{\cF_H}^T) \cmp  \twistO{n,p}([T,\id]_T^H)$:
\begin{align*}
& [(\Z/p^e)^n\x C_T(\tup h),\id]_{(\Z/p^e)^n\x C_T(\tup h)}^{(\Z/p^e)^n\x C_H(\tup h)}\cmp \Bigl(\twistO{n,p}([T,\id]_H^T) \cmp \twistO n((\omega_{\cF_H})_T^{\cF_H})\cmp \twistO n(H_{\cF_H}^{\cF_H})^{-1}\Bigr)_{\tup h,\tup b}
\\ ={}& \Bigl(\twistO n((\omega_{\cF_H})_{\cF_H}^T) \cmp  \twistO{n,p}([T,\id]_T^H)\cmp \twistO{n,p}([T,\id]_H^T) \cmp \twistO n((\omega_{\cF_H})_T^{\cF_H})\cmp \twistO n(H_{\cF_H}^{\cF_H})^{-1}\Bigr)_{\tup h,\tup b}
\\ ={}& \Bigl(\twistO n(H_{\cF_H}^{\cF_H})\cmp \twistO n(H_{\cF_H}^{\cF_H})^{-1}\Bigr)_{\tup h,\tup b}
\\ ={}& \Bigl(\twistO n((\omega_{\cF_H})_{\cF_H}^{\cF_H})\Bigr)_{\tup h,\tup b},
\end{align*}
which is the identity on $(\Z/p^e)^n\x C_{\cF_H}(\tup b)$ when $\tup h=\tup b$ and $0$ otherwise.

Because $\Bigl(\twistO{n,p}([T,\id]_H^T) \cmp \twistO n((\omega_{\cF_H})_T^{\cF_H})\cmp \twistO n(H_{\cF_H}^{\cF_H})^{-1}\Bigr)_{\tup h,\tup b}$ is uniquely determined by its restriction to $(\Z/p^e)^n\x C_T(\tup h)$, we conclude that it is zero unless $\tup h =\tup b$.

At the same time, we can also consider the restriction of 
\[[(\Z/p^e)^n\x C_T(\tup b),\id]_{(\Z/p^e)^n\x C_H(\tup b)}^{(\Z/p^e)^n\x C_T(\tup b)}\cmp \Bigl(((\Z/p^e)^n\x C_H(\tup b))_{(\Z/p^e)^n\x C_{\cF_H}(\tup b)}^{(\Z/p^e)^n\x C_{\cF_H}(\tup b)}\Bigr)^{-1}\]
from $(\Z/p^e)^n\x C_H(\tup b)$ to $(\Z/p^e)^n\x C_T(\tup b)$ on the left:
\begin{align*}
& [(\Z/p^e)^n\x C_T(\tup b),\id]_{(\Z/p^e)^n\x C_T(\tup b)}^{(\Z/p^e)^n\x C_H(\tup h)}\cmp [(\Z/p^e)^n\x C_T(\tup b),\id]_{(\Z/p^e)^n\x C_H(\tup b)}^{(\Z/p^e)^n\x C_T(\tup b)}
\\* &\qquad\cmp \Bigl(((\Z/p^e)^n\x C_H(\tup b))_{(\Z/p^e)^n\x C_{\cF_H}(\tup b)}^{(\Z/p^e)^n\x C_{\cF_H}(\tup b)}\Bigr)^{-1}
\\ ={}& ((\Z/p^e)^n\x C_H(\tup b))_{(\Z/p^e)^n\x C_{\cF_H}(\tup b)}^{(\Z/p^e)^n\x C_{\cF_H}(\tup b)}\cmp \Bigl(((\Z/p^e)^n\x C_H(\tup b))_{(\Z/p^e)^n\x C_{\cF_H}(\tup b)}^{(\Z/p^e)^n\x C_{\cF_H}(\tup b)}\Bigr)^{-1}
\\ ={}& \id_{(\Z/p^e)^n\x C_{\cF_H}(\tup b)}.
\end{align*}
We again get the identity on $(\Z/p^e)^n\x C_{\cF_H}(\tup b)$, and again since the bifree bisets are uniquely determined by their restrictions, we conclude that the matrix $\twistO{n,p}([T,\id]_H^T) \cmp \twistO n((\omega_{\cF_H})_T^{\cF_H})\cmp \twistO n(H_{\cF_H}^{\cF_H})^{-1}$ is diagonal with entries
\begin{multline}\label{eqPCompletionPostComposition}
 \Bigl(\twistO{n,p}([T,\id]_H^T) \cmp \twistO n((\omega_{\cF_H})_T^{\cF_H})\cmp \twistO n(H_{\cF_H}^{\cF_H})^{-1}\Bigr)_{\tup b,\tup b} 
 \\= [(\Z/p^e)^n\x C_T(\tup b),\id]_{(\Z/p^e)^n\x C_H(\tup b)}^{(\Z/p^e)^n\x C_T(\tup b)}\cmp \Bigl(((\Z/p^e)^n\x C_H(\tup b))_{(\Z/p^e)^n\x C_{\cF_H}(\tup b)}^{(\Z/p^e)^n\x C_{\cF_H}(\tup b)}\Bigr)^{-1}.
\end{multline}
This is the second main ingredient in the proof of the Theorem.

The claim of the theorem that \eqref{eqCompleteThenTwist} equals \eqref{eqTwistThenComplete} is now a straightforward check of matrix entries:
\begin{align*}
    & \twistO n ((X)^\wedge_p)_{\tup a,\tup b} 
\\ \overset{\eqref{eqCompleteThenTwist}}={}& \Bigl(\twistO n(X_{\cF_G}^{\cF_H}) \cmp \twistO n (H_{\cF_H}^{\cF_H})^{-1}\Bigr)_{\tup a,\tup b}
\\ ={}& \Bigl(\twistO n((\omega_{\cF_G})_{\cF_G}^S) \cmp  \twistO{n,p}([S,\id]_S^G) \cmp \twistO {n,p} (X_G^H) \cmp \twistO{n,p}([T,\id]_H^T) \cmp \twistO n((\omega_{\cF_H})_T^{\cF_H}) \cmp \twistO n (H_{\cF_H}^{\cF_H})^{-1}\Bigr)_{\tup a,\tup b}
\\ \overset{\eqref{eqPCompletionPreComposition}}={}& [(\Z/p^e)^n\x C_S(\tup a), \id]_{(\Z/p^e)^n\x C_S(\tup a)}^{(\Z/p^e)^n\x C_G(\tup a)} \cmp \Bigl(\twistO {n,p} (X_G^H) \cmp \twistO{n,p}([T,\id]_H^T) \cmp \twistO n((\omega_{\cF_H})_T^{\cF_H}) \cmp \twistO n (H_{\cF_H}^{\cF_H})^{-1}\Bigr)_{\tup a,\tup b}
\\ \overset{\eqref{eqPCompletionPostComposition}}={}& [(\Z/p^e)^n\x C_S(\tup a), \id]_{(\Z/p^e)^n\x C_S(\tup a)}^{(\Z/p^e)^n\x C_G(\tup a)} \cmp \Bigl(\twistO {n,p} (X_G^H) \Bigr)_{\tup a,\tup b} \cmp [(\Z/p^e)^n\x C_T(\tup b),\id]_{(\Z/p^e)^n\x C_H(\tup b)}^{(\Z/p^e)^n\x C_T(\tup b)}
\\* &\qquad\cmp \Bigl(((\Z/p^e)^n\x C_H(\tup b))_{(\Z/p^e)^n\x C_{\cF_H}(\tup b)}^{(\Z/p^e)^n\x C_{\cF_H}(\tup b)}\Bigr)^{-1}
\\ ={}& \Bigl(\twistO {n,p} (X_G^H)_{\tup a,\tup b} \Bigr)_{(\Z/p^e)^n\x C_{\cF_G}(\tup a)}^{(\Z/p^e)^n\x C_{\cF_H}(\tup b)}\cmp \Bigl(((\Z/p^e)^n\x C_H(\tup b))_{(\Z/p^e)^n\x C_{\cF_H}(\tup b)}^{(\Z/p^e)^n\x C_{\cF_H}(\tup b)}\Bigr)^{-1}
\\ \overset{\eqref{eqTwistThenComplete}}={}& ((\twistO {n,p} (X))^\wedge_p)_{\tup a,\tup b}.
\end{align*}
Finally we conclude that the composite functors $\twistO n \circ (-)^\wedge_p$ and $(-)^\wedge_p\circ \twistO{n,p}$ agree as functors $\AG\to \AF$.
\end{proof}

\appendix
\section{Auxiliary results about characteristic idempotents}
In this appendix we gather a few utility results about characteristic idempotents and ``virtual double cosets.'' The results are mostly elementary. However, to our knowledge, they do not appear in the literature.

First we observe a simple principle for recognizing fusion preserving maps in terms of characteristic idempotents:
\begin{lemma}\label{lemmaFusionPreserving}
Let $\cE$ and $\cF$ be saturated fusion systems over $p$-groups $R$ and $S$ respectively. Then the following statements are equivalent for any homomorphism $f\colon R\to S$.
\begin{enumerate}
\renewcommand{\theenumi}{$(\roman{enumi})$}\renewcommand{\labelenumi}{\theenumi}
    \item\label{itemFusionPreserving} $f$ is a fusion preserving map from $\cE$ to $\cF$.
    \item\label{itemAbsorbsLeftIdempotent} The equation
    \(
    \omega_\cE \cmp [R,f]_R^S\cmp \omega_\cF = [R,f]_R^S\cmp \omega_\cF
    \)
    holds in $\AF(R,S)$.
    \item\label{itemMapIsLeftStable} The virtual biset $[R,f]_R^S\cmp \omega_\cF\in \AF(R,S)$ is left $\cE$-stable.
\end{enumerate}
\end{lemma}

\begin{proof}
The properties \ref{itemAbsorbsLeftIdempotent} and \ref{itemMapIsLeftStable} are immediately seen to be equivalent since an element $X\in \AF(R,S)$ is left $\cE$-stable if and only if $\omega_\cE\cmp X=X$.

We will show that \ref{itemFusionPreserving} implies \ref{itemMapIsLeftStable}. Suppose that $f$ is fusion preserving. We will prove that $[R,f]_R^S\cmp \omega_\cF$ is left $\cE$-stable. Suppose that $P\leq R$ and $\ph\in \cE(P,R)$ is any morphism in the fusion system $\cE$. Since $f$ is fusion preserving, we have a corresponding map $\tilde\ph\in \cF(f(P),S)$ such that $f\circ \ph = \tilde\ph\circ f$ as maps $P\to S$. If we now restrict $[R,f]_R^S\cmp \omega_\cF$ along $\ph$ on the left, we get
\begin{multline*}
    [P,\ph]_P^R\cmp [R,f]_R^S\cmp \omega_\cF = [P,f\circ \ph]_R^S \cmp \omega_\cF 
    \\= [P,\tilde\ph\circ f]_R^S \cmp \omega_\cF = [P,f]_P^{f(P)}\cmp [f(P),\tilde\ph]_{f(P)}^S \cmp \omega_\cF.
\end{multline*}
Because $\omega_\cF$ is left $\cF$-stable, restricting $\omega_\cF$ along $\tilde\ph$ is the same as restricting along the inclusion:
\begin{multline*}
    [P,f]_P^{f(P)}\cmp [f(P),\tilde\ph]_{f(P)}^S \cmp \omega_\cF = [P,f]_P^{f(P)}\cmp [f(P),\incl]_{f(P)}^S \cmp \omega_\cF 
    \\= [P,\incl]_P^{R}\cmp [R,f]_R^S \cmp \omega_\cF.
\end{multline*}
Consequently, $[R,f]_R^S \cmp \omega_\cF$ is left $\cE$-stable as required.

Conversely, suppose that $f$ satisfies the equation \ref{itemAbsorbsLeftIdempotent}. We will then prove that $f$ is fusion preserving. To that purpose, suppose $P\leq R$ and $\ph\in\cE(P,R)$; we shall prove that there is a corresponding map $f(P)\to S$ inside $\cF$. Consider the basis element of $\AF(P,\cF)$ given by the composite $f\circ \ph$:
\[
[P,f\circ \ph]_P^\cF = [P,\ph]_P^R\cmp [R,f]_R^S\cmp \omega_\cF.
\]
By \ref{itemAbsorbsLeftIdempotent} we can add in $\omega_\cE$ in the middle and we get
\begin{multline*}
    [P,\ph]_P^R\cmp [R,f]_R^S\cmp \omega_\cF = [P,\ph]_P^R\cmp \omega_\cE\cmp  [R,f]_R^S\cmp \omega_\cF = [P,\incl]_P^R\cmp \omega_\cE\cmp  [R,f]_R^S\cmp \omega_\cF
    \\ = [P,\incl]_P^R\cmp [R,f]_R^S\cmp \omega_\cF = [P,f]_P^\cF.
\end{multline*}
These calculations show that the two basis elements $[P,f\circ \ph]_P^\cF$ and $[P,f]_P^\cF$ in $\AF(P,\cF)$ are equal. Consequently, $f\circ\ph$ arises from $f$ by simultaneously precomposing with some conjugation $c_x$ in $P$ and postcomposing with some map $\psi$ in $\cF$, i.e. we have:
\[f\circ \ph= \psi\circ f\circ c_x\]
as maps $P\to S$ with $x\in N_R(P)$ and $\psi\in\cF(f(P),S)$. We can now apply $f$ to the element $x$ as well, resulting in
\[f\circ \ph= \psi\circ f\circ c_x = \psi\circ c_{f(x)}\circ f.\]
Hence $\tilde\ph:=\psi\circ c_{f(x)}$ serves the purpose, proving that $f$ is in fact fusion preserving from $\cE$ to $\cF$.
\end{proof}

Next we will work towards writing out a double coset formula for composing virtual bisets of saturated fusion systems.
\begin{lemma}\label{lemmaBisetDoubleCosetFormula}
Let $R$, $S$, and $T$ be finite $p$-groups, and consider transitive bisets $[H,\ph]_R^S$ and $[K,\psi]_S^T$. Suppose $X$ is a bifree $(S,S)$-biset. For each $x\in X$ the stabilizer of $x$ is given as a pair of a subgroup $S_x\leq S$ and a homomorphism $c_x\colon S_x\to S$. That is, for all $s\in S_x$, we have $s x = x c_x(s)$.
The composition $[H,\ph]_R^S\cmp X\cmp [K,\psi]_S^T$ can then be calculated in terms of double cosets of $X$:
\[[H,\ph]_R^S\cmp X\cmp [K,\psi]_S^T = \sum_{x\in \ph H\backslash X / K} [\ph^{-1}(S_x\cap c_x^{-1}(K)), \psi\circ c_x\circ \ph]_R^T.\]
\end{lemma}

\begin{proof}
First write the composition $[H,\ph]_R^S\cmp X\cmp [K,\psi]_S^T$ as
\[[H,\ph]_R^{\ph H} \cmp [\ph H,\id]_{\ph H}^S\cmp X\cmp [K,\id]_S^K \cmp [K,\psi]_K^T.\]
The middle part, $[\ph H,\id]_{\ph H}^S\cmp X\cmp [K,\id]_S^K$, is then simply the $(S,S)$-biset $X$ restricted to the subgroup $\ph H$ on the left and $K$ on the right.

Now the orbits of the restriction ${}_{\ph H} X_K= [\ph H,\id]_{\ph H}^S\cmp X\cmp [K,\id]_S^K$ are precisely the double cosets $\ph H \backslash X / K$. Given a representative $x\in \ph H \backslash X / K$ for any of the orbits, the stabilizer of $x\in {}_{\ph H} X_K$ is obtained by restricting the stabilizer of $x$ inside the $(S,S)$-biset with the smaller subgroups $\ph H$ and $K$. Since the stabilizer of $x$ inside ${}_S X_S$ is given by $S_x\leq S$ and $c_x\colon S_x\to S$, the restriction to $(\ph H, K)$ becomes the subgroup $\ph H\cap S_x \cap c_x^{-1}(K)\leq \ph H$ and we may restrict $c_x$ to the homomorphism $c_x\colon \ph H\cap S_x \cap c_x^{-1}(K) \to K$.

When we take representatives for all the double cosets $\ph H \backslash X / K$, we get the following expression for the restriction ${}_{\ph H} X_K$:
\[{}_{\ph H} X_K = \sum_{x\in \ph H\backslash X / K} [\ph H \cap S_x\cap c_x^{-1}(K), c_x]_{\ph H}^K.\]
It remains to compose ${}_{\ph H} X_K$ with $[H,\ph]_R^{\ph H}$ on the left and $[K,\psi]_K^T$ on the right. This amounts to composing with the isomorphisms $\ph\colon H\to \ph H$ and $\psi\colon K\to \psi K$ and then inducing the actions from $H$ and $\psi K$ up to $R$ and $T$ respectively. Since induction doesn't change point stabilizers, we get the final formula:
\begin{align*}
&{} [H,\ph]_R^S\cmp X\cmp [K,\psi]_S^T
\\ ={}& [H,\id]_R^H \cmp [H,\ph]_H^{\ph H} \cmp \Bigl(\sum_{x\in \ph H\backslash X / K} [\ph H \cap S_x\cap c_x^{-1}(K), c_x]_{\ph H}^K\Bigr) \cmp [K,\psi]_K^{\psi K} \cmp [\psi K, \id]_{\psi K}^T
\\ ={}& [H,\id]_R^H \cmp \Bigl(\sum_{x\in \ph H\backslash X / K} [\ph^{-1}(S_x\cap c_x^{-1}(K)), \psi\circ c_x\circ \ph]_H^{\psi K}\Bigr) \cmp [\psi K, \id]_{\psi K}^T
\\ ={}& \sum_{x\in \ph H\backslash X / K} [\ph^{-1}(S_x\cap c_x^{-1}(K)), \psi\circ c_x\circ \ph]_R^T. \qedhere
\end{align*}
\end{proof}

\begin{remark}\label{remarkVirtualDoubleCosetFormula}
Lemma \ref{lemmaBisetDoubleCosetFormula} is linear in $X$ with respect to addition/disjoint union of bifree $(S,S)$-biset. We can therefore extend the lemma to cover all virtual bifree $(S,S)$-bisets $X\in \AG(S,S)^\wedge_p$, as long as we extend the notation
\[\sum_{x\in \ph H\backslash X / K} \dotsb\]
linearly to virtual bisets as well. The convention will be that each virtual point $x\in X$ is assigned a weight $\e_x$ equal to the coefficient of the virtual orbit in $X$ that contains $x$. Explicitly this means that if we decompose $X$ in terms of orbits
\[X = \sum_{\substack{L\leq S\\ \rho\colon L\to S}} \e_{L,\rho}\cdot [L,\rho]_S^S,\]
then each point of $[L,\rho]_S^S$ is assigned the weight $\e_{L,\rho}$, and the formula of Lemma \ref{lemmaBisetDoubleCosetFormula} becomes
\begin{align*}
&{}[H,\ph]_R^S\cmp X\cmp [K,\psi]_S^T
\\={}& \sum_{x\in \ph H\backslash X / K} [\ph^{-1}(S_x\cap c_x^{-1}(K)), \psi\circ c_x\circ \ph]_R^T
\\ \overset{\text{def}}={}& \sum_{\substack{L\leq S\\ \rho\colon L\to S}} \sum_{x\in \ph H\backslash [L,\rho]_S^S / K} \e_{L,\rho}\cdot[\ph^{-1}(S_x\cap c_x^{-1}(K)), \psi\circ c_x\circ \ph]_R^T.
\end{align*}
This formula is then valid for every bifree virtual biset $X\in \AG(S,S)^\wedge_p$.
\end{remark}

\begin{prop}\label{propFusionDoubleCosetFormula}
Suppose that $\cE$, $\cF$, and $\cG$ are saturated fusion systems on $p$-groups $R$, $S$, and $T$, respectively, and consider virtual bisets $[H,\ph]_\cE^\cF$ and $[K,\psi]_\cF^\cG$, where $H\leq R$, $\ph\colon H\to S$, $K\leq S$, and $\psi\colon K\to T$.
The composition of these bisets satisfies a double coset formula following the conventions of Lemma \ref{lemmaBisetDoubleCosetFormula} and Remark \ref{remarkVirtualDoubleCosetFormula}:
\[[H,\ph]_\cE^\cF\cmp [K,\psi]_\cF^\cG = \sum_{x\in \ph H\backslash \omega_\cF / K} [\ph^{-1}(S_x\cap c_x^{-1}(K)), \psi\circ c_x\circ \ph]_\cE^\cG.\]
In the special case where $K=S$ and $\psi\colon S\to T$ is a fusion preserving map from $\cF$ to $\cG$, the composition simplifies to
\[[H,\ph]_\cE^\cF \cmp [S,\psi]_\cF^\cG = [H,\psi\circ \ph]_\cE^\cG.\]
\end{prop}

\begin{proof}
By definition of the basis elements $[H,\ph]_\cE^\cF$ and $[K,\psi]_\cF^\cG$ (see the discussion following Definition \ref{defFusionBisetModule}) the expression $[H,\ph]_\cE^\cF\cmp [K,\psi]_\cF^\cG$ is shorthand for
\[[H,\ph]_\cE^\cF\cmp [K,\psi]_\cF^\cG = \omega_{\cE} \cmp [H,\ph]_R^S \cmp \omega_{\cF} \cmp [K,\psi]_S^T \cmp \omega_{\cG}.\]
We apply Lemma \ref{lemmaBisetDoubleCosetFormula} to the middle triple $[H,\ph]_R^S \cmp \omega_{\cF} \cmp [K,\psi]_S^T$ and get
\begin{align*}
{}& [H,\ph]_\cE^\cF\cmp [K,\psi]_\cF^\cG
\\ ={}& \omega_{\cE} \cmp \Bigl(\sum_{x\in \ph H\backslash \omega_\cF / K} [\ph^{-1}(S_x\cap c_x^{-1}(K)), \psi\circ c_x\circ \ph]_R^T \Bigr) \cmp \omega_{\cG}
\\ ={}& \sum_{x\in \ph H\backslash \omega_\cF / K} [\ph^{-1}(S_x\cap c_x^{-1}(K)), \psi\circ c_x\circ \ph]_\cE^\cG.
\end{align*}
This covers the general case of the proposition.

In the special case where $K=S$ and $\psi\colon S\to T$ is fusion preserving from $\cF$ to $\cG$, we refer to \cite{RSS_p-completion}*{Lemma 4.6} which states that $\omega_{\cF} \cmp [S,\psi]_S^T \cmp \omega_{\cG} = [S,\psi]_S^T \cmp \omega_{\cG}$ when $\psi$ is fusion preserving. This in turn gives us the required formula:
\begin{align*}
{}& [H,\ph]_\cE^\cF\cmp [S,\psi]_\cF^\cG
\\ ={}& \omega_{\cE} \cmp [H,\ph]_R^S \cmp \omega_{\cF} \cmp [S,\psi]_S^T \cmp \omega_{\cG}
\\ ={}& \omega_{\cE} \cmp [H,\ph]_R^S \cmp [S,\psi]_S^T \cmp \omega_{\cG}
\\ ={}& \omega_{\cE} \cmp [H,\psi \circ \ph]_R^T \cmp \omega_{\cG}
\\ ={}& [H,\psi \circ \ph]_\cE^{\cG}.\qedhere
\end{align*}
\end{proof}

\begin{prop}\label{propCentralizerIdempotent}
Let $\cF$ be a saturated fusion system on $S$ and suppose $a\in S$ is central in $\cF$, meaning that $C_\cF(a) = \cF$. Then $\omega_\cF$ only contains orbits of the form $[P,\ph]$, where $a\in P$ (and $\ph(a)=a$ since $a$ is $\cF$-central).
\end{prop}

\begin{proof}
Suppose $a\in S$ is central in the saturated fusion system $\cF$, and let $\omega_\cF$ be the characteristic idempotent.
In particular this means that if $a\leq P\leq S$ and $\ph\in \cF(P,S)$, then $\ph(a)=a$.

Furthermore, we let $\bar \omega_\cF$ be the $a$-centralizing part of $\omega_\cF$ meaning that
\[
\bar \omega_\cF := \sum_{\substack{a\in P\leq S\\ \ph\in \cF(P,S)\\ \text{up to $(\cF,\cF)$-conj.}}} c_{P,\ph} \cdot [P,\ph]_S^S,
\]
where $c_{P,\ph}\in \Z_p$ is the coefficient of $[P,\ph]_S^S$ in the orbit decomposition of $\omega_\cF$. We wish to prove that $\omega_\cF$ does not contain any additional orbits with $a\not\in P$, so we can equivalently prove that $\bar \omega_\cF = \omega_\cF$.

The algorithm for constructing $\omega_\cF$ in \cite{ReehIdempotent} starts with the transitive biset $[S,\id]_S^S$ and then proceeds one $(\cF,\cF)$-conjugacy class of pairs $(P,\ph)$ at a time (in decreasing order) and adds/subtracts orbits in the conjugacy class of $(P,\ph)$ to make the biset stable at that conjugacy class and above.

The virtual biset $\bar \omega_\cF$ is then the intermediate result of this algorithm where we have stabilized $[S,\id]_S^S$ only across all pairs $(P,\ph)$ with $a\leq P$. If we can prove that $\bar\omega_\cF$ is in fact fully $(\cF,\cF)$-stable instead of just stable for pairs containing $a$, this would mean that no further stabilization is needed and the algorithm stops with $\omega_\cF=\bar\omega_\cF$.

For readers who don't feel comfortable with the argument above, we will bring a more explicit way of finishing the proof below -- once we have proved that $\bar\omega_\cF$ is $(\cF,\cF)$-stable.

\textbf{$\bar\omega_\cF$ is $(\cF,\cF)$-stable:} This is essentially \cite{GelvinReeh}*{Proposition 9.10}, using that $\cF$ is the centralizer fusion system of $a$. The difference is that \cite{GelvinReeh} deals with actual bisets, while $\bar\omega_\cF$ here is a virtual biset. We adapt the proof used in \cite{GelvinReeh}.

One of the equivalent ways to state $(\cF,\cF)$-stability for an $(S,S)$-biset $X$ is the property that the number of fixed points for a pair $(Q,\psi)$ only depends on $(Q,\psi)$ up to $(\cF,\cF)$-conjugation:
\[\abs{X^{(Q,\psi)}} = \abs{X^{(Q',\psi')}}\]
whenever $(Q',\psi')$ is $(\cF,\cF)$-conjugate to $(Q,\psi)$. The fixed point set $X^{(Q,\psi)}$ is the set of $x\in X$ such that $q x = x\psi(q)$ for all $q\in Q$. Working with virtual fixed points, the same characterization of $(\cF,\cF)$-stability works for virtual bisets $X\in \AG(S,S)^\wedge_p$.

If $X\in \AG(S,S)^\wedge_p$ is $\cF$-generated, such as $\omega_\cF$ or $\bar\omega_\cF$, then it is sufficient to check fixed points for $(Q,\psi)$ with $\psi\in \cF(Q,S)$, where we ask that
\[\abs{X^{(Q,\psi)}} = \abs{X^{(Q',\id)}}\]
for any $\cF$-conjugate $Q'$ to $Q$.

If $(Q,\psi)$ has $a\in Q$, then the fixed points $\abs{\omega_\cF^{(Q,\psi)}}$ only depends on the orbits with stabilizers containing $Q$ and hence $a$. We thus have
\[
\abs{\bar\omega_\cF^{(Q,\psi)}} = \abs{\omega_\cF^{(Q,\psi)}} =  \abs{\omega_\cF^{(Q',\id)}} = \abs{\bar\omega_\cF^{(Q',\id)}},
\]
whenever $Q'$ is $\cF$-conjugate to $Q$, $Q$ contains $a$, and $\psi\in \cF(Q,S)$.

Given any pair $(Q,\psi)$ with $\psi\in \cF(Q,S)$ and $a\not\in Q$, the assumption that $a$ is central in $\cF$ means that $\psi$ extends (uniquely) to a homomorphism $\hat\psi\colon Q\gen a \to S$ such that $\hat\psi|_Q=\psi$ and $\hat\psi(a)=a$.

Because $\bar\omega_\cF$ only has orbits with stabilizers that contain $a$, and since the extension $\hat\psi$ is unique, we get
\[\abs{\bar\omega_\cF^{(Q,\psi)}} = \abs{\bar\omega_\cF^{(Q\gen a,\hat\psi)}}.\]
We already know that the fixed points are $(\cF,\cF)$-conjugation invariant for $(Q\gen a,\hat\psi)$, so the same is true for $(Q,\psi)$. Hence $\bar\omega_\cF$ is in fact $(\cF,\cF)$-stable.

This completes the proof that the algorithm for constructing $\omega_\cF$ stops once we have $\bar\omega_\cF$. Alternatively, though it boils down to the same formulas, we can complete the proof of the lemma as follows.

\textbf{Alternative proof that $\bar\omega_\cF = \omega_\cF$:}
Let $c_{P,\ph}$ be the coefficient of $[P,\ph]_S^S$ in the orbit decomposition of $\bar\omega_{\cF}$. Hence $c_{P,\ph}=0$ unless $a\in P\leq S$ and $\ph\in \cF(P,S)$.

Because $\omega_\cF$ is obtained by stabilizing the transitive biset $[S,\id]_S^S$. Remark 4.7 of \cite{ReehIdempotent} states that the coefficients satisfy
\[
\sum_{\substack{\text{$(P',\ph')$ up to $(S,S)$-conj.}\\ \text{s.t. $(P',\ph')$ is $(\cF,\cF)$-conj. to $(P,\id)$}}}\hspace{-1cm} c_{P',\ph'} = \begin{cases} 1 &\text{if $P=S$} \\ 0 &\text{otherwise},\end{cases}
\]
when $a$ is in $P$ and $P'$. When $a$ is not in the subgroups, all the coefficients are zero so the formula still holds.

If we multiply $\bar\omega_\cF$ with $\omega_\cF$ from both sides, the formula above gives us
\begin{align*}
{}& \omega_\cF\cmp \bar\omega_\cF \cmp \omega_\cF
\\ ={}& \sum_{\text{$(P,\ph)$ up to $(S,S)$-conj.}} \hspace{-1cm} c_{P,\ph}\cdot (\omega_\cF \cmp [P,\ph]_S^S \cmp \omega_\cF)
\\ ={}& \sum_{\text{$P\leq S$ up to $\cF$-conj.}} \hspace{-.8cm} [P,\id]_\cF^\cF \cdot \Bigl( \sum_{\substack{\ph'\in \cF(P',S)\\\text{up to $(S,S)$-conj.}\\\text{s.t. $P'$ is $\cF$-conj to $P$}}} c_{P',\ph'}\Bigr)
\\ ={}& [S,\id]_\cF^\cF \cdot 1
\\ ={}& \omega_\cF.
\end{align*}
At the same time, $\bar\omega_\cF$ is $(\cF,\cF)$-stable as proven earlier, so
\[\omega_\cF\cmp \bar\omega_\cF\cmp \omega_\cF=\bar\omega_\cF.\]
This finishes the alternative proof that $\omega_\cF=\bar\omega_\cF$.
\end{proof}

\makeatletter
\def\eprint#1{\@eprint#1 }
\def\@eprint #1:#2 {%
    \ifthenelse{\equal{#1}{arXiv}}%
        {\href{http://front.math.ucdavis.edu/#2}{arXiv:#2}}%
        {\href{#1:#2}{#1:#2}}%
}
\makeatother

\end{document}